\documentclass[11pt,a4paper,reqno]{amsart}
\usepackage{amsfonts}
\usepackage{amsmath,amssymb,amsthm,amsxtra,mathabx}
\usepackage{float}
\usepackage[colorlinks, linkcolor=blue,anchorcolor=Periwinkle,
    citecolor=red,urlcolor=Emerald]{hyperref}
\usepackage[usenames,dvipsnames]{xcolor}
\usepackage{geometry,array} 
\usepackage{graphicx,subfigure,bookmark,tikz,pgfplots,url,transparent,pdflscape,transparent}
\makeatletter
\newcommand{\addresseshere}{%
  \enddoc@text\let\enddoc@text\relax
}
\makeatother
\title{Geometric classification of total stability conditions}
\author{Yu Qiu}
\address{Qy:
	Yau Mathematical Sciences Center and Department of Mathematical Sciences,
	Tsinghua University,
    100084 Beijing,
    China.
    \&
    Beijing Institute of Mathematical Sciences and Applications, Yanqi Lake, Beijing, China}
\email{yu.qiu@bath.edu}

\author{Xiaoting Zhang}
\address{Zx:
    Beijing Advanced Innovation Center for Imaging Theory and Technology, Academy for Multidisciplinary Studies, Capital Normal University, Beijing 100048, China}
\email{xiaoting.zhang09@hotmail.com}
\dedicatory{Dedicated to Alastair King on the occasion of his sixtieth birthday}
\date{\today}

\usetikzlibrary{matrix,positioning,decorations.markings,arrows,decorations.pathmorphing,
    backgrounds,fit,positioning,shapes.symbols,chains,shadings,fadings,calc}

\tikzset{->-/.style={decoration={  markings,  mark=at position #1 with
    {\arrow{>}}},postaction={decorate}}}
\tikzset{-<-/.style={decoration={  markings,  mark=at position #1 with
    {\arrow{<}}},postaction={decorate}}}

\usepackage[all]{xy}

\usepackage{mathrsfs}
\usepackage[frak=esstix]{mathalfa}
\def\XX{\mathbb{X}}
\def\Serre{\mathbb{S}}

\theoremstyle{plain}
\newtheorem{theorem}{Theorem}[section]
\newtheorem{thm}{Theorem}
\newtheorem{lemma}[theorem]{Lemma}
\newtheorem{corollary}[theorem]{Corollary}
\newtheorem{proposition}[theorem]{Proposition}
\newtheorem{conjecture}[theorem]{Conjecture}
\theoremstyle{definition}
\newtheorem{definition}[theorem]{Definition}

\newtheorem{example}[theorem]{Example}

\newtheorem{remark}[theorem]{Remark}

\numberwithin{equation}{section}

\def\hh{\mathcal}

\def\<{\langle}
\def\>{\rangle}

\def\ZZ{\mathbb{Z}}

\def\RR{\mathbb{R}}

\def\CC{\mathbb{C}}
\def\Add{\operatorname{Add}}
\def\Aut{\operatorname{Aut}}
\def\Ind{\operatorname{Ind}}

\def\Hom{\operatorname{Hom}}

\def\Stab{\operatorname{Stab}}

\def\Br{\operatorname{Br}}

\def\arg{\operatorname{arg}}
\def\rank{\operatorname{rank}}
\newcommand{\h}{\hh{H}}            

\renewcommand{\k}{\mathbf{k}}
\renewcommand{\mod}{\operatorname{mod}}

\renewcommand{\Re}{\operatorname{Re}}

\newcommand{\id}{\operatorname{id}}

\newcommand{\D}{\operatorname{\hh{D}}}

\def\arrow{red}
\def\surf{\mathbf{S}}                       

\def\TT{\mathbf{T}}

\newtheorem{construction}[theorem]{Construction}
\newtheorem{setup}[theorem]{Set-up}

\def\gldim{\operatorname{gldim}}

\def\sli{\mathcal{P}}

\def\hua{\mathcal}

\def\tshift{\overline{\tau}}
\def\AR{\operatorname{AR}}
\def\Dwq{\D_\infty}
\def\DQ{\D_\infty(Q)}

\def\nn{node{$\bullet$}}
\def\ww{node[white]{$\bullet$}node[red]{$\circ$}}

\def\sun{Emerald}

\def\fblue{blue!20}
\def\fgreen{green!30}
\def\forange{red!30!orange!30!white}
\def\fcyan{Emerald!49!cyan!20!white}

\def\fire{orange!40}
\def\fires{yellow!20}
\def\ice{blue!15!cyan!50!white}
\def\ices{blue!30!cyan!10!white}

\def\parity{\varrho}

\def\roots{\Lambda}

\def\gms{\surf^\lambda}

\def\dynkin{\Delta}
\begin{document}

\def\Sth{\operatorname{Stgon}}
\def\hgon{\mathbf{V}}
\def\ihgon{\hgon_{\mathrm{ice}}}
\def\fhgon{\hgon_{\mathrm{fire}}}
\def\icore{\ihgon^{\text{\tiny$\copyright$}}}
\def\fcore{\fhgon^{\text{\tiny$\copyright$}}}
\def\uQ{\underline{Q}}
\def\ToSt{\operatorname{ToSt}}
\begin{abstract}
We construct a geometric model for the root category $\mathcal{D}_\infty(Q)/[2]$ of any Dynkin quiver $Q$, which is an $h_Q$-gon $\mathbf{V}_Q$ with cores, where $h_Q$ is the Coxeter number and $\mathcal{D}_\infty(Q)=\mathcal{D}^b(Q)$ is the bounded derived category associated to $Q$. As an application, we describe all spaces $\mathrm{ToSt}\mathcal{D}$ of total stability conditions on triangulated categories $\mathcal{D}$, where $\mathcal{D}$ must be of the form $\mathcal{D}_\infty(Q)$. More precisely, we prove that $\mathrm{ToSt}\mathcal{D}_\infty(Q)/[2]$ is isomorphic to a suitable moduli space of stable $h_Q$-gons of type $Q$.

In particular, an $h_Q$-gon $\mathbf{V}$ of type $D_n$ is a (centrally) symmetric doubly punctured $2(n-1)$-gon. $\mathbf{V}$ is stable if it is positively convex and the punctures are inside the level-$(n-2)$ diagonal-gon. Another interesting case is $E_6$, where the (stable) $h_Q$-gon (dodecagon) can be realized as a pair of planar tiling pattern.

    \vskip .3cm
    {\parindent =0pt
    \it Key words:
    stability conditions, root system, Dynkin diagram, geometric model}

\end{abstract}
\maketitle
\tableofcontents\addtocontents{toc}{\setcounter{tocdepth}{1}}


\section{Introduction}
The notion of stability conditions on a triangulated category $\D$ is introduced by Bridgeland \cite{B1}, whose motivation is $\Pi$-stability in string theory. They measure certain stability structure, as the name suggested, in physics as well as in mathematics.
A stability condition $\sigma=(Z,\sli)$ on $\D$ consists of a group homomorphism $Z\colon K\D\to\CC$, known as the central charge,
and an $\RR$-collection of abelian subcategories $\sli(\phi)$,
known as the slicing.
The main theorem in \cite{B1} states that all stability conditions on a triangulated category form a complex manifold.
Moreover, one recent breakthrough in this area
is the identification of stability spaces with the moduli spaces of quadratic differentials, cf. \cite{BS,HKK}.
Although there are already many studies, the known examples of global structures of stability spaces are still very limited.
\subsection{Total stability}
One key concept in the theory of stability conditions is stable object,
i.e., simple object in some $\sli(\phi), \phi\in\RR$.
It goes back to geometric invariant theory and King's $\theta$-stability,
which plays an important role in the study of
Donaldson-Thomas invariant, as well as cluster theory, cf. \cite{K11}.
When passing to the dynamical system side, i.e. identifying a stability condition with a quadratic differential on some Riemann/marked surface,
stable objects correspond to non-broken geodesics (connecting zeroes) or saddles.
An interesting question proposed in \cite{Q11} is to find stability conditions that make all indecomposable objects in the triangulated category stable.
Such stability conditions are called total stable.
The abelian version of this question is proposed by Reineke \cite{R} when studying quantum dilogarithm identities of
the module category $\mod \k Q$ of an ADE quiver $Q$.

Motivated by $q$-deformation of stability conditions in \cite{IQ1,IQ2},
Ikeda-Qiu introduce the notion of global dimension ($\RR_{\ge0}$-valued) function $\gldim$ of a stability condition $\sigma$
as a piecewise Morse-ish function. 
Qiu \cite{Q18} shows that $\sigma$ is totally stable if and only if $\gldim\sigma<1$, which is very rare.
In fact, such an existence of $\sigma$ implies that $\D$ must be the bounded derived category $\Dwq(Q)=\D^b(Q)$ of a Dynkin diagram (cf. \cite{KOT,Q20}).
In this paper, we classify all spaces $\ToSt(Q)$ of total stability conditions on the Dynkin category $\Dwq(Q)$ by constructing their geometric models.

\subsection{Geometric model for Dynkin categories}
For type $A_n$, \cite{Q18} gives a geometric description of $\ToSt\Dwq(A_n)$,
i.e. as the moduli space of positively convex $(n+1)$-gons.
It can be viewed as a variation/consequence of the fact
that $\D=\Dwq(Q)$ is the topological Fukaya category of $(n+1)$-gons as follows.

Let $\gms$ be a graded marked surface, that is a topological surface $\surf$ with marked points on its boundary $\partial\surf$ and a grading $\lambda\in\mathrm{H}_1(\mathbb{P}T\surf)$.
Given an `$\infty$-angulation' $\mathbf{A}$ of $\gms$, one obtains a graded gentle algebra $\Lambda_\mathbf{A}$ and the topological Fukaya category $\Dwq(\gms)$ can be constructed as the bounded derived category $\D^b(\Lambda_\mathbf{A})$ of $\Lambda_\mathbf{A}$.
The arcs correspond to indecomposable objects in $\D^b(\Lambda_\mathbf{A})$, cf. \cite{HKK}.
We have $\Dwq(\gms)=\Dwq(A_n)$ for a disk $\gms$ with $(n+1)$ marked points.
A geometric model for module category $\h$ of an $A_n$ quiver (with any orientation)
is given in \cite{BBMS} (and they construct a ToSt with heart $\h$ as a consequence).

If $\gms$ has punctures, which carry additional $\ZZ_2$-symmetry see \cite{QZZ}, cf. \cite{FST,S, QZ},
the story still works. Then $\Lambda_\mathbf{A}$ is a graded skew-gentle algebra.
In particular, $\Dwq(\gms)=\Dwq(D_n)$ for a once-punctured disk $\surf$ with $n-1$ marked points on $\partial\surf$.

In this paper, we introduce another model for $\Dwq(D_n)$, featuring double cover, which is similar, but different from the ones in \cite{AP,AB}.
More precisely, they take the double cover branching at punctures but we take the double cover of a once-punctured disk branching at a point other than the puncture.
What we get is a doubly punctured disk with $2(n-1)$ marked points on its boundary.
As the other models mentioned above, we can also realized objects as arcs (without tagging).
Moreover the two punctures in our model naturally correspond to the two ways of tagging in \cite{FST, QZ, QZZ}.
For instance, the tagged(-switching) rotation (introduced in \cite{BQ}) corresponding to the Auslander-Reiten (AR) translation $\tau\in\Aut\Dwq(Q)$ becomes the puncture-switching rotation in our setting. See Figure~\ref{fig:tagged} that $\tau(B_-V_3)=B_+V_2$ in both cases, where the subscripts $+/-$ denote untagged/tagged respectively in the left picture.
\begin{figure}[h]   \centering
\def\xquan{circle(.5)}
\begin{tikzpicture}[scale=.9,rotate=0]
\begin{scope}[shift={(-6,0)}]
\draw[thick,blue](0.3,0)node[right]{$\times$} to (2,0);
\draw[thick,red](0.3,0) to (-120:2)(0,0)\ww;
\foreach \j in {1,...,3}{
  \draw[thick](360/3*\j:2)\nn(360/3*\j:2.4)node{$V_\j$};
}
\draw[very thick](0,0)circle(2)(.3,0)\nn node[above]{$B$};
\end{scope}
\draw[thick,blue](.2,.4)to (2,0);
\draw[dotted](-.2,-.4)to(.2,.4);
\draw[thick,red](-0.2,-.4) to (-60:2);
\foreach \j in {1,...,6}{
  \draw[thick](360/6*\j:2)\nn(360/6*\j-180:2.4)node{$V_\j$};
}\draw(0,0)\ww;
\draw[very thick,font=\tiny](0,0)circle(2)(-.2,-.4)\nn node[left]{$B_{+}$}
    (.2,.4)\nn node[above]{$B_{-}$};
\end{tikzpicture}
\caption{tagged/puncture-switching rotation as $\tau$, type $D_4$}
\label{fig:tagged}
\end{figure}

Moreover, after straightening the model and making it (centrally) symmetric (then the geometric center coincide with the branching point),
the (oriented) arcs corresponding to the objects naturally become their central charges.
Furthermore, we can easily describe total stability using such a model, which we will mention in Section~\ref{sec:summ}.

One step further, we manage to find similar geometric models (to describe objects, central charge and total stability conditions) for all exceptional type categories, i.e. $\Dwq(E_{6,7,8})$ (see Remark~\ref{rem:sum}).
Such models shed lights on understanding the geometry of $\Stab\Dwq(Q)$ and its connection to the root systems/Kleinian singularities/Calabi-Yau categories.

\subsection{Summary of notations and results}\label{sec:summ}
A (labelled) $h$-gon $\hgon$ (in $\CC$) consists of vertices $V_0,V_1,\cdots,V_h=V_0$
and (oriented) edges $0\neq z_j=V_{j-1}V_j, \forall j\in\ZZ_h$.
We will usually use $VW$ for the (oriented) edges and $z=\overrightarrow{VW}$ for the corresponding vectors.
The $h$-gon $\hgon$ is \emph{positively convex} if all other $V_i$ is on the left hand side of
the edge $V_{j-1}V_j$.
A polygon is symmetric if it is centrally symmetric, which implies the number of vertices/edges is even.

For $1\le s\le h/2$, the (oriented) \emph{length-$s$ diagonals} of an $h$-gon are $V_jV_{j+s}$.
For instance, length-$1$ diagonals are just edges.
For a positively convex $h$-gon $\hgon$, the \emph{level-$s$ diagonal-gon}
is the convex polygon bounded by its length-$s$ diagonals (i.e. on the left hand side of).
See the hexagon bounded by the orange diagonals in Figure~\ref{fig:tt} as a level-2 diagonal-gon.

Let $Q$ be a Dynkin quiver.
\begin{definition}\label{def:h-gon}
An \emph{$h$-gon of type $Q$} is defined respectively as
\begin{description}
\item[$A_n$] an $(n+1)$-gon.
\item[$B_n$] a symmetric $2n$-gon with one puncture at its geometric center.
\item[$C_n$] a symmetric $2n$-gon.
\item[$D_n$] a symmetric doubly punctured $2(n-1)$-gon.
\item[$E_6$] a $12$-gon satisfying (4 triangle relations and 3 square relations):
\begin{gather}\label{eq:E6-rel}
    \begin{cases}
     z_j+z_{j+4}+z_{j+8}=0,\\
     z_{j}-z_{j-3}+z_{j-6}-z_{j-9}=0,
    \end{cases}\quad \forall j\in\ZZ_{12}.
\end{gather}
Note that the rank of the 7 relation equations is actually 6.
\item[$E_7$] a symmetric $18$-gon satisfying
(3 hexagon relations):
\begin{gather}\label{eq:E7-rel}
     z_{j}+z_{j+1}+z_{j+6}+z_{j+7}+z_{j+12}+z_{j+13}=0,\qquad\forall j\in\ZZ_{18}.
\end{gather}
Note that after setting $z_{j+9}=-z_j$ (by the central symmetry),
the rank of the 3 relation equations is actually 2.
\item[$E_8$] a symmetric $30$-gon satisfying
(5 triangle relations and 3 pentagon relations):
\begin{gather}\label{eq:E8-rel}
    \begin{cases}
     z_{j}+z_{j+10}+z_{j+20}=0,\\
     z_{j}+z_{j+6}+z_{j+12}+z_{j+18}+z_{j+24}=0,
    \end{cases}\quad \forall j\in\ZZ_{30}.
\end{gather}
Note that after setting $z_{j+15}=-z_j$ (by the central symmetry),
the rank of the 8 relation equations is actually 7.
\item[$F_4$] a symmetric $12$-gon satisfying \eqref{eq:E6-rel}.
Note that after setting $z_{j+6}=-z_j$ (by the central symmetry),
the rank of the 7 relation equations in \eqref{eq:E6-rel} further reduces to 4.
\item[$G_2$] a symmetric 6-gon satisfying (2 triangle relations)
\begin{gather}\label{eq:G2-rel}
     z_{j}+z_{j+2}+z_{j+4}=0,\qquad\forall j\in\ZZ_{6}.
\end{gather}
Note that after setting $z_{j+3}=-z_j$ (by the central symmetry),
the two relation equations are equivalent.
\end{description}
Note that the edge vectors also need to satisfy the condition $\sum_{j=1}^{h} z_h=0$,
which is implied by the symmetric condition or the equations above except for type $A_n$.
In each of the types except $D_n$, one can choose/fix $n$ linear independent $0\neq z_j\in\CC$ to be the local coordinates
and thus makes the space of $h_Q$-gons of $Q$ a complex submanifold of $\CC^n$.
In type $D_n$, one can freely choose all $n-1$ $\{z_j=-z_{j+n-1}\mid 1\leq j\leq n-1\}$'s
together with one vector (say $V_0B_+$ to determine the position of the two punctures) to form the $2(n-1)$-gon of type $D_n$.
Thus, the space of $h_Q$-gons of $Q$ is still a complex submanifold of $\CC^n$ in this case.
\end{definition}
In particular, an $h$-gon of type $Q$ is an $h_Q$-gon (satisfying extra conditions),
where $h_Q$ is the Coxeter number associated to the Dynkin diagram $\dynkin$.

For an exceptional case (or $F_4$ case), the relation \eqref{eq:E6-rel}, \eqref{eq:E7-rel} or \eqref{eq:E8-rel} induces a pair of $h_Q/2$-gons, called the \emph{ice} and \emph{fire cores}.
See Figure~\ref{fig:E6+-} for the two cores in type $E_6$, the upper picture of Figure~\ref{fig:E7} for the fire core in type $E_7$ and Figure~\ref{fig:E8-sthgon} for the fire core in type $E_8$.
The ice cores are the symmetric mirror of the fire ones in type $E_7$/$E_8$, respectively.

\begin{definition}\label{def:sth-gon}
An $h_Q$-gon of type $Q$ is \emph{stable} if it is positively convex and moreover:
\begin{description}
\item[$D_n$] the punctures are inside the level-$(n-2)$ diagonal-gon.
\item[$E_n$] the ice and fire cores are inside the level-$(n-3)$ diagonal-gon, for $n\in\{6,7,8\}$.
\item[$F_4$] same as $E_6$.
Note that the pair of 6-gons are symmetric to each other.
\end{description}
In particular, a stable $h$-gon means a positively convex $h$-gon without referring the types.
\end{definition}
Two key observations are: (I) the central charges of objects in any $\tau$-orbit of the root category form an $h_Q$-gon (see Proposition~\ref{pp:sthgon}) and (II) there are distinguished $\tau$-orbits, called the far-end $\tau$-orbits (cf. Lemma~\ref{lem:far-end}, as well as the mid-end/near-end $\tau$-orbits).
We always fix a far-end $\tau$-orbit and call the corresponding $h_Q$-gon the far-end one.

\begin{definition}
Denote by $\Sth(Q)$ the moduli space of stable $h_Q$-gons of type $Q$, where two stable $h_Q$-gons are equivalent
if and only if they are parallel, i.e. are related by a translation of $\CC$.

Note that the stability condition for an $h_Q$-gon is an open condition.
Hence, $\Sth(Q)$ is an open subspace of the moduli space of $h_Q$-gon,
which is a complex manifold of dimension $n$.
For instance, in type $E_6$, $\{z_j\mid 1\leq j\leq 6\}$ is a choice of local coordinate for $\Sth(E_6)$
(satisfying certain conditions, cf. Definition.~\ref{def:sth E6}).
\end{definition}

\begin{thm}\label{thm:0}
An $h_Q$-gon of type $Q$ provides a geometric model for
the root category $\Dwq(Q)/[2]$,
in the sense that it induces a central charge naturally.
Moreover, we have the isomorphism (between complex manifolds)
\begin{gather}
    Z_h\colon\ToSt(Q)/[2]\cong\Sth(Q),
\end{gather}
sending a total stability condition to the far-end stable $h_Q$-gon (of type $Q$).
\end{thm}

If further quotient by $\CC^*=\{e^{\mathbf{i} \pi s}\mid s\in\CC\}\cong\CC/[2]$, where ($s\in$)$\CC$ is the natural $\CC$-action on $\Stab$ in \eqref{eq:C-action},
then we have
\[
  Z_h\colon\ToSt(Q)/\CC\cong\mathrm{Sth}(Q),
\]
where $\mathrm{Sth}(Q)$ denotes the space of stable $h$-gons up to similarity.

A side product is that we have a simple description of the global dimension function on $\ToSt(Q)$
(Theorem~\ref{thm:gldim}).

\subsection{Connection to root systems}
\def\csub{\mathfrak{h}}
\def\hreg{\csub^{\mathrm{reg}}}
In this subsection, we explain the connection between $\ToSt(Q)$, the root system $\roots(Q)$ and
the space $\Stab\D_2(Q)$ of stability conditions on $\D_2(Q)$.
Here, $\D_2(Q)$ is the finite dimensional derived category of the Calabi-Yau-2 Ginzburg dg algebra/derived preprojective algebra associated to $Q$. It can also be constructed from coherent sheaves for the Kleinian singularity $\CC^2/G$, where $G$ is the finite subgroup of $\mathrm{SL}_2(\CC)$ corresponding to $Q$ (i.e. the McKay correspondence).

Let $\mathfrak{g}$ be the finite-dimensional complex simple Lie algebra associated to $Q$, with Cartan subalgebra $\csub$ and $\roots(Q)$ the corresponding root system.
Let $\hreg$ be the complement of the root hyperplanes in $\csub$:
\[
    \hreg=\{v\in\csub \mid v(\alpha)\neq 0, \forall \alpha\in\roots(Q) \}.
\]
The Weyl group $W$, generated by reflections of the root hyperplanes, acts freely on $\hreg$.
In \cite{B2}, Bridgeland shows
\[\Stab\D_2(Q)/\Br_Q\cong \hreg/W_Q.\]
On the other hand, one expects that
\begin{gather}\label{eq:h/W}
    \Stab\Dwq(Q)\cong \csub/W_Q
\end{gather}
and it is confirmed in \cite{HKK} for type $A$.
Clearly, there is a close relation between these two results.
For instance, one can found a detailed case study in \cite{BQS} for $A_2$ case.
For all Dynkin case, there is also a conjectural description on the almost Frobenius structure on $\Stab\Dwq(Q)$
and on the Frobenius structure on $\Stab\D_2(Q)$, which has been proved in \cite{IQ2} for type $A$.
More precisely, \cite{IQ1} introduces the Calabi-Yau-$\mathbb{X}$ category to link the Calabi-Yau-$\infty$ category and the Calabi-Yau-2 one.
By the induction-reduction procedure there, each stability condition $\sigma$ in $\ToSt(Q)$ induces a stability condition in $\Stab\D_2(Q)$.

Another connection between $\ToSt(Q)$ and the root system is the Gepner point $\sigma_G$, a stability condition with extra symmetry.
In the Dynkin case, such a point is the (unique up to $\CC$-action) solution to the equation (found by \cite{KST})
\begin{gather}\label{eq:Gepner}
    \tau(\sigma)=(-2/h_Q)\cdot\sigma,
\end{gather}
where $(-2/{h_Q})$ is the $\CC$-action on the right hand side.
The stability condition $\sigma_G$ is also the most stable one, in the sense that it is the minimal point of the global dimension function $\gldim$, cf. \cite[Thm.~4.7]{Q18}.
Interestingly, the central charge at the Gepner point for a Dynkin quiver $Q$ is given by the projection of the root system on the Coxeter plane (Lemma~\ref{lem:L.H}).
This was pointed out by Lutz Hille to Qy and Alastair King in Oberwolfalch, Jan. 2020.
We draw many pictures in Appendix~\ref{app} showing such projections (together with features of characterization of total stability conditions, i.e. the stable $h_Q$-gon).
Here is the trailer:
\[
\makebox[\textwidth][c]{
\includegraphics[width=1.5in]{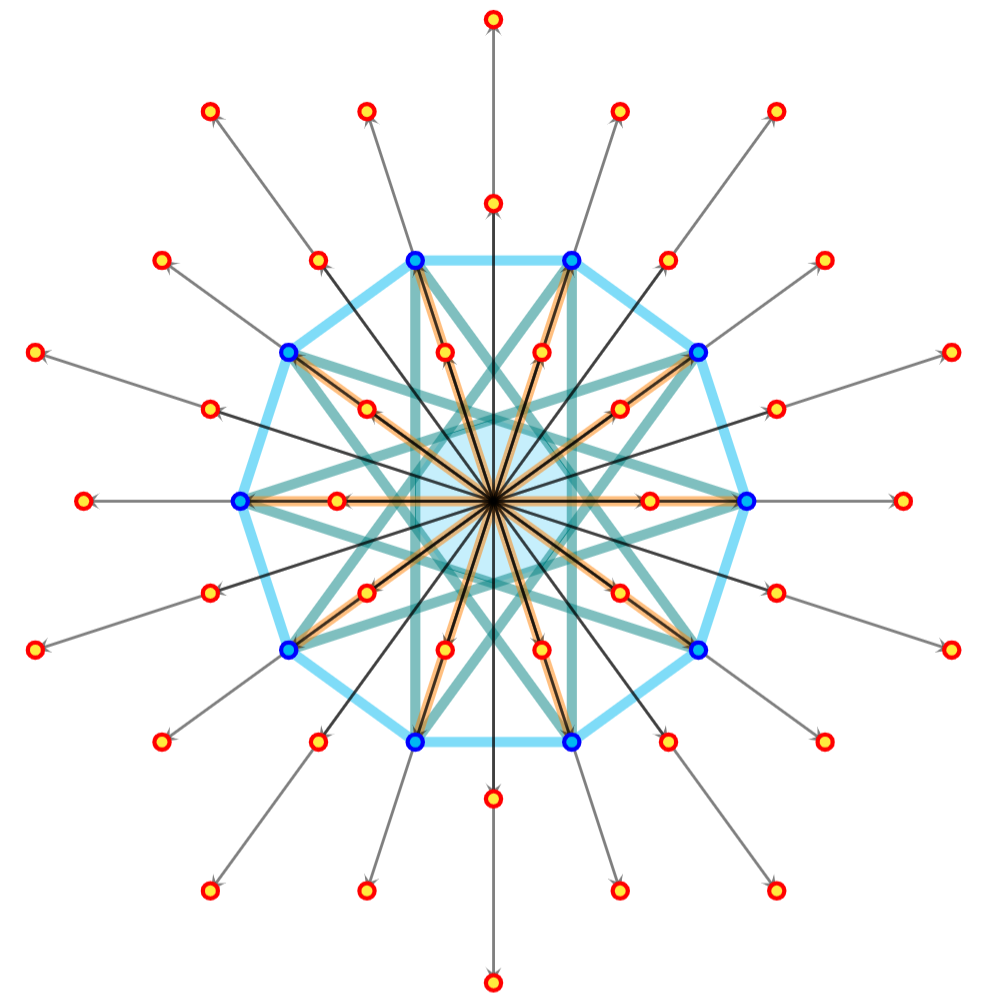}
\includegraphics[width=1.5in]{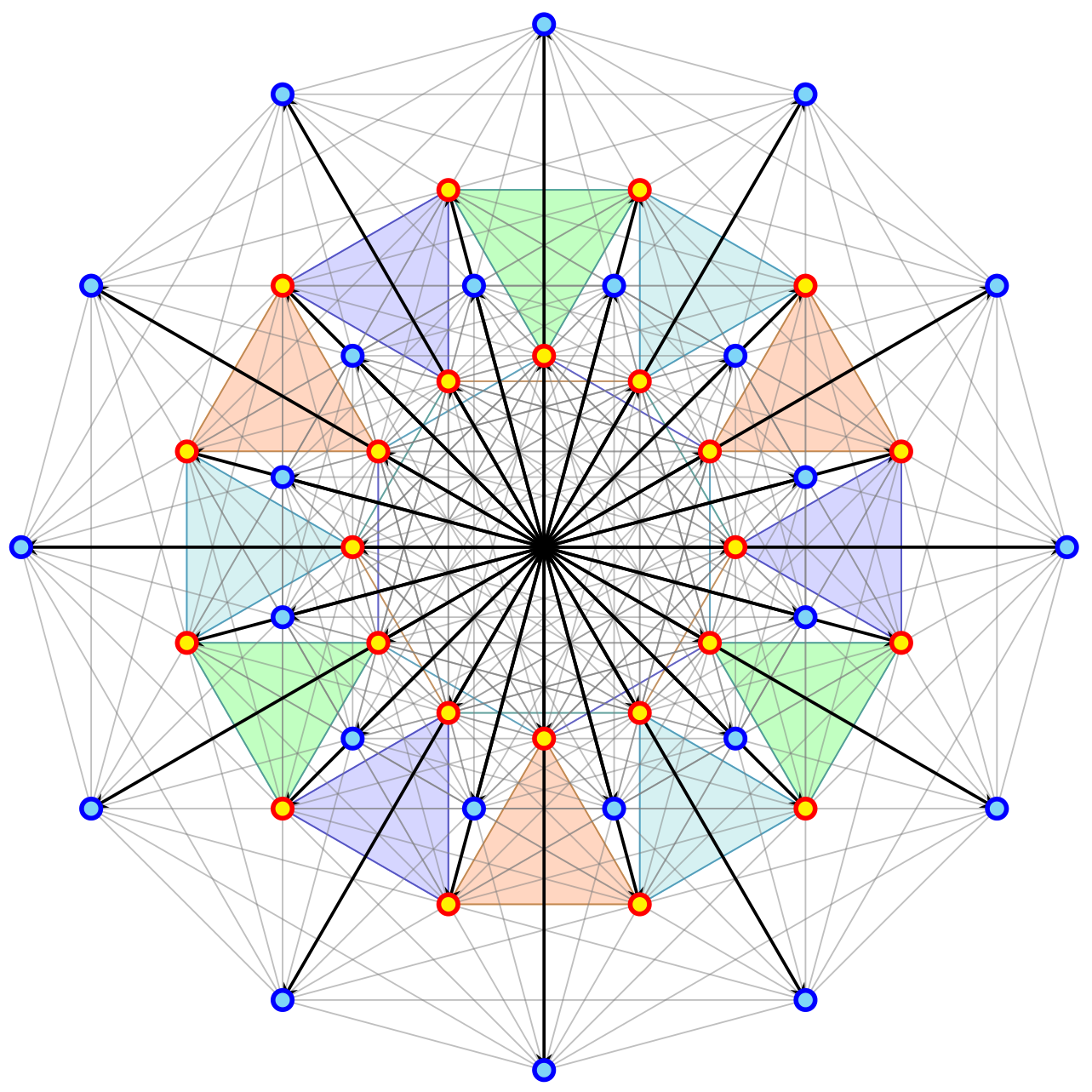}
\includegraphics[width=1.5in]{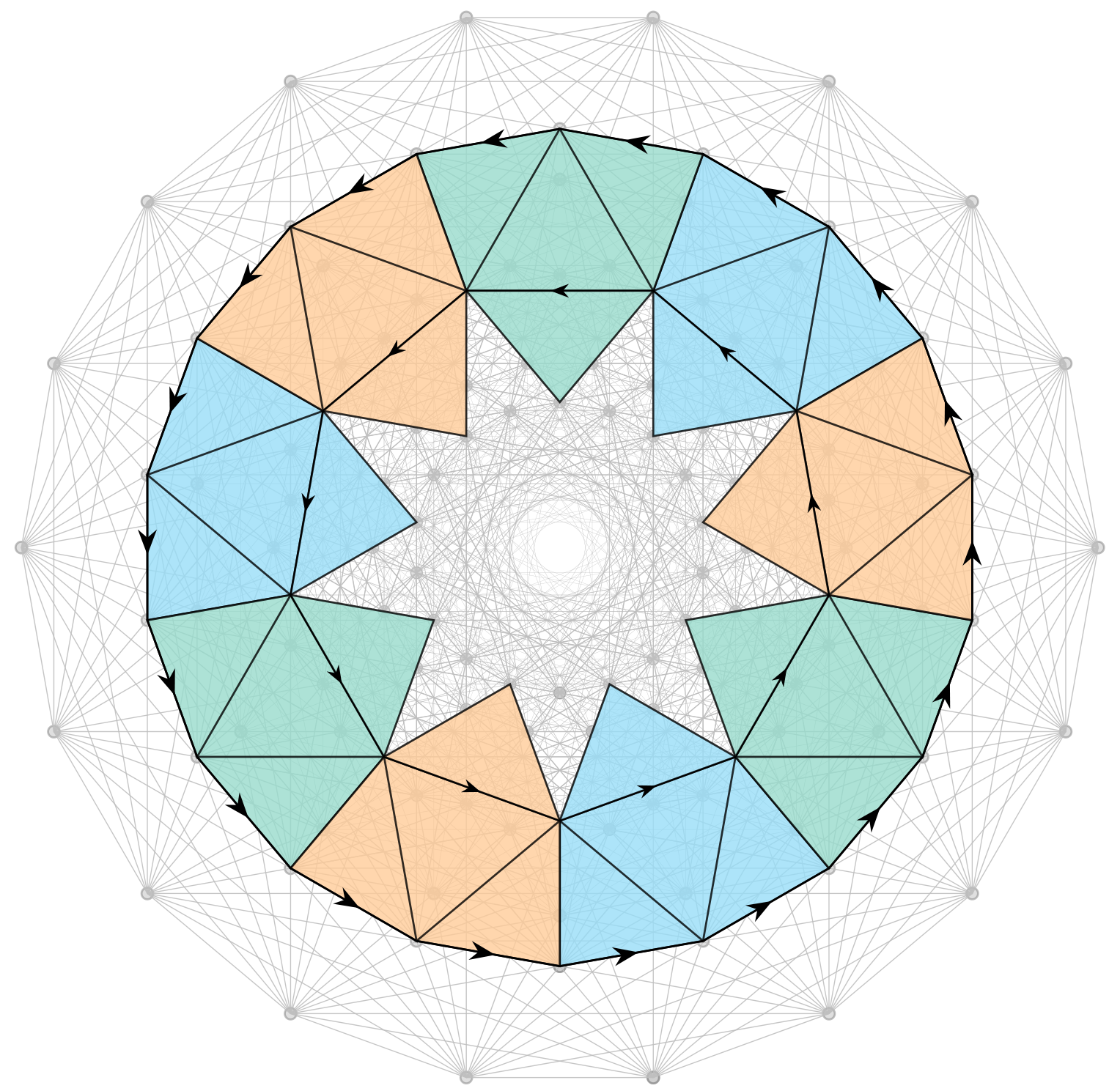}
\includegraphics[width=1.5in]{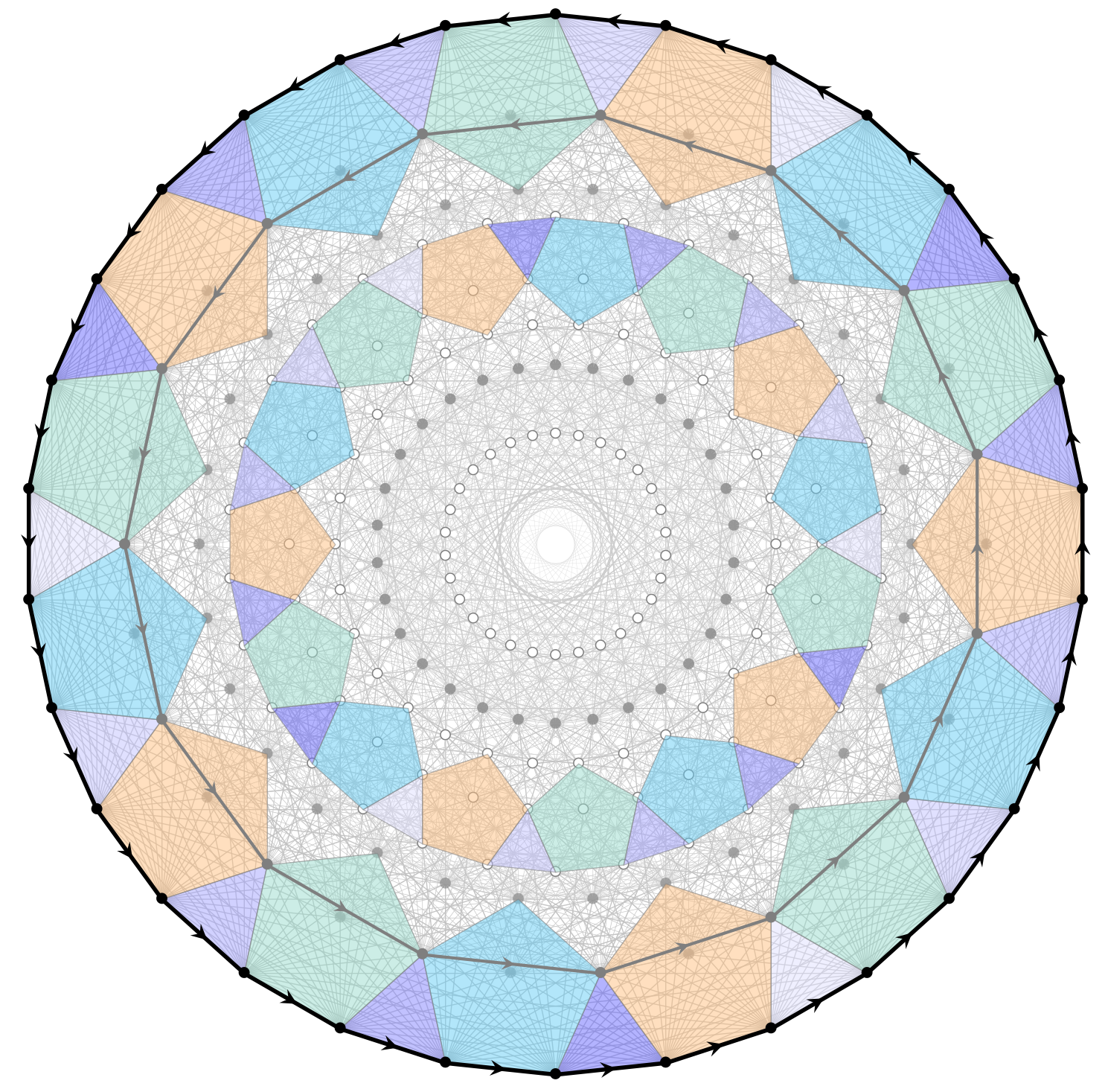}
}\]

\subsection{Further studies}
A direct byproduct of our description of $\ToSt(Q)$ is to prove Reineke's conjecture.
Note that as we describe all total stability functions on any heart of $\DQ$,
which is a global result (as finding one total stability function for a chosen heart is really a local result).
\begin{conjecture}\cite{R,K11}
For any Dynkin quiver $Q$, there is a stability function $Z$
on $\h_Q=\mod \k Q$ such that all indecomposable objects are stable.
\end{conjecture}
Notice that a total stability condition induces a stability function on its heart
that makes all indecomopsable objects in the heart stable.
Thus, to prove this conjecture, we only need to understand the heart of a total stability condition and
construct a stable $h_Q$-gon so that its heart is equivalent to any given $\h_Q$.
This is done in the sequel \cite{CQZ}.
Very recently, an equivalent conditions for a total stability function on module category of a Dynkin quiver is given in \cite{DGK}.

Another follow-up work is to attack the following conjecture
(we state this for $\Dwq(Q)$ but one could consider general case).
\begin{conjecture}\cite{Q18,Q20}
For any Dynkin quiver $Q$,
the global dimension function $\gldim$ is a piece-wise Morse-ish function on $\Stab\Dwq(Q)/\CC$.
Its only critical point is the Gepner point and hence induces a contractible flow on $\Stab\Dwq(Q)$.
\end{conjecture}\emph{}
Note that (cf. \cite{Q18}) we have $\ToSt(Q)=\gldim^{-1}[1-2/h_Q,1)$.
So $\ToSt(Q)$ is the core of $\Stab\Dwq(Q)$.
The conjecture above decomposes into two parts:
\begin{enumerate}
    \item $\gldim$ is piece-wise Morse-ish on $\Stab\Dwq(Q)\setminus\ToSt(Q)$ and
thus $\Stab\Dwq(Q)$ contracts to $\ToSt(Q)$.
    \item $\gldim$ is piece-wise Morse-ish on $\ToSt(Q)$.
\end{enumerate}
The first part is proved in \cite{Q20} for graded gentle algebras (and in particular for type A).
Presumably, the method there should generalize to graded skew-gentle algebras (with extra effort to dealing with punctures, and in particular for type D) once the geometric model sets up (e.g. \cite{QZZ}).
Now with our geometric model for all Dynkin types, one may try to prove the first part for type E as well.
Combining with the expectation above, i.e. \eqref{eq:h/W},
the first part says that there is a canonical contraction:
$\mathfrak{h}/W_Q\to\ToSt(Q)$.

The second part is much more tricky, even in type A. The statement (say in type A) is easy to understand: the derivative of the function $\gldim$ gives a canonical way to deform any convex polygon into a regular polygon.
This part is specific to the Dynkin case, while the phenomenon of the first part is much more general (e.g. for coherent sheaves on a complex line in \cite{Q18} and for coherent sheaves on a complex plane in \cite{FLLQ}).

Another potential application of our geometric model is to study tilted algebras of type E.
The tilted algebras of type A and D lead to gentle and skew-gentle algebras, where the corresponding geometric models can be glued together to produce the topological Fukaya categories of (graded marked) surfaces, as mentioned above. It would be interesting to see how could these work in type E case.
Also, ToSt corresponds to $\gldim<1$ and one may want to study the ToSS=total semi-stability case,
which corresponds to $\gldim\le1$.

Finally, we mention a related work \cite{H},
which gives also relates polytopes with tilting modules (but restricted to type A).
One could try to generalize the results there using our geometric model.

\subsection*{Acknowledgments}
Qy would like to thank Alastair King for pushing him to think carefully about total stability and suggesting the abbreviation ToSt among numerous comments.
We would also like to thank Yu Zhou and Shiquan Ruan for pointing out references for Coxeter elements and Wen Chang for proofreading.
We are supported by National Natural Science Foundation of China 
(Grant No.12425104, No.12101422 and No.12031007)
and National Key R\&D Program of China (No. 2020YFA0713000).


\section{Preliminaries}\label{sec:pre}
\subsection{Stability conditions}\label{sec:BSC}
Following \cite[Def.~1.1]{B1}, we recall the notion of stability conditions on a triangulated category.
In this paper, $\D$ is a triangulated category with Grothendieck group
$K\D$ and assume that $K\D\cong\ZZ^n$ for some $n$.
Denote by $\Ind\D$ the set of (isomorphism classes of) indecomposable
objects in $\D$.

\begin{definition}
\label{def:stab}
A {\it stability condition} $\sigma = (Z, \sli)$ on $\D$ consists of
a group homomorphism $Z \colon K\D \to \CC$, called the {\it central charge}, and
a family of full additive subcategories $\sli (\phi) \subset \D$ for $\phi \in \RR$, called the {\it slicing},
satisfying the following conditions:
\begin{itemize}
\item[(a)]
if  $0 \neq E \in \sli(\phi)$,
then $Z(E) = m(E) e^{ \mathbf{i} \pi \phi }$ for some $m(E) \in \RR_{>0}$,
\item[(b)]
for all $\phi \in \RR$, $\sli(\phi + 1) = \sli(\phi)[1]$,
\item[(c)]if $\phi_1 > \phi_2$ and $A_i \in \sli(\phi_i)\,(i =1,2)$,
then $\Hom(A_1,A_2) = 0$,
\item[(d)]for each object $0 \neq E \in \D$, there is a finite sequence of real numbers
\begin{equation}\label{eq:>}
\phi_1 > \phi_2 > \cdots > \phi_l
\end{equation}
and a collection of exact triangles (known as the \emph{HN-filtration})
\begin{equation*}
0 =
\xymatrix @C=5mm{
 E_0 \ar[rr]   &&  E_1 \ar[dl] \ar[rr] && E_2 \ar[dl]
 \ar[r] & \dots  \ar[r] & E_{l-1} \ar[rr] && E_l \ar[dl] \\
& A_1 \ar@{-->}[ul] && A_2 \ar@{-->}[ul] &&&& A_l \ar@{-->}[ul]
}
= E
\end{equation*}
with $A_i \in \sli(\phi_i)$ for all $i$.
\item[(e)] a technique condition, known as the support property, which holds automatically in our setting.
\end{itemize}
The categories $\sli(\phi)$ are then abelian. Their non zero
objects are called {\it semistable of phase $\phi$} and simple objects
{\it stable of phase $\phi$}.
For a semistable object $E\in\sli(\phi)$,
denote by $\phi_\sigma(E)\colon\!\!=\phi$ its phase.
\end{definition}

There is a natural $\CC$-action on the set $\Stab\D$ of all stability conditions on $\D$, namely:
\begin{equation}\label{eq:C-action}
    s \cdot (Z,\hh{P})=(Z \cdot e^{-\mathbf{i} \pi s},\hh{P}_{\Re(s)}),
\end{equation}
where $\hh{P}_x(\phi)=\hh{P}(\phi+x)$ and $s\in\CC$.
There is also a natural action on $\Stab\D$ by the group of
autoequivalences $\Aut\D$, namely:
$$\Phi  (Z,\hh{P})=\big(Z \circ \Phi^{-1}, \Phi (\hh{P}) \big),$$
where $\Phi\in\Aut\D$.

The famous result in \cite{B1} states that $\Stab\D$ is a complex manifold with dimension $\rank K\D$
and the local coordinate is provided by the central charge $Z$.

\subsection{Total (semi)stability via global dimension function}
Let $\sigma=(Z,\sli)\in\Stab\D$.
\begin{definition}
The {\em global dimension of $\sigma$} is defined as
\begin{gather}\label{eq:geq}
\gldim\sigma=\sup\{ \phi_2-\phi_1 \mid
    \Hom(\hh{P}(\phi_1),\hh{P}(\phi_2))\neq0\}\in\RR_{\ge0}\cup\{+\infty\},
\end{gather}
which is a continuous function and an invariant under both the $\CC$-action and $\Aut\D$.
The global dimension of $\D$ is defined to be the $\inf$ of $\gldim\sigma$ for all $\sigma\in\Stab\D$.
\end{definition}
Note that the notion generalizes the global dimension of an algebra/abelian category.

We recall the notion of total (semi)stability on triangulated categories,
whose abelian version is due to Reineke \cite{R},
cf. \cite[Conjecture~7.13]{Q11} and comments there.
\begin{definition}
A stability condition $\sigma$ is called {\em totally (semi)stable},
if any indecomposable object in $\D$ is (semi)stable with respect to $\sigma$.
We will call it a total stability condition for short.
\end{definition}

\subsection{Dynkin diagrams and KOT-Q classification}\label{sec:KOT-Q}
A \emph{simply-laced Dynkin quiver} $Q$ is an oriented graph whose underlying graph $\dynkin$ is one of ADE Dynkin diagram.
Explicitly, $\dynkin$ is of the form $T_{p,q,r}$ (with $p+q+r-2$ vertices, cf. \cite{Gab}):
\begin{gather}\label{eq:pqr}
\begin{tikzpicture}[scale=1.5]
\draw(0,0)node[right]
    {$\underbrace{\bullet-\bullet-\cdots-\bullet-\circ}_r$ \hskip .5pc the far-end vertex}
    (-.64,.72)node[rotate=150]
    {$\underbrace{\bullet-\bullet-\cdots-\bullet-\bullet}$}
    node[above right]{$^p$}
    (-.727,-.25)node[rotate=-150]
    {$\underbrace{\bullet-\bullet-\cdots-\bullet-\bullet}$}
    node[above left]{$^q$};
\end{tikzpicture}
\end{gather}

with $1\le p\le q \le r$ and
\[
    \frac{1}{p}+\frac{1}{q}+\frac{1}{r}>1.
\]
Then we have ($n=p+q+r-2$)
\begin{itemize}
\item $Q$ is of type $A_n$ if $\dynkin=T_{1,q,r}$ for $q+r=n+1$.
\item $Q$ is of type $D_{n}$ if $\dynkin=T_{2,2,n-2}$ for $n\ge4$.
\item $Q$ is of type $E_{n}$ if  $\dynkin=T_{2,3,n-3}$ for $n\in\{6,7,8\}$.
\end{itemize}

\begin{definition}\label{def:leave}
For a simply-laced Dynkin quiver $Q$,
we call its leaves (i.e. univalent vertices) the \emph{boundary vertices}.
Moreover, we call
\begin{itemize}
  \item the leaf at the end of length $p$ branch of $Q$ the \emph{near-end vertex}.
  \item the leaf at the end of length $q$ branch of $Q$ the \emph{mid-end vertex}.
  \item the leaf at the end of length $r$ branch of $Q$ the \emph{far-end vertex}.
\end{itemize}
Note that there is a choice involved for fixing the far-end vertex in type $A_n$, $D_4$ and $E_6$.
\end{definition}

Denote by $\Dwq(Q)=\D^b(\k Q)$ the bounded derived category of the path algebra $\k Q$ for a simply-laced Dynkin quiver $Q$.
Denote by $\AR\Dwq(Q)$ the \emph{AR quiver} of $\Dwq(Q)$,
which is isomorphic to $\ZZ Q$.
Each vertex of $Q$ canonically corresponds to a \emph{$\tau$-orbit} (containing the corresponding projective).
The $\tau$-orbit that corresponds to an xx vertex
is call xx \emph{$\tau$-orbit}, for xx being boundary/near-end/mid-end/far-end.

The following is well-known.
\begin{lemma}\label{lem:far-end}
For type $A_n$ and $E_6$, there are two choices of far-end $\tau$-orbits
of $\AR\Dwq(Q)$, which are shift $[1]$ of each other.
For type $D_4$, there are three choices of far-end $\tau$-orbit of $\AR\Dwq(Q)$
(i.e. any boundary $\tau$-orbit).
For other cases, there is a unique choice of the far-end $\tau$-orbit.
\end{lemma}
For the orbit categories of $\Dwq(Q)$ (e.g. the root category $\Dwq(Q)/[2]$ we are about to mention in particular), we keep the same notions.
We will always choose a preferred far-end vertex/$\tau$-orbit
in each case.

A non-simply-laced Dynkin quiver, denoted by $R^\iota$, is defined to be
a simply-laced Dynkin quiver $R$ together with an automorphism $\iota$ of $R$.
All possible cases are listed as follows (we omit the orientations but they should be compatible with $\iota$):
\def\S{R^\iota}
\begin{description}
\item[$B_n$]
    $R$ is of type $D_{n+1}$ and $\S$ is of type $B_n$ while
    $\iota$ exchanges the black bullets in the same column.
\vskip.05cm
\[
    D_{n+1} \quad
    \xymatrix@R=0.1pc@C=0.7pc{
        &&&& \bullet &\ar@/^/@{<->}[dd]^{\iota}\\
        \circ \ar@{-}[r]& \circ \ar@{-}[r]& \cdots \ar@{-}[r]& \circ \ar@{-}[dr]\ar@{-}[ur]\\
        &&&& \bullet &\\
    }
    \qquad
    B_n \quad
    \xymatrix@R=0.1pc@C=0.7pc{
        \\
        \circ \ar@{-}[r]& \circ \ar@{-}[r]& \cdots \ar@{-}[r]& \circ \ar@{-}[r]& \bullet\\
        _1 & _1 && _1 & _\mathbf{2}
    }
\]
\item[$C_n$]
    $R$ is of type $A_{2n-1}$ and $\S$ is of type $C_n$ while
    $\iota$ exchanges the black bullets in the same column.
\vskip.05cm
\[
    A_{2n-1} \quad
    \xymatrix@R=0.1pc@C=0.7pc{
        \ar@/_/@{<->}[dd]_{\iota} & \bullet  \ar@{-}[r]& \bullet \ar@{-}[r]&
            \cdots \ar@{-}[r]& \bullet \ar@{-}[dr]\\
        &&&&& \circ \\
        &\bullet \ar@{-}[r]& \bullet \ar@{-}[r]&
            \cdots \ar@{-}[r]& \bullet \ar@{-}[ur]
    }
    \qquad
    C_n \quad
    \xymatrix@R=0.1pc@C=0.7pc{
        \\
        \bullet \ar@{-}[r]& \bullet \ar@{-}[r]& \cdots \ar@{-}[r]& \bullet \ar@{-}[r]& \circ\\
        _\mathbf{2} & _\mathbf{2} && _\mathbf{2} & _1
    }
\]
\item[$F_4$]
    $R$ is of type $E_6$ and $\S$ is of type $F_4$ while
    $\iota$ exchanges the black bullets in the same column.
\vskip.05cm
\[
    E_6 \quad
    \xymatrix@R=0.1pc@C=0.7pc{
        && \bullet \ar@{-}[r]& \bullet & \ar@/^/@{<->}[dd]^{\iota}\\
        \circ \ar@{-}[r]& \circ \ar@{-}[dr]\ar@{-}[ur]\\
        && \bullet \ar@{-}[r]& \bullet &\\
    }
    \qquad
    F_4 \quad
    \xymatrix@R=0.1pc@C=0.7pc{
        \\
        \circ \ar@{-}[r]& \circ \ar@{-}[r]& \bullet \ar@{-}[r]& \bullet\\
        _1 & _1 & _\mathbf{2} & _\mathbf{2}
    }
\]
\item[$G_2$]
    $R$ is of type $D_4$ and $\S$ is of type $G_2$ while
    $\iota$ rotates the black bullets in the same column.
\vskip.05cm
\[
    D_4 \quad
    \xymatrix@R=0.1pc@C=0.7pc{
        & \bullet & \ar@/^/@{<->}[dd]^{\iota}\\
        \circ \ar@{-}[r]\ar@{-}[dr]\ar@{-}[ur]& \bullet\\
        & \bullet &\\
    }
    \qquad
    G_2 \quad
    \xymatrix@R=0.1pc@C=0.7pc{
        \\
        \circ \ar@{-}[r]& \bullet\\
        _1 & _\mathbf{3}
    }
\]
\end{description}
The label below any black bullet in $\S$ is its weight (i.e. number of vertices in that orbit).

In the non-simply-laced case $Q=R^\iota$, 
$\DQ$ means the $\iota$-stable category $\Dwq(R)^\iota$
and the boundary/far-end vertices (of $Q$) and $\tau$-orbits (of $\AR\Dwq(Q)$) are induced from $R$.
Lemma~\ref{lem:far-end} implies that there is a unique far-end $\tau$-orbit in any non-simply-laced case.

Denote by $\ToSt\D$ the subspace of $\Stab\D$
consisting of all total stability conditions on a triangulated category $\D$.
The classification theorem below is a combination of the following two results:
\begin{itemize}
    \item \cite[Prop.~3.5]{Q18}: $\sigma$ is totally stable if and only if $\gldim\sigma<1$.
    \item \cite[Thm.~3.2]{Q20}: $\D$ admits a stability condition $\sigma$ with $\gldim\sigma<1$ if and only if $\D$ equals $\Dwq(Q)$ for a Dynkin diagram $\dynkin$.
\end{itemize}
Note that the second result is essentially due to Kikuta-Ouchi-Takahashi, where they impose a mild condition that excludes the non simply-laced case.

\begin{theorem}
$\ToSt\D$ is non-empty if and only if $\D$ equals $\Dwq(Q)$ for a Dynkin quiver $Q$.
In such a case, we write $\ToSt(Q)$ for $\ToSt\Dwq(Q)$.
\end{theorem}

\subsection{Root system}
Let $\mathfrak{g}$ be the complex simple Lie algebra of type $Q$
and $\mathfrak{h}$ its Cartan subalgebra.
Denote by $\roots(Q)$ the associated \emph{root system} and $h_Q$ the \emph{Coxeter number}.

The \emph{Coxeter element} $w$ is the product of all simple reflections.
Although it depends on the order of the product,
all such elements are conjugate to each other.
Thus, up to symmetry, there is essentially one Coxeter element.
After we fix the Coxeter element $w$, there is a unique plane $P_w$,
known as the \emph{Coxeter plane},
on which $w$ acts by rotation by $2\pi/h_Q$.

For instance, the projection of the root system on the Coxeter plane $P_w$ of type $D_5$, $E_6$, $E_7$ and $E_8$ are shown in Figures~\ref{fig:D5},~\ref{fig:E6-RS},~\ref{fig:E7-1} and~\ref{fig:E8-RS} respectively.

\subsection{Root categories}
Note that the dimension function $\dim$ depends on the orientation of Q but the root category does not.
Recall the famous Gabriel's theorem that
the root category categorifies the root system,
in the sense that there is a bijection
\begin{gather}\label{eq:G}
    \dim\colon\Ind\Dwq(Q)/[2]\to \roots(Q).
\end{gather}
Set $\tshift=\tau^{-1}$.
It is well-known that
\[
    \tshift^{h_Q}=[2]
\]
and $\Serre=[1]\circ\tau$ is the Serre functor.
The following observation is due to L. Hille.

\begin{lemma}\label{lem:L.H}
Under the bijection in \eqref{eq:G},
the projection of the root system on the Coxeter plane
gives the central charge of the Gepner point $\sigma_G$ of $\Dwq(Q)$.
\end{lemma}
\begin{proof}
As the Coxeter element $w$ acts by rotation by $2\pi/h_Q$,
the roots with the same length spread evenly in $P_w$.
This matches the Gepner equation \eqref{eq:Gepner}
and implies the lemma. See the top pictures in Figures~\ref{fig:E6-RS} and~\ref{fig:E7-1}.
\end{proof}

\section{Prototype of $\ToSt$ (i.e. type $A_n$)}\label{sec:A}
\subsection{The $h_Q$-gons induced by $\tau$-orbits}
\begin{lemma}\label{lem:z}
$\id-\tau$ (or $\id+\Serre$) is a non-singular linear map on $\RR^n=K\Dwq(Q)\otimes \RR$.
\end{lemma}
\begin{proof}
Taking the basis $\{ \dim P_i \mid i\in Q_0 \}$ of $\RR^n$ for $P_i$ being the projectives of $Q$.
As $\tau(P_i)=I_i[-1]$ for $I_i$ being the injectives,
the matrix presentation $C(Q)$ of $\tau$, as a linear transformation on $K\Dwq(Q)$ is the so-called \emph{Coxeter transformation} of $Q$.
It is well-known that the eigenvalues of the Coxeter transformation is not 1 for Dynkin type (cf. \cite{L}).
Thus the lemma follows.
\end{proof}

For any $M$ in $\Ind\Dwq(Q)/[2]$, define
\[g(M)\colon=\sum_{j=1}^{h_Q} [\tau^{j} M],\]
to be the sum of all objects in the $\tau$-orbit of $M$ in $K\Dwq(Q)$.

\begin{proposition}\label{pp:sum=0}
The function $g\equiv0$.
\end{proposition}
\begin{proof}
Note that $\tau$ is a group automorphism of $K\Dwq(Q)$.
Then $\tau^{h_Q}=[-2]=\id$ implies
\[
    (\id-\tau)(\id+\tau+\cdots\tau^{h_Q-1})=0.
\]
Then Lemma~\ref{lem:z} implies the proposition.

An alternating categorical proof is as follows, which fits better in our content.
We claim that any $g(P_i), i\in Q_0$ are proportional.
Note that we only need to check this for the simply laced case,
as the other cases can be obtained from folding and hence hold also.
Then
\[ \sum_{i\in Q_0} g(P_i)=\sum_{M\in\Ind\Dwq(Q)/[2]} [M]
=\sum_{M\in\Ind\hua{H}_Q} \Big([M]+[M[1]] \Big)=0\]
implies the proposition.
For the claim, let $i$ be a leaf of $Q$
(e.g. the boundary vertex, cf. Definition~\ref{def:leave})
and $j$ be its neighbour vertex.
Summing all the mesh relations in $K\Dwq(Q)$ at the $\tau$-orbit of $P_i$,
which are triangle relations, we have
\[ 2g(P_i)=g(P_j). \]
If $j$ is not the tri-valent vertex in \eqref{eq:pqr},
let $k$ be its neighbour vertex other than $i$.
Summing all the mesh relation in $K\Dwq(Q)$ at the $\tau$-orbit of $P_i$,
which are square relations, we have
\[ 2g(P_j)=g(P_i)+g(P_k) \]
and hence $g(P_k)=3g(P_i)$.
Do all such calculations, we see the claim holds.
\end{proof}

Let $Z$ be any central charge and recall that $\tshift=\tau^{-1}$. Then we have the following fact:
\begin{corollary}\label{cor:hgon}
For any $M\in\Ind\Dwq(Q)$, the vectors $z_1=Z(\tshift M),z_2=Z(\tshift^2 M),\ldots,$ $z_{h_Q}=Z(\tshift^{h_Q} M)$ form an $h_Q$-gon.
\end{corollary}

\begin{definition}\label{def:far-end hgon}
The $h_Q$-gon (with respect to any given central charge $Z$) induced by a boundary/$?$-end $\tau$-orbit (of $\AR\Dwq(Q)/[2]$) is called the \emph{boundary/$?$-end $h_Q$-gon}, for $?$ being mid/near/far.
\end{definition}

Once we choose our favorite far-end vertex/$\tau$-orbit,
we will have the corresponding far-end $h_Q$-gon.
It will not make much difference in type $A_n$ and $E_6$.
But as we will explain in Section~\ref{sec:D4},
the choice matters in type $D_4$.

\subsection{Positive convexity as stability for $h_Q$-gons}
\begin{lemma}\label{lem:easy}
Let $\sigma=(Z,\sli)\in\ToSt(Q)$.
If there is a (non-trivial) path from $M$ to $L$ in the AR quiver of $\DQ$,
then $\phi_\sigma(M)<\phi_\sigma(L)$.
\end{lemma}
\begin{proof}
Any path in the AR-quiver consists of irreducible morphisms (which are non-zero).
Hence, the total stability implies that the phases of (stable) objects in the path are increasing.
\end{proof}

\begin{proposition}\label{pp:sthgon}
If $\sigma=(Z,\sli)\in\ToSt(Q)$,
then, for any $M\in\Ind\Dwq(Q)$, the vectors $Z(\tshift M),Z(\tshift^2 M),\ldots,Z(\tshift^{h_Q} M)$
form a positively convex $h_Q$-gon.
\end{proposition}
\begin{proof}
In the $A_1$ case, the statement is trivial. Assume that $n=|Q_0|\ge2$ in the following.

In the AR quiver $\Dwq(Q)$,
there is a path from $M$ to $\tshift M$ and hence $\phi_\sigma(M)<\phi_\sigma(\tshift M)$ by Lemma~\ref{lem:easy}.
Then
\[
    \phi_\sigma(M)<\phi_\sigma(\tshift M)<
        \cdots <\phi_\sigma(\tshift^{h_Q-1} M) < \phi_\sigma(\tshift^{h_Q} M)=\phi_\sigma(M[2])=\phi_\sigma(M)+2.
\]
Let $z_i=Z(\tshift^i M), \forall i\in\ZZ_{h_Q}$. Then $\arg z_0=\pi\cdot\phi_\sigma (M)$.
So we have
\[
    \arg z_0 < \arg z_1 < \cdots < \arg z_{h_Q-1} <
    \arg z_{h_Q}+2\pi =\arg z_0 + 2\pi,
\]
where $\arg$ takes values in $[\pi\cdot\phi_\sigma (M) , \pi\cdot\phi_\sigma (M)+2\pi)$.
It implies the $h_Q$-gon with edges $z_i$ is indeed positively convex.
\end{proof}

Up to translation, each $\tau$-orbit induces one positively convex $h_Q$-gon for
a chosen $\sigma\in\ToSt(Q)$.
We call all such $h_Q$-gons the \emph{$\sigma$-induced positively convex $h_Q$-gons}.

A direct consequence is the following characterization of $\gldim$ on $\ToSt(Q)$, which generalizes \cite[Prop.~3.6]{Q18} for type $A_n$ to all Dynkin cases.
\begin{theorem}\label{thm:gldim}
If $\sigma\in\ToSt(Q)$,
then $\pi\cdot\gldim\sigma$ equals the maximal angle among any interior angles of $\sigma$-induced positively convex $h_Q$-gons.
\end{theorem}
\begin{proof}
By AR duality, we have $\Hom(E,F)=D\Hom(F,\tau(E[1]))$.
So in particular, $\Hom(E,\tau(E[1]))=D\Hom(E,E)\neq0$.
Moreover, by the hammock property of the AR-quiver of $\DQ$, we know that $\Hom(E,F)\neq0$ implies that
there is a path from $E$ to $F$ to $\tau(E)[1]$ in the AR-quiver, 
for any indecomposables $E,F$.
By Lemma~\ref{lem:easy},
\[
    \phi_\sigma(E)\le \phi_\sigma(F)\le \phi_\sigma\big(\tau(E[1])\big),
\]
and thus $\phi_\sigma(F)-\phi_\sigma(E)\le \phi_\sigma\big(\tau(E[1])\big)-\phi_\sigma(E)$.
Hence,
\[\begin{array}{rl}
    \gldim\sigma&=\sup\{ \phi_\sigma\big(\tau(E[1])\big)-\phi_\sigma(E) \mid E\in\Ind\Dwq(Q) \}\\
    &=\max\{ \phi_\sigma\big(\tau(E[1])\big)-\phi_\sigma(E) \mid E\in\Ind\Dwq(Q)/[2] \}.
\end{array}\]
Finally, notice that $\pi\cdot\left(\phi_\sigma(\tau(E[1]))-\phi_\sigma(E)\right)$
is indeed an interior angle of the $h_Q$-gons between the diagonals corresponding to $E$ and $\tau(E[1])$.
Thus, the statement follows.
\end{proof}

\subsection{$\ToSt$ of type $A_n$}
Although the orientation of the quiver does not matter when considering $\ToSt(Q)$, we will choose our favorite one for convenience (e.g. to describe objects).
\begin{setup}\label{setup:A}
Take the $A_n$ quiver with straight orientation
\begin{gather}\label{eq:An}
\begin{tikzpicture}[scale=1,xscale=1,ar/.style={->,>=stealth,thick}]
\foreach\j in {0,2,4} {\draw[ar](\j+.5,0)to++(-1,0);}
\draw(-1,0)node{$1$}(1,0)node{$2$}(3,0)node{$\cdots$}
    (5,0)node{$n$};
\end{tikzpicture}
\end{gather}
We declare vertex 1 as the fixed far-end vertex
and label the indecomposable objects in $\Dwq(A_n)$ by
\begin{gather}\label{eq:Pji}
    P^{j}_{i}\colon=\tshift^{j}P_i,\quad j\in\ZZ,1\le i\le n,
\end{gather}
where each $P_i$ is the projective corresponding to vertex $i$.
The AR quiver of $\Dwq(A_n)/[2]$, for $n=5$, is illustrated as in Figure~\ref{fig:AR An},
where $j (j+i)$ denotes $P^j_i$ and $(j+i) j$ denotes $P^j_i[1]$ for $1\le i\le n$ and $0\le j\le n-i$.
\end{setup}
Note that the oriented diagonals of an $(n+1)$-gon is a model for indecomposable objects in the root category of type $A_n$ is well-known, cf. e.g. \cite{BM}.
\begin{figure}
\def\xquan{circle(.5)}
\begin{tikzpicture}[xscale=.7,yscale=.7,rotate=-45]
\clip[rotate=-45](-.8,2)rectangle(6.3,19);
\foreach \j in {-3,...,8}{
\foreach \k in {0,...,8}{
    \draw[teal,->,>=stealth,very thick]
        (2*\k,2*\j+.7)to++(0,.6);
    \draw[teal,->,>=stealth,very thick]
        (2*\k+.7,2*\j)to++(.6,0);
}}
\draw[rotate=-45,white,fill=white](0,2)rectangle(-2,19)(8,2)rectangle(6,19);
\foreach \j in {-3,...,8}{
\foreach \k in {0,...,8}{
    \draw[white](2*\k,2*\j) node[circle,draw=teal] (x\k\j) {$00$};
}}
\foreach \j in {0,...,4}{
\foreach \k in {\j,...,5}{
    \draw[](2*\k,2*\j+2)node[circle,draw=teal,fill=teal!10] {$\j\k$};
}}
\draw[gray!80]
    (4,0)node{$52$}(6,0)node{$53$}(8,0)node{$54$}(6,-2)node{$43$}
    (12,12)node{$50$}    (12,10)node{$40$}    (12,8)node{$30$}
    (12,6)node{$20$}    (12,4)node{$10$}    (14,12)node{$51$}    (14,10)node{$41$}
    (14,8)node{$31$}    (14,6)node{$21$}    (16,10)node{$42$}    (16,8)node{$32$};
\end{tikzpicture}
\caption{AR quiver of $\Dwq(A_n)/[2]$, for $n=5$}\label{fig:AR An}
\end{figure}
\begin{figure}[ht]\centering
\begin{tikzpicture}[xscale=.7,yscale=.5]
\path (-1,0) coordinate (t1) node[below]{$V_0$};
\path (1,0) coordinate (t2) node[below]{$V_1$};
\path (3,3) coordinate (t3) node[right]{$V_2$};
\path (3.5,7) coordinate (t4)node[right]{$V_3$};
\path (-1,10) coordinate (t5)node[above]{$V_4$};
\path (-4,5) coordinate (t6)node[left]{$V_5$};
\foreach \j in {1,...,6}{
\foreach \k in {\j,...,6}{\draw[blue!50,ultra thick,opacity=.3](t\j)to(t\k);}}
\draw[fill=orange!10](t1)to(t3)to(t5);
\draw[fill=orange!10](t2)to(t4)to(t6);
\draw[fill=white](t1)to(t2)to(t3) (t4)to(t2)to(t3) (t5)to(t4)to(t3)
    (t4)to(t5)to(t6) (t1)to(t5)to(t6) (t1)to(t2)to(t6);

\foreach \j/\i in {1/2,2/3,3/4,4/5,5/6,6/1}
{\draw[->,>=stealth,very thick,blue] (t\j) to (t\i);}

\foreach \j/\i in {1/3,2/4,3/5,4/6,5/1,6/2}
{\draw[orange,very thick] (t\j) to (t\i);}
\foreach \j in {1,2,3,4,5,6}
{\draw (t\j) node {$\bullet$};}
\draw[blue!50]
    ($(t1)!.5!(t2)$) node[below] {$z_1$}
    ($(t3)!.5!(t2)$) node[below right] {$z_2$}
    ($(t3)!.5!(t4)$) node[right] {$z_3$}
    ($(t5)!.5!(t4)$) node[above] {$z_4$}
    ($(t5)!.5!(t6)$) node[above left] {$z_5$}
    ($(t6)!.5!(t1)$) node[below left] {$z_0$};
\end{tikzpicture}
\caption{The $6$-gon of type $A_5$ with level-2 diagonal-gon shaded}\label{fig:tt}
\end{figure}

Recall that $\Sth(A_n)$ denotes the moduli space of stable $(n+1)$-gons on $\CC$ up to translation.
The following proposition is a $\CC^*$-covering version of \cite[Prop.~3.6]{Q18}.
But the proof is essentially the same.

\begin{proposition}\label{pp:A}
There is a natural isomorphism
\[
    Z_h\colon\ToSt(A_n)/[2]\to\Sth(A_n)
\]
sending a total stability condition $\sigma$ to
the far-end $(n+1)$-gon.
\end{proposition}

More precisely, we have $z_j=Z(P^{j-1}_{1})$ for $j\in \ZZ_{n+1}$
for the edges of the far-end $(n+1)$-gon and
\[
    Z(P^{j}_i)=V_jV_{j+i},\quad 1\le i\le n, 0\le j\le n-i.
\]

A study of total stability for the module category $\h$ of type A is in \cite{Ki},
which can be deduced from the description above
after fixing the heart of a total stability condition to be $\h$.

\section{$\ToSt$ of type $D_n$}\label{sec:D}
\begin{setup}\label{setup:D}
Set $m=n-1$ and then $h_Q=2m$ for $Q=D_n$.
We choose an orientation as follows:
\[
\begin{tikzpicture}[scale=1,xscale=1,ar/.style={->,>=stealth,thick}]
\foreach\j in {0,2,4} {\draw[ar](\j+.5,0)to++(-1,0);}
\draw[ar](6.5,.5)to++(-1,-.3);
\draw[ar](6.5,-.5)to++(-1,+.3);
\draw(-1,0)node{$1$}(1,0)node{$2$}(3,0)node{$\cdots$}
    (5,0)node{\small{$n-2$}}(7,.5)node{\small{$m$}}(7,-.5)node{\small{$n$}};
\end{tikzpicture}
\]
Moreover, we fix 1 as our favorite far-end vertex for $n=4$
and it is the unique far-end vertex for $n\ge5$.
Label the indecomposable objects in the $\Dwq(D_n)$ still by \eqref{eq:Pji}.
Thus the $i^{\mathrm{th}}$-$\tau$-orbit consisting of objects $\{P^j_i\}_{j\in\ZZ}$ is the one corresponding to vertex $i$.

As $\tshift^{2m}=[2]$, the indecomposables in the root category
are $\{P^{j}_{i}\mid 1\leq i\leq n,  j\in\ZZ_{2m}\}$.
For instance, see the AR quiver $\AR\Dwq(D_{5})/[2]$ in the low picture of Figure~\ref{fig:AR.D5},
where $j (j+i)$ denotes $P^{j}_{i}$, $j\;\parity(j)$ denotes $P^j_{m}$
and $j\;\parity(j+1)$ denotes $P^j_{n}$.
Here
$$
    \parity(x)=\mathrm{sign}(-1)^x\colon\ZZ\to \{\pm\}
$$
is the parity function.
\end{setup}

\subsection{Symmetric doubly punctured $h$-gons}\label{sec:D0}
\begin{definition}
A \emph{symmetric} $h$-gon $\hgon$ is an $h$-gon (cf. notations in Section~\ref{sec:summ}) such that
$h$ is even and
\begin{gather}\label{eq:z h}
    z_{j+h/2}=-z_{j}.
\end{gather}
Let $O$ be its geometric center.
It is \emph{doubly punctured} if there is a pair of punctures $B_\pm$
such that $O$ is the middle point of them.
\end{definition}

Given a central charge $Z\colon K\Dwq(D_n)\to\CC$, we construct a symmetric doubly punctured $2(n-1)$-gon $(\hgon_Z,B_\pm)$ as follows:
\begin{gather}\label{eq:ZtoP}
\begin{cases}\displaystyle
  V_j=-\big(\sum_{t=j}^{j+m-1}Z(P^{t}_{1}) \big)/2,\\
  \displaystyle B_+=\big( Z(P_{m})-Z(P_{n} ) \big)/2\\
  \displaystyle B_-=\big( Z(P_{n})-Z(P_{m} ) \big)/2.
\end{cases}
\end{gather}
So $z_j=V_{j-1}V_j=Z(P^{j-1}_{1})$ noticing $Z(P^{j+m-1}_{1})=-Z(P^{j-1}_{1})$.
More precisely,
\begin{itemize}
\item $\hgon=\hgon_Z$ is the far-end $2(n-1)$-gon,
which is indeed symmetric since
\[\tshift^m=[1]\]
holds on the far-end $\tau$-orbit $\{P^{t}_{1}\mid t\in\ZZ_{2m}\}$
that implies \eqref{eq:z h}.
\item Moreover, the geometric center of $\hgon$ is at origin as
\begin{gather}\label{eq:=-=}
\begin{cases}
    V_j=-V_{j+m}, \quad j\in\ZZ_{2m} \\
    B_+=-B_-.
\end{cases}
\end{gather}
\end{itemize}

\begin{definition}
An \emph{$h$-gon of type $D_n$} is a symmetric doubly punctured $2(n-1)$-gon.
It is \emph{stable} (of type $D_n$) if it is positively convex and the punctures are inside the level-$(n-2)$ diagonal-gon.
Denote by $\Sth(D_n)$ the moduli space of stable $2(n-1)$-gons of type $D_n$ up to translation,
which has complex dimension $n$.
\end{definition}

Next we show that total stability implies the stability of $2(n-1)$-gons of type $D_n$.

\begin{proposition}\label{pp:Dn0}
If $Z$ is the central charge of some $\sigma\in\ToSt(D_n)$,
then the symmetric doubly punctured $2(n-1)$-gon $(\hgon,B_\pm)$ defined above is a stable $2(n-1)$-gon of type $D_n$.
We call it the far-end stable $2(n-1)$-gon with respect to $\sigma$.
\end{proposition}
\begin{proof}
We need to show that the total stability of $\sigma$ implies the stability of the corresponding far-end $h$-gon $\hgon$ together with punctures $B_\pm$.

Firstly, the positive convexity of $\hgon$ follows from Proposition~\ref{pp:sthgon}.

Secondly, we will check that $B_\pm$ is inside the level-$(m-1)$ diagonal-gon.
Namely, for any $j\in\ZZ_{2m}$, it suffices to show that $B_{\pm}$
is on the left hand side of $V_{j+1}V_{j+m}$:
\[\begin{tikzpicture}[scale=2.5,rotate=15]
  \draw[teal,very thick](-70:1)node[below]{$V_{j}$}
    (-20:1)node[right]{$V_{j+1}$}to(180-70:1)node[left]{$V_{j+m}$};
  \draw[dashed](-70:1)to(180-70:1) (-20:1)to(180-20:1)
        (0,0)node[above left]{$_O$}\nn;
  \draw[orange](-70:1)to(-160:.2)node[below]{$B_\pm$}\nn to(160:1);
  \draw[orange](180-70:1)to(180-160:.2)node[above]{$B_\mp$}\nn to(180+160:1);
  \draw[teal,very thick](-70:1)to(160:1)node[left]{$V_{j+m+1}$};
\end{tikzpicture}\]
Note that we have AR triangles ($i\in\{m,n\}$)
\[
    P^{j}_{i}\to P^{j+1}_{n-2}\to P^{j+1}_{i} \to P^{j}_{i}[1]
\]
with central charges
\[
    V_{j}B_{\pm}=B_{\mp}V_{j+m},\quad V_{j+1}V_{j+m},\quad
    V_{j+1}B_{\mp}\quad\text{and}\quad -V_{j}B_{\pm}
\]
Here $\pm=\parity(j)$ if $i=m$ and
$\pm=\parity(j+1)$ if $i=n$, to be more precise.
Since all $P^?_?$ are stable objects, their phases are increasing.
Therefore we have
\[
    \arg B_{\mp}V_{j+m} < \arg V_{j+1}V_{j+m} <
    \arg V_{j+1}B_{\mp} < \arg B_{\mp}V_{j+m}+\pi,
\]
where $\arg$ takes values in $[\pi\cdot\phi_\sigma (P^{j}_{i}), \pi\cdot\phi_\sigma (P^{j}_{i})+2\pi)$.
This completes the proof.
\end{proof}

\subsection{Geometric model for root category of type $D_n$}
Suppose that we have a $2(n-1)$-gon $(\hgon,B_\pm)$ of type $D_n$ (with vertices $V_{j}$ and punctures $B_\pm$).
Up to translation, we may assume its geometric center is at the origin,
i.e. \eqref{eq:=-=} holds. Then we have
\begin{gather}\label{eq:sym}
\begin{cases}
  V_{j}V_k=V_{k+m}V_{j+m},\\
  V_{j}B_+=B_-V_{j+m}.
\end{cases}
\end{gather}

\begin{figure}[hb]\centering
\makebox[\textwidth][c]{
\begin{tikzpicture}[scale=2, rotate=-30, xscale=1,arrow/.style={->,>=stealth,thick}]
\clip[rotate=30] (-2,-2) rectangle (2,2);
\path (0,0) coordinate (O)
      (40:1.12) coordinate (v1)
      (78:1.3) coordinate (v2)
      (123:1.2) coordinate (v3)
      (172:1.7) coordinate (v4);
\draw[cyan, ultra thick]
      ($(O)-.5*(v1)-.5*(v2)-.5*(v3)-.5*(v4)$) coordinate (V0) coordinate (V8)  coordinate (W5)
      ($(V0)+(v1)$) coordinate (V1) coordinate (W6)
      ($(V1)+(v2)$) coordinate (V2) coordinate (W7)
      ($(V2)+(v3)$) coordinate (V3) coordinate (W8) coordinate (W0)
      ($(V3)+(v4)$) coordinate (V4) coordinate (W1);
\path ($(O)-(V1)$) coordinate (V5) coordinate (W2)
      ($(O)-(V2)$) coordinate (V6) coordinate (W3)
      ($(O)-(V3)$) coordinate (V7) coordinate (W4);
\draw[cyan, ultra thick, fill=cyan!20](V0)
    \foreach \j in {1,...,8} {to(V\j)};
\foreach \j in {0,...,7}
    {\draw[white,fill=white](V\j)to($3*(V\j)$)to($3*(W\j)$)to(W\j);}
\foreach \j in {0,...,7} {\draw[Green, very thick](V\j)to(W\j);
\draw[thin,gray,dotted](V\j)to($(O)-(V\j)$)  ($1.1*(V\j)$)node[blue]{$V_\j$};}
\draw[cyan, ultra thick](V0)
    \foreach \j in {1,...,8} {to(V\j)\nn };
\draw[dotted](125:.3) coordinate
    (p) \nn ($1.5*(p)$) node{\footnotesize{$B_+$}}to[thin,gray,dotted]
    ($(O)-(p)$)\nn ($(O)-1.5*(p)$) node{\footnotesize{$B_-$}};
\end{tikzpicture}
\begin{tikzpicture}[scale=2, rotate=-30, xscale=1,arrow/.style={->,>=stealth,thick}]
\clip[rotate=30] (-2,-2) rectangle (2,2);
\path (0,0) coordinate (O)
      (40:1.12) coordinate (v1)
      (78:1.3) coordinate (v2)
      (123:1.2) coordinate (v3)
      (172:1.7) coordinate (v4);
\draw[cyan, ultra thick]
      ($(O)-.5*(v1)-.5*(v2)-.5*(v3)-.5*(v4)$) coordinate (V0) coordinate (V8)  coordinate (W5)
      ($(V0)+(v1)$) coordinate (V1) coordinate (W6)
      ($(V1)+(v2)$) coordinate (V2) coordinate (W7)
      ($(V2)+(v3)$) coordinate (V3) coordinate (W8)coordinate (W0)
      ($(V3)+(v4)$) coordinate (V4) coordinate (W1);
\path ($(O)-(V1)$) coordinate (V5) coordinate (W2)
      ($(O)-(V2)$) coordinate (V6) coordinate (W3)
      ($(O)-(V3)$) coordinate (V7) coordinate (W4);
\draw (125:.3) coordinate (p);
\foreach \j in {1,3,5,7} {
  \draw[orange,ultra thick](V\j)to(p);
  \draw[orange,dashed,thick](V\j)to($(O)-(p)$);}
\foreach \j in {2,4,6,8} {
  \draw[orange,dashed,thick](V\j)to(p);
  \draw[orange,ultra thick](V\j)to($(O)-(p)$);}
\foreach \j in {1,...,8} {
\draw[dashed,gray,dashed,thin](V\j)($(O)-(V\j)$)  ($1.1*(V\j)$)node[blue]{$V_\j$};}
\draw[cyan,opacity=0.5](V0)
    \foreach \j in {1,...,8} {to(V\j)\nn };
\draw(p) \nn ($1.7*(p)$) node{\footnotesize{$B_+$}}
    ($(O)-(p)$)\nn ($(O)-1.7*(p)$) node{\footnotesize{$B_-$}};
\foreach \j in {1,...,8} {\draw[Green,thin,opacity=0.5](V\j)to(W\j);}
\draw[blue!50,ultra thick]
    (V1)to(V3)to(V5)to(V7)to(V1)(V2)to(V4)to(V6)to(V8)to(V2);
\end{tikzpicture}}
\makebox[\textwidth][c]{
\def\quan{circle(.6)}
\begin{tikzpicture}[scale=.7, rotate=0,ar/.style={->,>=stealth,thick,teal},
ra/.style={<-,>=stealth,thick,teal}]
\tikzstyle{every node}=[font=\small,circle]
\clip[rotate=0] (1,-3) rectangle (21,5);
\draw[dashed,thin,cyan](-2,4)to(21,4);
\draw[dashed,thin,blue!50](-2,2)to(21,2);
\draw[dashed,thin,Green](-2,0)to(21,0);
\draw[dashed,thin,orange](-2,-2.3)to(21,-2.3)
    (-2,-.7)to(21,-.7);
\foreach \j in {0,...,5}{
\begin{scope}[shift={(\j*4,0)}]
\draw[white,draw=white]
    (0,0) node[fill=white] (a\j) {$00$}
    (4,0) node[fill=white] (A\j) {$00$}
    (2,-2-.3) node[fill=white] (p\j) {$00$}
    (2,-1+.3) node[fill=white] (m\j) {$00$}
    (2,2) node[fill=white] (b\j) {$00$}
    (0,4) node[fill=white] (c\j) {$00$}
    (4,4) node[fill=white] (C\j) {$00$};
\end{scope}}
\foreach \j in {0,...,5}{
  \draw[Green,fill=Green,opacity=.32](a\j)\quan;
  \draw[blue!50,fill=blue!50,opacity=.32](b\j)\quan;
  \draw[cyan,fill=cyan,opacity=.32](c\j)\quan;
  \draw[orange,fill=orange,opacity=.25](p\j)\quan(m\j)\quan;}
\draw[gray,dashed,fill=Green!15](a5)\quan;
\draw[gray,dashed,fill=blue!8](b0)\quan;
\draw[gray,dashed,fill=cyan!15](c5)\quan;
\draw[gray,dashed,fill=orange!15](p0)\quan(m0)\quan;
\foreach \j in {0,...,5}{
\begin{scope}[shift={(\j*4,0)}]
\draw[ar](a\j)edge(p\j)edge(m\j)edge(b\j);
\draw[ra](A\j)edge(p\j)edge(m\j)edge(b\j);
\draw[ar](c\j)to(b\j);\draw[ra](C\j)to(b\j);
\end{scope}}
\draw
  (a1)node{$\frac{72}{63}$}
  (a2)node{$\frac{03}{74}$}
  (a3)node{$\frac{14}{05}$}
  (a4)node{$\frac{25}{16}$}
  (b1)node{$\frac{02}{64}$}
  (b2)node{$\frac{13}{75}$}
  (b3)node{$\frac{24}{06}$}
  (b4)node{$\frac{35}{17}$}
  (c1)node{$\frac{01}{54}$}
  (c2)node{$\frac{12}{65}$}
  (c3)node{$\frac{23}{76}$}
  (c4)node{$\frac{34}{07}$}
  (p1)node{$\frac{7+}{-3}$}
  (p2)node{$\frac{0-}{+4}$}
  (p3)node{$\frac{1+}{-5}$}
  (p4)node{$\frac{2-}{+6}$}
  (m1)node{$\frac{7-}{+3}$}
  (m2)node{$\frac{0+}{-4}$}
  (m3)node{$\frac{1-}{+5}$}
  (m4)node{$\frac{2+}{-6}$};
\draw[gray!80]
  (a5)node{$\frac{36}{27}$}
  (c5)node{$\frac{45}{10}$}
  (b0)node{$\frac{71}{53}$}
  (p0)node{$\frac{6-}{+2}$}
  (m0)node{$\frac{6+}{-2}$};
\end{tikzpicture}
}
\caption{The $8$-gon of type $D_5$ and AR quiver of $\Dwq(D_{5})$.}
\label{fig:AR.D5}
\end{figure}

\begin{theorem}\label{thm:geo.Dn}
A $2(n-1)$-gon of type $D_n$ is a geometric model for the root category $\Dwq(D_n)/[2]$ in the sense that by setting
\begin{gather}\label{eq:Z Dn}
\begin{cases}
    Z(P^{j}_{i})=V_{j}V_{j+i}=V_{j+i+m}V_{j+m},& j\in\ZZ_{2m},1\le i\le m-1,\\
    Z(P^{j}_{m})=V_{j}B_{\parity(j)}=B_{\parity(j+1)}V_{j+m},& j\in\ZZ_{2m},\\
    Z(P^{j}_{n})=V_{j}B_{\parity(j+1)}=B_{\parity(j)}V_{j+m},& j\in\ZZ_{2m},
\end{cases}
\end{gather}
we obtain a central charge $Z\colon K\Dwq(D_n)\to\CC$.
\end{theorem}
\begin{proof}
To show that $Z$ is a group homomorphism,
we only need to check all the mesh relations in the root category still hold after applying $Z$.
To start with, notice that
\begin{gather}\label{eq:sum mn}
    Z(P^{j}_{m})+Z(P^{j}_{n})=V_{j}V_{j+m}
\end{gather}
by \eqref{eq:sym}.
The rest of the proof is just a direct checking:
\begin{enumerate}
 \item At the $\tau$-orbit of $i=1$
 the mesh relation has the form
 \[ [P^{j}_{1}]+[P^{j+1}_{1}]=[P^{j}_{2}] \]
 and indeed we have
 \[ V_{j}V_{j+1}+V_{j+1}V_{j+2}=V_{j}V_{j+2}.\]
 \item At the $\tau$-orbit of $1<i<n-2$
 the mesh relation has the form
 \[ [P^{j}_{i}]+[P^{j+1}_{i}]=[P^{j+1}_{i-1}]+[P^{j}_{i+1}] \]
 and indeed we have
 \[ V_{j}V_{j+i}+V_{j+1}V_{j+i+1}
    =V_{j+1}V_{j+i}+V_{j}V_{j+i+1}.\]
 \item At the $\tau$-orbit of $i=n-2$
 the mesh relation has the form
 \[ [P^{j}_{n-2}]+[P^{j+1}_{n-2}]=[P^{j+1}_{n-3}]+[P^{j}_{m}]+[P^{j}_{n}] \]
 and indeed we have, using \eqref{eq:sum mn},
 \[ V_{j}V_{j+n-2}+V_{j+1}V_{j+n-1}
    =V_{j+1}V_{j+n-2}+V_{j}V_{j+n-1}.\]
 \item At the $\tau$-orbit of $i\in\{m,n\}$,
 the mesh relation has the form
 \[ [P^{j}_{i}]+[P^{j+1}_{i}]=[P^{j+1}_{n-2}] \]
 and indeed we have, using \eqref{eq:sym},
 \[ V_{j}B_{\pm}+V_{j+1}B_{\mp}=
  B_{\mp}V_{j+m}+V_{j+1}B_{\mp}=
  V_{j+1}V_{j+m}.\]
  Note that here the sign $\pm$ depends on $i\in\{m,n\}$ and the parity function $\rho$.\qedhere
\end{enumerate}
\end{proof}

\begin{example}[Type $Q=D_5$]
In Figure~\ref{fig:AR.D5}, we have
\begin{itemize}
\item The objects drawn in blue/violet/green circle, in the $\tau$-orbit of $\AR\Dwq(Q)$, correspond to
length-1/2/3 diagonals drawn in blue/violet/green respectively as in the upper $8$-gons.
\item The objects drawn in orange circle, in the upper/lower $\tau$-orbit of $\AR\Dwq(Q)$, correspond to
solid/dashed orange line segments respectively as in the right upper $8$-gon.
\end{itemize}

Note that all the (direct checking) calculations in the proof of Theorem~\ref{thm:geo.Dn}
can be easily read off from Figure~\ref{fig:AR.D5} for case $n=5$.
\end{example}

\subsection{From stable $h$-gon to total stability conditions of type $D_n$}\label{sec:recon:Dn}
Given a stable $2(n-1)$-gon $\hgon$ of type $D_n$, we can construct a total stability condition $\sigma=(Z,\sli)$ with $Z$ defined as in Theorem~\ref{thm:geo.Dn}. 

\begin{construction}\label{con:sli}
Let us construct a slicing $\sli$ as follows.

{\bf Step~1}: Assign a real number $\phi(M)$ for each object $M$ in the far-end $\tau$-orbit of $\AR \Dwq(D_n)$:
\begin{itemize}
\item Let
\[
    \phi(P^{}_{1}[k])\colon= \arg Z(P^{}_{1}) /\pi +k,
\]
for any $k\in\ZZ$, where $\arg$ takes values in $[0,2\pi)$.
\item For $1\leq j\leq m-1$ and $k\in\ZZ$, let
\[
    \phi(P^{j}_{1}[k])\colon= \arg Z(P^{j}_{1}) /\pi +k,
\]
where $\arg$ takes values in $[\arg Z(P_{1}),Z(P_{1})+2\pi)$.
\end{itemize}
So we have the monotonicity and periodicity:
\begin{gather}\label{eq:ineq1}
\begin{cases}
    \phi(P^{j}_{1})<\phi(P^{j+1}_{1}),\\
    \phi(P^{j+m}_{1})=\phi(P^{j}_{1}[1])=\phi(P^{j}_{1})+1,
\end{cases}\quad \forall j\in\ZZ.
\end{gather}

{\bf Step~2}: For any $2\leq i\leq n-2<m$ and $j\in\ZZ$,
we have
\[
    Z(P^{j}_{i})=V_{j} V_{j+i}=\sum_{t=j}^{j+i-1}V_tV_{t+1}=Z(P^{j}_{1})+\cdots Z(P^{j+i-1}_{1}).
\]
Note that by Step~1, the positive convexity of $\hgon$ gives
\[
    \phi(P^{j}_{1})<\phi(P^{j+1}_{1})<\cdots<\phi(P^{j+i-1}_{1})< \phi(P^{j+m}_{1})=\phi(P^{j}_{1}[1])=\phi(P^{j}_{1})+1
\]
Then let
\begin{gather}\label{eq:arg j}
    \phi(P^{j}_{i}[k])\colon= \arg Z(P^{j}_{i}) /\pi +k,
\end{gather}
where $\arg$ takes values in $[\pi\cdot\phi(P^{j}_{1}),\pi\cdot\phi(P^{j}_{1})+2\pi)$.
It is straightforward to check that
the monotonicity and periodicity of \eqref{eq:ineq1} are inherited:
\[\begin{cases}
    \phi(P^{j}_{i})<\phi(P^{j+1}_{i}),\\
    \phi(P^{j+m}_{i})=\phi(P^{j}_{i}[1])=\phi(P^{j}_{i})+1,
\end{cases}\quad \forall j\in\ZZ.\]
Note that, locally,
the additive subcategory of $\AR\Dwq(D_n)$ generated by
\begin{gather}\label{eq:sub Aj}
\begin{cases}
    P^{j}_{1}, \ldots ,P^{j+i-2}_{1}, P^{j+i-1}_{1}\\
    P^{j}_{2}, \ldots , P^{j+i-2}_{2}\\
    \cdots\\
    P^{j}_{i}
\end{cases}
\end{gather}
is isomorphic to $\mod \k A_{i}$ for an $A_i$ quiver with straight orientation.
Inductively, the positive convexity implies that
\begin{gather}\label{eq:ineq2}
\begin{cases}
    \phi(P^{s}_{t})<\phi(P^{s+1}_{t-1})<\phi(P^{s+1}_{t})
        & 1<t\le i, j\le s<j+i-t\\
    \phi(P^{s}_{t})<\phi(P^{s}_{t+1})<\phi(P^{s+1}_{t})
        & 1\le t< i, j\le s<j+i-t.
\end{cases}
\end{gather}

{\bf Step~3}:
For $i\in\{m,n\}$ and $j\in\ZZ$, let
\begin{gather}\label{eq:ineq3}
    \phi(P^{j}_{i})\colon= \arg Z(P^{j}_{i}) /\pi,
\end{gather}
where $\arg$ takes value in $[\pi\cdot\phi(P^{j}_{n-2}),\pi\cdot\phi(P^{j}_{n-2})+2\pi)$.
Note that
$$\phi(P^{j+m}_{n-2})=\phi(P^{j}_{n-2}[1])=\phi(P^{j}_{n-2})+1,$$
and
\[\begin{cases}
    \phi(P^{j+m}_{m})=\phi(P^{j}_{m}[1])=\phi(P^{j}_{m})+1,\\
    \phi(P^{j+m}_{n})=\phi(P^{j}_{n}[1])=\phi(P^{j}_{n})+1,
\end{cases}\]
if $m$ is even and
\[\begin{cases}
    \phi(P^{j+m}_{m})=\phi(P^{j}_{n}[1])=\phi(P^{j}_{n})+1,\\
    \phi(P^{j+m}_{n})=\phi(P^{j}_{m}[1])=\phi(P^{j}_{m})+1,
\end{cases}\]
if $m$ is odd.
This completes the assigning $\phi$.

{\bf Step~4}:
Define
\begin{gather}\label{eq:sli=}
    \sli(\varphi)=  \Add \Big(
        \bigoplus_{ \substack{P\in\Ind\Dwq(D_n) \\ \phi(P)=\varphi}} P\Big)
\end{gather}
\end{construction}

\begin{proposition}\label{pp:Dn}
$\sigma=(Z,\sli)$ defined as above is a total stability condition on $\Dwq(D_n)$.
\end{proposition}
\begin{proof}
By construction, we already have
\begin{itemize}
\item all indecomposable objects are in some $\sli(\phi)$ for $\phi\in\RR$ and
\item $\sli$ is compatible with the central charge $Z$ as well as the shift $[1]$.
\end{itemize}

Thus, what is left to show is $\phi(M)<\phi(L)$
whenever there is a non-zero map from $M$ to $L$
for any indecomposable objects $M,L$.
This amounts to check that, for any arrow
$M\to L$ in $\AR\Dwq(D_n)$, we have $\phi(M)<\phi(L)$.
There are three cases:
\begin{itemize}
\item If the arrow is between the far-end $\tau$-orbit
(corresponding to vertex 1) and the double-trivalent $\tau$-orbit (corresponding to vertex $n-2$),
then it is a type $A_n$ issue,  cf. \eqref{eq:sub Aj},
and \eqref{eq:ineq2} gives the inequality.
\item If the arrow is from the double-trivalent $\tau$-orbit
to the boundary $\tau$-orbits (corresponding to vertex $m$ or $n$),
then \eqref{eq:ineq3} implies the inequality.
\item If the arrow is from the boundary $\tau$-orbits (corresponding to vertex $m$ or $n$) to the double-trivalent $\tau$-orbit,
then we need to examine the triangle
\[\xymatrix@R=1pc@C=1pc{
& P^{j+1}_{n-3}\ar[dr]\\
P^{j}_{n-2}\ar[ur]\ar[r]\ar[dr]& P^{j}_{m}\ar[r] & P^{j+1}_{n-2}\\
& P^{j}_{n} \ar[ur] }\]
that corresponds to
\[\begin{tikzpicture}[scale=2,rotate=10]
  \draw[teal,very thick](-70:1)node[below]{$V_{j}$}to
    (70:1)node[right]{$V_{j+m-1}$}
    (-20:1)node[right]{$V_{j+1}$}to
    (180-70:1)node[left]{$V_{j+m}$};
  \draw[dashed](-70:1)to(180-70:1) (-20:1)to(180-20:1)
        (0,0)node[above]{$O$}\nn;
  \draw[orange](-70:1)to(-160:.2)node[below]{$B_?$}\nn to(110:1);
  \draw[teal,very thick](-70:1)to(160:1)node[left]{$V_{j+m+1}$}   (110:1)to(-110:1)node[below]{$V_{j-1}$};
\end{tikzpicture}\]
Then, from the stability condition of $(\hgon,B_\pm)$, we know that
$B_\pm$ is bounded by solid lines in the above picture.
By \eqref{eq:arg j} and \eqref{eq:ineq3},
we deduce that both $\phi(P^{j}_{i})$, $i\in\{m,n\}$, and $\phi(P^{j+1}_{n-2})$ are in
$[\phi(P^{j}_{n-2}),\phi(P^{j}_{n-2})+2)$.
Thus we have
\[\arraycolsep=1.4pt\def\arraystretch{2.2}
\begin{array}{rl}
    \phi(P^{j}_{n-2})&=\displaystyle\frac{\arg V_{j}V_{j+m-1}}{\pi} \\
    &<\displaystyle\frac{\arg V_{j}B_\pm}{\pi} =\phi(P^{j}_{i}) \\
    &<\displaystyle\frac{\arg V_{j+1}V_{j+m}}{\pi}=\phi(P^{j+1}_{n-2})
\end{array}\]
where $\pm$ in the second row depends on $i\in\{m,n\}$ and the parity function $\rho$. Note that $\arg$ takes values in
$[\pi\cdot\phi(P^{j}_{n-2}),\pi\cdot\phi(P^{j}_{n-2})+2\pi)$.
\end{itemize}
In all, the first case uses positive convexity of $\hgon$
and the last two cases uses the extra condition of stability of $(\hgon,B_\pm)$.
\end{proof}

\begin{remark}
Alternatively, one can define $\sli$ as follows:
\begin{itemize}
\item For any $1\leq i\leq n$, let
\[
    \phi(P_i)=\arg Z(P_i) /\pi,
\]
where $\arg$ takes value in $[\arg Z(P_1),Z(P_1)+2\pi)$.
\item For fixed $1\leq i\leq n$, let
\[
    \phi(P^{j}_{i})=\arg Z(P^{j}_{i}) /\pi
\]
for $1\leq j\leq m$,
where $\arg$ takes values in $[\arg Z(P_i),Z(P_i)+2\pi)$.
\item Define $\sli$ as in \eqref{eq:sli=}.
\end{itemize}
But the simplicity of this construction will trade off more complexity of the proof above.
\end{remark}

\begin{theorem}\label{thm:D}
There is a natural isomorphism
\[ Z_h\colon\ToSt(D_n)/[2]\to\Sth(D_n),\]
sending a total stability condition $\sigma$ to the far-end stable $2(n-1)$-gon.
\end{theorem}
\begin{proof}
Note that in Construction~\ref{con:sli} one can shift all slices in a slicing simultaneously by $[2]$.
Then the theorem follows from combining Propositions~\ref{pp:Dn0} and~\ref{pp:Dn}.
\end{proof}

\subsection{Example: $D_4$ with three far-end stable $6$-gons}\label{sec:D4}
\begin{figure}[h]\centering
\includegraphics[width=5in]{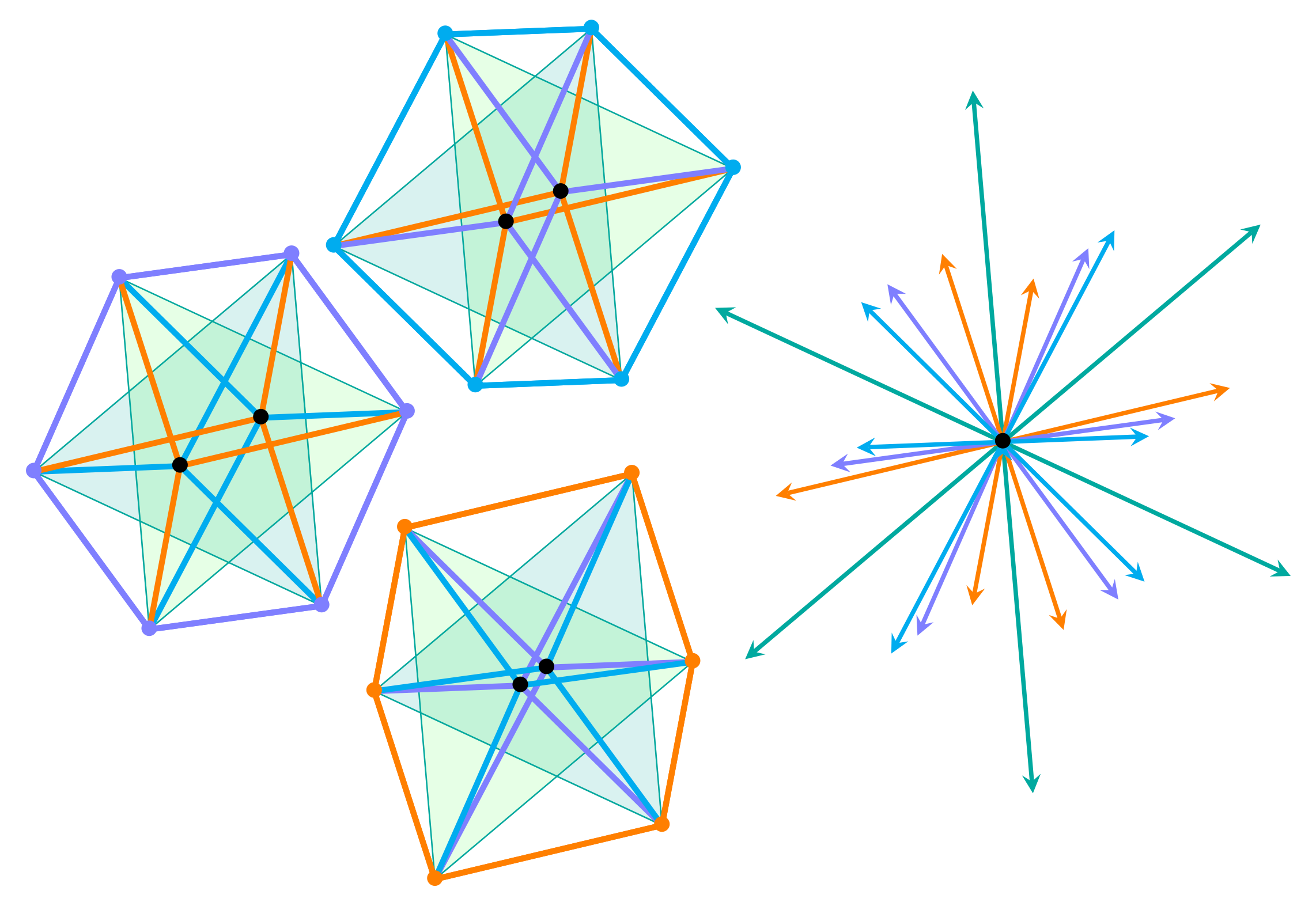}
\caption{Three far-end $6$-gons of type $D_4$ for a given central charge}
\label{fig:D4}
\end{figure}
As mentioned above a couple of times,
there are three choices of far-end vertices/$\tau$-orbits/$6$-gons.
The isomorphism in Theorem~\ref{thm:D} depends on such a choice.
Here, we give an example, as shown in Figure~\ref{fig:D4},
of how the three (far-end) stable $6$-gons look like (3 hexagons on the left),
with respect to a fixed $\sigma\in\ToSt(D_4)$ (whose central charge is on the right).

\section{Exceptional type: Preview}\label{sec:exp}
For an exceptional quiver $Q=E_n$, where $n\in\{6,7,8\}$,
we fix the following orientation
\begin{gather}\label{eq:En}
\begin{tikzpicture}[scale=.5,xscale=1,ar/.style={->,>=stealth,thick}]
\draw(0,0)node(v1){$1$}(2,0)node(v2){$2$}(4,0)node(v3){$4$}(6,0)node(v5){}(8,0)node(v6){}(10,0)node(v7){$n$}
    (4,2)node(v4){$3$}
    (7,0)node{$\cdots$};
\draw[ar](v2)edge(v1) (v3)edge(v2) (v3)edge(v4) (v6)edge(v7)  (v3)edge(v5);
\end{tikzpicture}
\end{gather}
and still label the indecomposable objects in the $\Dwq(E_n)$ by \eqref{eq:Pji}.
The objects in the boundary $\tau$-orbits play the key role in constructing geometric model and we label them as
\begin{gather}\label{eq:137}
\begin{cases}
  \text{$1^{\mathrm{st}}$-$\tau$-orbit: Denote $C_{j}=P^{j-1}_{1}$.}&\text{
  This is the mid-end $\tau$-orbit. }\\
  \text{$3^{\mathrm{rd}}$-$\tau$-orbit: Denote $M_{j}=P^{j-1}_{3}$.}&\text{
  This is the near-end $\tau$-orbit. }\\
  \text{$n^{\mathrm{th}}$-$\tau$-orbit: Denote $B_{j}=P^{j-1}_{n}$.}&\text{
  This is the far-end $\tau$-orbit. }\\
\end{cases}
\end{gather}
In Figure~\ref{fig:AR E6}, \ref{fig:E7} and \ref{fig:AR E8}, the objects in the mid-end/near-end/far-end orbits will be drawn in green/yellow/blue in the AR-quiver, respectively.

Moreover, we will use the following family of triangles
\begin{gather}\label{eq:common1}
 B_{j} \to C_{j+2} \to B_{j+n-2} \to B_{j}[1]
\end{gather}
and octahedral diagrams
\begin{gather}\label{eq:common2}
  \xymatrix{
    & B_j \ar[d]\ar@{=}[r] & B_j \ar[d]\\
    C_{j} \ar@{=}[d] \ar[r] & M_{j+1} \ar[r] \ar[d]
        & B_{j+n-3} \ar[d] \ar[r] & C_{j}[1] \ar@{=}[d]\\
    C_{j} \ar[r] & C_{j+3} \ar[d] \ar[r] &
        L_j \ar[d]\ar[r] & C_{j}[1]\\
    & B_j[1] \ar@{=}[r]& B_j[1]\\
  }
\end{gather}
in $\DQ$, for
\begin{gather}L_j=\begin{cases}
    M_{j+4} & n=6,\\
    C_{j+6} & n=7,\\
    B_{j+10} & n=8.
\end{cases}\end{gather}

\begin{remark}\label{rem:sum}
We will construct a geometric model $\hgon$ for $\DQ/[2]$ of a type $E_n$ quiver $Q$ consisting of the following data:
\begin{itemize}
\item $\hgon$ is an $h_Q$ gon whose edges are given by some central charge of the indecomposable objects in the far-end $\tau$-orbit of $\Dwq(Q)$.
\item its ice/fire core is an $h_Q/2$ gon whose edges are given by the (same) central charge of half of the indecomposable objects in the mid-end $\tau$-orbit of $\Dwq(Q)$.
\end{itemize}
We say such an $h_Q$-gon is stable if it is positively convex and both the ice and fire cores are inside the length-$(n-3)$ diagonal-gon of $\hgon$
(see Definition~\ref{def:sth E6}).
Moreover, the key of realizing the root category is knowing the central charges of objects in the boundary $\tau$-orbits,
i.e. we will prove:
\begin{gather}\label{eq:common}
 \begin{cases}
   Z(P^{j-1}_{1})=Z(C_j)=W_{j-1}W_{j+1}\\
   Z(P^{j-1}_{3})=Z(M_j)=V_{j-2}W_{j}=W_{j+1}V_{j+n-4},\\
   Z(P^{j-1}_{n})=Z(B_j)=V_{j-1}V_{j}.
 \end{cases}
\end{gather}
in the following sections.
We will also show that the edges of the ice/fire cores are given by
\begin{gather}\label{eq:w_j}
    w_{j+2}\colon=z_{j}+z_{j+n-2},\quad j\in\ZZ_{h}
\end{gather}
for $h=h_{Q}$.

Furthermore, the usual types (A and D) can also be thought as a degeneration of the model above, in the sense that:
\begin{itemize}
  \item for $Q=D_n$, the ice/fire core shrinks to the puncture $B_\pm$.
  \item for $Q=A_n$, the ice/fire core vanishes.
\end{itemize}
\end{remark}

\section{$\ToSt$ of exceptional type $E_6$}\label{sec:E6}
\begin{figure}[hb]\centering
\begin{tikzpicture}[scale=2, rotate=15, arrow/.style={->,>=stealth,thick}]
\clip[rotate=-15] (-3,-.3)rectangle(3,5);
\path (0,0) coordinate (o)
      (1.15,0) coordinate (v1)
      (30+2:1.23) coordinate (v2)
      (60-4:1.2) coordinate (v3)
      (90-1:1.2) coordinate (v4)
      (120+4:1.15) coordinate (v5)
      (150+3:1.15) coordinate (v6)
      ($(o)-(v1)-(v5)$) coordinate (v9)
      ($(o)-(v2)-(v6)$) coordinate (v10)
      ($(v4)+(v10)-(v1)$) coordinate (v7)
      ($(o)-(v3)-(v7)$) coordinate (v11)
      ($(v5)+(v11)-(v2)$) coordinate (v8)
      ($(v3)+(v9)-(v6)$) coordinate (v12);
\path ($(v1)+(v2)+(v3)$) coordinate (v123)
      ($(v2)+(v4)+(v3)$) coordinate (v234)
      ($(v4)+(v5)+(v3)$) coordinate (v345)
      ($(v4)+(v5)+(v6)$) coordinate (v456)
      ($(v5)+(v6)+(v7)$) coordinate (v567)
      ($(v8)+(v6)+(v7)$) coordinate (v678);
\draw[font=\scriptsize](o)node[below]{$V_0$}++(v1)node[below]{$V_1$}
    ++(v2)node[below right]{$V_2$}++(v3)node[right]{$V_3$}
    ++(v4)node[right]{$V_4$}++(v5)node[above]{$V_5$}
    ++(v6)node[above]{$V_6$}++(v7)node[above left]{$V_7$}
    ++(v8)node[left]{$V_8$}++(v9)node[left]{$V_9$}
    ++(v10)node[left]{$V_{10}$}++(v11)node[below]{$V_{11}$};
\draw[white,fill=\fire,opacity=.11]
    ($(o)-(v9)$)to++(v6)to++(v4)to++(v2)to++(v12)to++(v10)to++(v8);
\draw[white,fill=\ice,opacity=.14]
    ($(v4)$)to++(v5)to++(v3)to++(v1)to++(v11)to++(v9)to++(v7);
\draw[ultra thick]($(v1)-(v10)$)edge[\sun]+(v7)edge[orange]+($(o)-(v9)$);
\draw[ultra thick]($(v1)+(v345)$)edge[\sun]+(v11)edge[orange]+($(o)-(v1)$);
\draw[ultra thick]($(v4)+(v5)$)edge[\sun]+(v3)edge[orange]+($(o)-(v5)$);
\draw[ultra thick]($(o)-(v9)$)edge[cyan]+(v6)edge[blue]+($(o)-(v8)$);
\draw[ultra thick]($(v234)$)edge[cyan]+(v10)edge[blue]+($(o)-(v12)$);
\draw[ultra thick]($(v456)+(v1)$)edge[cyan]+(v2)edge[blue]+($(o)-(v4)$);
\draw[font=\Large](o)
    to node[below,orange]{$z_1$}++(v1)
    to node[below right,cyan]{$z_2$}++(v2)
    to node[right,teal]{$z_3$}++(v3)
    to node[right,blue]{$z_4$}++(v4)
    to node[above right,orange]{$z_5$}++(v5)
    to node[above,cyan]{$z_6$}++(v6)
    to node[above,teal]{$z_7$}++(v7)
    to node[above left,blue]{$z_8$}++(v8)
    to node[left,orange]{$z_9$}++(v9)
    to node[left,cyan]{$z_{10}$}++(v10)
    to node[below left,teal]{$z_{11}$}++(v11)
    to node[below,blue]{$z_{0}$}++(v12);
\draw[cyan, ultra thick,fill=\fcyan]
    (v1)\nn \foreach\j in {2,6,10}{-- ++(v\j)\nn}
    ($(v123)+(v456)+(v10)$)\nn \foreach\j in {2,6,10}{-- ++(v\j)\nn}
    ($(v1)+(v567)$)\nn \foreach\j in {2,6,10}{-- ++(v\j)\nn};
\draw[blue, ultra thick,fill=\fblue]
    ($(v123)$)\nn \foreach\j in {4,8,12}{-- ++(v\j)\nn}
    ($(v123)+(v567)$)\nn \foreach\j in {4,8,12}{-- ++(v\j)\nn}
    (o)\nn \foreach\j in {4,8,12}{-- ++(v\j)\nn};
\draw[\sun, ultra thick,fill=\fgreen]
    ($(v123)-(v3)$)\nn \foreach\j in {3,7,11}{-- ++(v\j)\nn}
    ($(v123)+(v456)-(v3)$)\nn \foreach\j in {3,7,11}{-- ++(v\j)\nn}
    ($(o)-(v12)$)\nn \foreach\j in {3,7,11}{-- ++(v\j)\nn};
\draw[orange, ultra thick,fill=\forange]
    ($(v234)$)\nn \foreach\j in {1,5,9}{-- ++(v\j)\nn}
    ($(o)+(v456)$)\nn \foreach\j in {1,5,9}{-- ++(v\j)\nn}
    ($(o)$)\nn \foreach\j in {1,5,9}{-- ++(v\j)\nn};
%
\draw[gray](o)\foreach \j in {1,...,12}{++(v\j)\nn};
\draw[font=\tiny]
      (v4)node[above,blue]{$W_{0}$}
    ++(v5)node[right,cyan]{$\;\;W_{10}$}
    ++(v3)node[below right,blue]{$W_8$}
    ++(v1)node[right,cyan]{$W_6$}
    ++(v11)node[left,blue]{$W_4$}
    ++(v9)++(-.05,.05)node[above left,cyan]{$W_2$};
\draw[font=\tiny]
      ($(v1)+(v5)$)+(0,-.1)node[left,orange]{$W_1$}
    ++(v6)node[above left,\sun]{$W_{11}$}
    ++(v4)+(0,.01)node[above,orange]{$W_{9}$}
    ++(v2)node[left,\sun]{$W_7$}
    ++(v12)+(-.1,0)node[below,orange]{$W_5$}
    ++(v10)+(0,.01)node[below,\sun]{$W_3$};
\end{tikzpicture}

\begin{tikzpicture}\draw
    (0,0)node[draw=red]{$z_j+z_{j+4}+z_{j+8}=0,\quad \forall j\in\ZZ_{12}$}
    (0,-.8)node[draw=red]{$z_{j}-z_{j-3}+z_{j-6}-z_{j-9}=0,\quad \forall j\in\ZZ_{12}$} (0,-1.2)node{};
\end{tikzpicture}

\begin{tikzpicture}[xscale=.8,yscale=.8,rotate=-45, font=\tiny,
arrow/.style={->,>=stealth,thick}]
\clip[rotate=-45](-.8,2)rectangle(6.3,19);
\foreach \j in {-3,...,8}{
\foreach \k in {0,...,8}{
    \draw[teal,->,>=stealth,very thick]
        (2*\k,2*\j+.5)to++(0,1);
    \draw[teal,->,>=stealth,very thick]
        (2*\k+.5,2*\j)to++(1,0);
}}
\draw[rotate=-45,white,fill=white](0,2)rectangle(-2,19)(8,2)rectangle(5.5,19);
\foreach \j in {-3,...,8}{
\foreach \k in {0,...,8}{
    \draw[white](2*\k,2*\j) node[circle,draw=teal] (x\k\j) {$000$};
}}
\foreach \j in {1,...,6}{
    \draw[yellow!10](2*\j+2+1+.005,2*\j-1-.005) node[circle,draw=teal,fill=yellow!15](m\j) {$000$};
    \draw(m\j)node{$M_{\j}$};
}
\foreach \j in {0,...,5}{
    \draw[teal,->,>=stealth,very thick] (2*\j+4+1+.005+.35,2*\j+1-.005+.35)to (2*\j+4+2-.35,2*\j+2-.35);
    \draw[teal,->,>=stealth,very thick] (2*\j+4+.35,2*\j+.35) to (2*\j+4+1+.005-.35,2*\j+1-.005-.35);
}
\foreach \j/\k in {0/3,1/4,2/5,3/6,4/7,5/8}{
    \draw[](2*\j+6,2*\j+2) node {$\j/\k$};
}
\draw[](-2+6,-2+2) node {$11/2$};
\foreach \j/\k in {0/2,1/3,2/4,3/5,4/6,5/7}{
    \draw[](2*\j+6,2*\j) node {$\j/\k$};
}
\foreach \j/\k in {8/6,9/7,10/8,11/9,12/10}{
    \draw[](2*\j+6-18,2*\j-14) node {$\j/\k$};
}
\foreach \j/\k in {13/11}{
    \draw[](2*\j+6-18,2*\j-14) node {$1/\k$};
}
\foreach \j/\k in {0/1,1/2,2/3,3/4,4/5,5/6}{
    \draw[](2*\j+6,2*\j-2) node[cyan!15,circle,draw=teal,fill=cyan!15] {$000$} node {$\j/\k$};
    \draw[](2*\j+6-.5,2*\j-2+.5) node {$B_\k$};
}
\foreach \j/\k/\l in {7/6/1,8/7/2,9/8/3,10/9/4,11/10/5,12/11/6}{
    \draw[](2*\j-12,2*\j-12) node[green!15,circle,draw=teal,fill=green!15] {$000$} node {$\j/\k$};
    \draw[](2*\j-12+.5,2*\j-12-.5) node {$C_\l$};
}
\end{tikzpicture}
\caption{The $12$-gon of type $E_6$ and AR quiver of $\Dwq(E_6)$}
\label{fig:AR E6}
\end{figure}

\subsection{The $h$-gon of type $E_6$}\label{sec:hgon E6}
\begin{definition}\label{def:hgon E6}
An \emph{$h$-gon $\hgon$ of type $E_6$} is a $12$-gon satisfying \eqref{eq:E6-rel},
i.e. the equations in Figure~\ref{fig:AR E6}.
\end{definition}

A direct calculation shows that:
\begin{lemma}\label{lem:E6}
The set of 7 equations \eqref{eq:E6-rel} has rank 6.
So the space of $12$-gons of type $E_6$ has complex dimension $12-6=6$.
\end{lemma}
By \eqref{eq:E6-rel}, we have the following hexagon relations.
\begin{gather}\label{eq:E6-8}
    z_{j}+z_{j+1}+z_{j+2}+z_{j+6}+z_{j+7}+z_{j+8}=0,
    \quad\forall j\in\ZZ_{12}.
\end{gather}

\begin{figure}[hb]\centering
\begin{tikzpicture}[scale=1.2, rotate=15, arrow/.style={->,>=stealth,thick}]
\clip[rotate=-15] (-3,-.2)rectangle(7.5,5);
\begin{scope}[shift={(-14:5)}]
\path (0,0) coordinate (o)
      (1.15,0) coordinate (v1)
      (30+2:1.23) coordinate (v2)
      (60-4:1.2) coordinate (v3)
      (90-1:1.2) coordinate (v4)
      (120+4:1.15) coordinate (v5)
      (150+3:1.15) coordinate (v6)
      ($(o)-(v1)-(v5)$) coordinate (v9)
      ($(o)-(v2)-(v6)$) coordinate (v10)
      ($(v4)+(v10)-(v1)$) coordinate (v7)
      ($(o)-(v3)-(v7)$) coordinate (v11)
      ($(v5)+(v11)-(v2)$) coordinate (v8)
      ($(v3)+(v9)-(v6)$) coordinate (v12);
\path ($(v1)+(v2)+(v3)$) coordinate (v123)
      ($(v2)+(v4)+(v3)$) coordinate (v234)
      ($(v4)+(v5)+(v3)$) coordinate (v345)
      ($(v4)+(v5)+(v6)$) coordinate (v456)
      ($(v5)+(v6)+(v7)$) coordinate (v567)
      ($(v8)+(v6)+(v7)$) coordinate (v678);
\draw[white,fill=\fires]
    (o)\foreach \j in {1,...,12}{to++(v\j)};
\draw[white,fill=\fires](o)to++($(o)-(v9)$)to++(v6)to++($(o)-(v3)$)--cycle;
\draw[white,fill=\fires]
    ($(v123)+(v4)-(v1)$)to++($(o)-(v9)$)to++(v6)to++($(o)-(v3)$)--cycle;
\draw[white,fill=\fires](v1)to++(v2)to++($(o)-(v11)$)to++(v8)--cycle;
\draw[white,fill=\fires]($(v456)+(v1)$)to++(v2)to++($(o)-(v11)$)to++(v8)--cycle;
\draw[white,fill=\fires](v123)to++(v4)to++($(o)-(v1)$)to++(v10)--cycle;
\draw[white,fill=\fires]($(v6)-(v9)$)to++(v4)to++($(o)-(v1)$)to++(v10)--cycle;
\draw[white,fill=\fire,opacity=.6]
    ($(o)-(v9)$)to++(v6)to++(v4)to++(v2)to++(v12)to++(v10)to++(v8);
\draw[very thick]($(o)-(v9)$)edge[cyan]+(v6)edge[blue]+($(o)-(v8)$);
\draw[very thick]($(v234)$)edge[cyan]+(v10)edge[blue]+($(o)-(v12)$);
\draw[very thick]($(v456)+(v1)$)edge[cyan]+(v2)edge[blue]+($(o)-(v4)$);
\draw[cyan, very thick,fill=\fcyan]
    (v1)to++(v2)++(v345)to++(v6) (v456)to++(v10);
\draw[blue, very thick,fill=\fblue]
    (o)to($(o)-(v12)$) (v123)to++(v4)++(v567)to+(v8);
\draw[Green, very thick,fill=\fgreen]
    ($(v123)-(v3)$)\nn \foreach\j in {3,7,11}{-- ++(v\j)\nn}
    ($(v123)+(v456)-(v3)$)\nn \foreach\j in {3,7,11}{-- ++(v\j)\nn}
    ($(o)-(v12)$)\nn \foreach\j in {3,7,11}{-- ++(v\j)\nn};
\draw[orange, very thick,fill=\forange]
    ($(v234)$)\nn \foreach\j in {1,5,9}{-- ++(v\j)\nn}
    ($(o)+(v456)$)\nn \foreach\j in {1,5,9}{-- ++(v\j)\nn}
    ($(o)$)\nn \foreach\j in {1,5,9}{-- ++(v\j)\nn};
\draw[red,very thick,dashed]
    (o)to++(v123)++(v4)to++(v567)++(v8)to++($(o)-(v345)$)
    (v456)to++(v123)(v1)to++(v567)   ($(v1)+(v2)$)to++(v345);
\draw[gray](o)\foreach \j in {1,...,12}{ ++(v\j)\nn};
\end{scope}
\path (0,0) coordinate (o)
      (1.15,0) coordinate (v1)
      (30+2:1.23) coordinate (v2)
      (60-4:1.2) coordinate (v3)
      (90-1:1.2) coordinate (v4)
      (120+4:1.15) coordinate (v5)
      (150+3:1.15) coordinate (v6)
      ($(o)-(v1)-(v5)$) coordinate (v9)
      ($(o)-(v2)-(v6)$) coordinate (v10)
      ($(v4)+(v10)-(v1)$) coordinate (v7)
      ($(o)-(v3)-(v7)$) coordinate (v11)
      ($(v5)+(v11)-(v2)$) coordinate (v8)
      ($(v3)+(v9)-(v6)$) coordinate (v12);
\path ($(v1)+(v2)+(v3)$) coordinate (v123)
      ($(v2)+(v4)+(v3)$) coordinate (v234)
      ($(v4)+(v5)+(v3)$) coordinate (v345)
      ($(v4)+(v5)+(v6)$) coordinate (v456)
      ($(v5)+(v6)+(v7)$) coordinate (v567)
      ($(v8)+(v6)+(v7)$) coordinate (v678);
\draw[white,fill=\ices]
    (o)\foreach \j in {1,...,12}{to++(v\j)};
\draw[white,fill=\ice,opacity=0.6]
    (v4)to++(v5)to++(v3)to++(v1)to++(v11)to++(v9)to++(v7);
\draw[very thick]($(v1)-(v10)$)edge[\sun]+(v7)edge[orange]+($(o)-(v9)$);
\draw[very thick]($(v1)+(v345)$)edge[\sun]+(v11)edge[orange]+($(o)-(v1)$);
\draw[very thick]($(v4)+(v5)$)edge[\sun]+(v3)edge[orange]+($(o)-(v5)$);
\draw[cyan, very thick,fill=\fcyan]
    (v1)\nn \foreach\j in {2,6,10}{-- ++(v\j)\nn}
    ($(v123)+(v456)+(v10)$)\nn \foreach\j in {2,6,10}{-- ++(v\j)\nn}
    ($(v1)+(v567)$)\nn \foreach\j in {2,6,10}{-- ++(v\j)\nn};
\draw[blue, very thick,fill=\fblue]
    ($(v123)$)\nn \foreach\j in {4,8,12}{-- ++(v\j)\nn}
    ($(v123)+(v567)$)\nn \foreach\j in {4,8,12}{-- ++(v\j)\nn}
    (o)\nn \foreach\j in {4,8,12}{-- ++(v\j)\nn};
\draw[Green, very thick,fill=\fgreen]
    (v123)to($(v123)-(v3)$)
    ($(v123)+(v456)$)to++(v7)
    ($(v1)+(v567)$)to+(v11);
\draw[orange, very thick,fill=\forange]
    (o)to(v1)++(v234)to++(v5) (v456)to++($(o)-(v9)$);
\draw[blue,very thick,dashed]
    (v1)to++(v234)++(v5)to++(v678)++(v9)to++($(o)-(v456)$)
    ($(v1)+(v2)$)to++(v678) (v123)to++(v456) ($(v456)+(v10)$)to++(v234);
\draw[gray](o)\foreach \j in {1,...,12}{ ++(v\j)\nn};
\end{tikzpicture}
\caption{Configurations Ice and Fire}
\label{fig:E6+-}
\end{figure}


\begin{construction}\label{con:E6}
Using the triangle relations in \eqref{eq:E6-rel},
we can draw triangles $$\TT_j\colon=V_{j-1}V_jW_{j}    \quad\forall j\in\ZZ_{12}$$ with edges
\[
    V_{j-1}V_{j}=z_j,\quad V_jW_{j}=z_{j+4},\quad
    W_{j}V_{j-1}=z_{j+8}.
\]
Note that
we have 4 sets of 3 parallel triangles, see Figure~\ref{fig:AR E6},
\[
    \{\TT_{j+4k}\mid k=0,1,2\}\quad \text{for $j\in\ZZ_4$},
\]
drawn in orange/blue/green/violet respectively.
Moreover, the square relations in \eqref{eq:E6-rel} correspond to the squares
$$\surf_j\colon=V_{j-1}V_jW_{j+1}W_{j-1}    \quad\forall j\in\ZZ_{12}$$
with edges
\[
    V_{j-1}V_{j}=z_j,\quad V_jW_{j+1}=-z_{j-3},\quad
    W_{j+1}W_{j-1}=z_{j-6},\quad W_{j-1}V_{j-1}=-z_{j-9}.
\]
Note that
we have 3 sets of 4 parallel/anti-parallel squares
\[
    \{\surf_{j+3k}\mid k=0,1,2,3\}\quad \text{for $j\in\ZZ_3$},
\]
see the squares in Figure~\ref{fig:E6+-}.

Set $w_j$ as \eqref{eq:w_j} for $n=6$ and $h=12$.
By the triangle relation in \eqref{eq:E6-rel}, we have $w_{j+2}=-z_{j+8}.$

It is better to have two copies of $\hgon$
(called \emph{ice} and \emph{fire}, cf. Figure~\ref{fig:E6+-}):
\begin{description}
\item[$\ihgon$] the one with 6 triangles $\{\TT_{2j}\mid j\in\ZZ_6\}$ and 6 squares $\{\surf_{2j+1}\mid j\in\ZZ_6\}$.
Then we have a hexagon \emph{ice core} $\icore$ of $\hgon$ with vertices $W_{2j}$ and edges $$w_{2j+1}=W_{2j}W_{2j+2}.$$
\item[$\fhgon$] the one with 6 triangles $\{\TT_{2j+1}\mid j\in\ZZ_6\}$ and 6 squares $\{\surf_{2j}\mid j\in\ZZ_6\}$.
Then we have a hexagon \emph{fire core} $\fcore$ of $\hgon$ with vertices $W_{2j+1}$ and edges $$w_{2j}=W_{2j-1}W_{2j+1}.$$
\end{description}
\end{construction}

\begin{remark}
An interesting feature is that both ice and fire configurations induce planar tiling patterns, see Figure~\ref{fig:A Tiling} in Appendix~\ref{app}.
In other words, a $12$-gon of type $E_6$ is equivalent to `A Tiling of Ice and Fire' as shown there.
\end{remark}

\subsection{Geometric model for the root category of type $E_6$}
\begin{setup}\label{setup:E6}
Take the orientation in \eqref{eq:En} with $n=6$.
Note that 1/3/6 are the mid-end/near-end/far-end vertices in this case.
The lower picture in Figure~\ref{fig:AR E6} is part of $\AR\Dwq(E_6)$ and
the labelling is given by
\begin{description}
  \item[$1^{\mathrm{st}}$-$\tau$-orbit] $j+6/j+5$ denotes $C_{j}\colon=P^{j-1}_{1}$ (drawn in green).
  \item[$2^{\mathrm{nd}}$-$\tau$-orbit] $j+7/j+5$ denotes $P^{j-1}_{2}$.
  \item[$3^{\mathrm{rd}}$-$\tau$-orbit] Let $M_{j}\colon=P^{j-1}_{3}$ (drawn in yellow).
  \item[$4^{\mathrm{th}}$-$\tau$-orbit] $j-1/j+2$ denotes $P^{j-1}_{4}$.
  \item[$5^{\mathrm{th}}$-$\tau$-orbit] $j-1/j+1$ denotes $P^{j-1}_{5}$.
  \item[$6^{\mathrm{th}}$-$\tau$-orbit] $j-1/j$ denotes $B_{j}\colon=P^{j-1}_{6}$ (drawn in blue).
\end{description}
for $1\leq j\leq 6$.
\end{setup}

Given any central charge $Z\colon K\Dwq(E_6)\to\CC$.
Let $\hgon_Z$ be the far-end $12$-gon of $\Dwq(E_6)$
with edges $z_{j}=Z(B_{j})$ for $1\le j\le 12$.
\begin{lemma}\label{lem:E6-1}
The set $\{ [B_{j}] \mid 1\leq j\leq 6\}$ spans $K \Dwq(E_6)$ (and in fact is a basis).
Moreover, $\hgon_Z$ is a $12$-gon of type $E_6$
and its ice/fire core is formed by (the odd/even) half of
$$\{ Z(C_{j}) \mid 1\leq j\leq 12\}$$
for $Z(C_{j})=w_{j}$.
\end{lemma}
A coincidence in this case is $w_j=-z_{j+6}$ as mentioned above.
\begin{proof}
The first statement can be checked directly, e.g. via dimension vectors.

Now  we will show that all $z_{j}$ satisfy \eqref{eq:E6-rel}.
On one hand, by \eqref{eq:common1} for $n=6$ we have
\[
    [B_{j}] + [B_{j+4}] = [C_{j+2}]
\]
and thus $Z(C_{j+2})=z_j+z_{j+4}=w_{j+2}$ by \eqref{eq:w_j}.
Due to $B_{j+8}=C_{j+2}[1]$, we have
\[
    [B_{j}] + [B_{j+4}] + [B_{j+8}] = 0.
\]
Thus the triangle relations in \eqref{eq:E6-rel} holds.
On the other hand, triangles in the second row/column of \eqref{eq:common2} for $n=6$ implies that
\[
    [B_{j}]+[C_{j+3}]=[M_{j+1}]
    =[C_{j}]+[B_{j+3}].
\]
Noticing that $C_{j+6}=B_{j}[1]$ and $B_{j+3}=B_{j-9}[2]$, the above equation becomes
\[
    [B_{j}]-[B_{j-3}]=-[B_{j-6}]+[B_{j-9}],
\]
again using.
Thus the square relations in \eqref{eq:E6-rel} holds.

In particular, we see that the edges of the ice/fire core are also $w_j$'s.
\end{proof}

\begin{figure}[hb]\centering
\begin{tikzpicture}[scale=1.8, rotate=15, arrow/.style={->,>=stealth,thick},font=\tiny]
\clip[rotate=-15] (-3,-.3)rectangle(3,5);
\path (0,0) coordinate (o)
      (1.15,0) coordinate (v1)
      (30+2:1.23) coordinate (v2)
      (60-4:1.2) coordinate (v3)
      (90-1:1.2) coordinate (v4)
      (120+4:1.15) coordinate (v5)
      (150+3:1.15) coordinate (v6)
      ($(o)-(v1)-(v5)$) coordinate (v9)
      ($(o)-(v2)-(v6)$) coordinate (v10)
      ($(v4)+(v10)-(v1)$) coordinate (v7)
      ($(o)-(v3)-(v7)$) coordinate (v11)
      ($(v5)+(v11)-(v2)$) coordinate (v8)
      ($(v3)+(v9)-(v6)$) coordinate (v12);
\path ($(v1)+(v2)+(v3)$) coordinate (v123)
      ($(v2)+(v4)+(v3)$) coordinate (v234)
      ($(v4)+(v5)+(v3)$) coordinate (v345)
      ($(v4)+(v5)+(v6)$) coordinate (v456)
      ($(v5)+(v6)+(v7)$) coordinate (v567)
      ($(v8)+(v6)+(v7)$) coordinate (v678);
\draw[](o)node[below]{$V_0$}++(v1)node[below]{$V_1$}
    ++(v2)node[below right]{$V_2$}++(v3)node[right]{$V_3$}
    ++(v4)node[right]{$V_4$}++(v5)node[above]{$V_5$}
    ++(v6)node[above]{$V_6$}++(v7)node[above left]{$V_7$}
    ++(v8)node[left]{$V_8$}++(v9)node[left]{$V_9$}
    ++(v10)node[left]{$V_{10}$}++(v11)node[below]{$V_{11}$};
\draw[fill=yellow,opacity=.5]
    (o)to++($(v1)-(v10)$)to++($(v3)-(v6)$)--cycle;
\draw[fill=yellow,opacity=.1]
    (v456)to++($(v3)-(v6)$)to++($(v1)-(v10)$)--cycle;
\draw[very thick]
    (o)to++($(v1)-(v10)$)
    (v456)++(v123)to++($(v10)-(v1)$)
    (v1)++(v234)to++($(v10)-(v1)$)
    (v456)++(v10)to++($(v1)-(v10)$);
%
%
\draw[thick]($(v1)-(v10)$)edge[\sun]+(v7)edge[orange]+($(o)-(v9)$);
\draw[thick]($(v1)+(v345)$)edge[\sun]+(v11)edge[orange]+($(o)-(v1)$);
\draw[thick]($(v4)+(v5)$)edge[\sun]+(v3)edge[orange]+($(o)-(v5)$);
\draw[thick]($(o)-(v9)$)edge[cyan]+(v6)edge[blue]+($(o)-(v8)$);
\draw[thick]($(v234)$)edge[cyan]+(v10)edge[blue]+($(o)-(v12)$);
\draw[thick]($(v456)+(v1)$)edge[cyan]+(v2)edge[blue]+($(o)-(v4)$);
\draw[cyan, thick]
    (v1)\nn \foreach\j in {2,6,10}{-- ++(v\j)\nn}
    ($(v123)+(v456)+(v10)$)\nn \foreach\j in {2,6,10}{-- ++(v\j)\nn}
    ($(v1)+(v567)$)\nn \foreach\j in {2,6,10}{-- ++(v\j)\nn};
\draw[blue, thick]
    ($(v123)$)\nn \foreach\j in {4,8,12}{-- ++(v\j)\nn}
    ($(v123)+(v567)$)\nn \foreach\j in {4,8,12}{-- ++(v\j)\nn}
    (o)\nn \foreach\j in {4,8,12}{-- ++(v\j)\nn};
\draw[Green, thick]
    ($(v123)-(v3)$)\nn \foreach\j in {3,7,11}{-- ++(v\j)\nn}
    ($(v123)+(v456)-(v3)$)\nn \foreach\j in {3,7,11}{-- ++(v\j)\nn}
    ($(o)-(v12)$)\nn \foreach\j in {3,7,11}{-- ++(v\j)\nn};
\draw[orange, thick]
    ($(v234)$)\nn \foreach\j in {1,5,9}{-- ++(v\j)\nn}
    ($(o)+(v456)$)\nn \foreach\j in {1,5,9}{-- ++(v\j)\nn}
    ($(o)$)\nn \foreach\j in {1,5,9}{-- ++(v\j)\nn};
\draw[gray, thick]
    (o)to++(v123) (v456)to++(v123) ;
\draw[gray](o)\foreach \j in {1,...,12}{++(v\j)\nn};
\draw[font=\tiny]
      (v4)node[above,blue]{$W_{0}$}
    ++(v5)node[right,cyan]{$\;\;W_{10}$}
    ++(v3)node[below right,blue]{$W_8$}
    ++(v1)node[right,cyan]{$W_6$}
    ++(v11)node[left,blue]{$W_4$}
    ++(v9)++(-.05,.05)node[above left,cyan]{$W_2$};
\draw[font=\tiny]
      ($(v1)+(v5)$)+(0,-.1)node[left,orange]{$W_1$}
    ++(v6)node[above left,\sun]{$W_{11}$}
    ++(v4)+(0,.01)node[above,orange]{$W_{9}$}
    ++(v2)node[above right,\sun]{$W_7$}
    ++(v12)+(-.1,0)node[below,orange]{$W_5$}
    ++(v10)+(0,.01)node[below,\sun]{$W_3$};
\end{tikzpicture}
\caption{The triangles/squares for $E_6$}
\label{fig:E6 tech}
\end{figure}

In this subsection, we will construct a central charge from a $12$-gon $\hgon$ of type $E_6$.
Recall that $\hgon$ has edges $z_{j}=V_{j-1}V_{j}$, $j\in\ZZ_{12}$,
as shown in the upper picture of Figure~\ref{fig:AR E6}.
Moreover, we have the associated ice/fire configuration/core according to Construction~\ref{con:E6}.

\begin{theorem}\label{thm:geo.E6}
A $12$-gon $\hgon$ of type $E_6$ is a geometric model for the root category $\Dwq(E_6)/[2]$
in the sense that by setting (for $n=6$)
\begin{gather}\label{eq:commonx}
 \begin{cases}
   Z(P_1^j)=W_{j}W_{j+2}\\
   Z(P_2^j)=V_{j-1}W_{j+2}=W_{j+1}V_{j+n-3}\\
   Z(P_3^j)=V_{j-1}W_{j+1}=W_{j+2}V_{j+n-3},\\
   Z(P_{i}^j)=V_{j}V_{j+1+(n-i)},& 4\le i\le n,
 \end{cases}
\end{gather}
we obtain a central charge $Z\colon K\Dwq(E_6)\to\CC$.
\end{theorem}
In particular, \eqref{eq:common} holds (for mid-end/near-end/far-end $\tau$-orbits).
Also, central charges of objects in the mid-end $\tau$-orbits are given by the edges of ice/fire core of $\hgon$
and central charges of objects in the $i^{\mathrm{th}}$-$\tau$-orbits are given by length-$(7-i)$ diagonals of $\hgon$ for $4\le i\le n$.

\begin{remark}
In fact, there are several ways to realize the indecomposable objects in the $\tau$-orbits. For instance,
for any object $E\in\Ind\mod\k E_6$ that is not in the $3^{\mathrm{rd}}$-$\tau$-orbit, it admits a labeling $j/k$ as in Set-up~\ref{setup:E6}, which is also the labelling corresponding to its image under the central charge, namely
\[
    Z(E)=V_{j}V_k \quad \text{and}\quad Z(E[1])=V_kV_{j};
\]
for each $M_j$ in the $3^{\mathrm{rd}}$-$\tau$-orbits, it can be realized in two more ways:
\begin{gather}\label{eq:4ways}
    W_{j+h/2}V_{j-2+h/2}=V_{j-2}W_{j}=W_{j+1}V_{j+n-4}=V_{j+n-4+h/2}W_{j+1+h/2}
\end{gather}
for $h=12$.
See the thick black line segments in Figure~\ref{fig:E6 tech} for $j=2$.

In type $E_6$, the far-end and mid-end $\tau$-orbits are shift of each other which enable us to realize the objects in many ways. However, in types $E_7$ and $E_8$, this is not the case.
Therefore, we write the statement in the above theorem of the form that can be generalized to types $E_7$ and $E_8$.
\end{remark}

\begin{proof}[Proof of Theorem~\ref{thm:geo.E6}]
The proof follows the same way as in the proof of Theorem~\ref{thm:geo.Dn}.
One has to show that the central charge $Z$ preserves the mesh relations of $\AR\Dwq(E_6)/[2]$.
Here we only point out the reason that $Z(M_{j})$ can be realized as in \eqref{eq:4ways} is due to
the triangles in the second row/column of \eqref{eq:common2} and Construction~\ref{con:E6}.
The remaining calculations are left to the reader.
\end{proof}

\subsection{Stability of $12$-gon for $E_6$}\label{sec:sth E6}
Suppose that $\hgon$ is a positively convex $12$-gon of type $E_6$.
We have the following observations.
\begin{itemize}
\item Its ice/fire core is convex.
This follows from the fact that the edges of ice/fire core are
exactly half of the edges of the $12$-gon induced by the mid-end $\tau$-orbit.
\item Any length-3 diagonal $V_{j-1}V_{j+2}$ is a long diagonal of a narrow hexagon
\[\mathbf{H}_j=V_{j-1}V_jV_{j+1}V_{j+2}W_{j+2}W_{j}\]
that corresponds to \eqref{eq:E6-8}.
Note that $\mathbf{H}_j$ and $\mathbf{H}_{j+6}$ are related by a translation $V_{j+2}V_{j+5}$.
Hence they inherit the convexity of $\hgon$.
\item See the left picture in Figure~\ref{fig:E6+-}, the configuration ice contains 6 narrow hexagons
$\{\mathbf{H}_{2j}\}$ and the corresponding length-3 diagonals (blue dashed lines in the left picture) bound another hexagon, we call the \emph{ice-boundary},
which contains the ice core automatically due to the convexity of the narrow hexagons.
\item Dually, see the right picture in Figure~\ref{fig:E6+-}, the configuration fire contains 6 narrow hexagons
$\{\mathbf{H}_{2j+1}\}$ and the corresponding length-3 diagonals (red dashed lines in the right picture) bound another hexagon, we call the \emph{fire-boundary},
which also contains the fire core automatically.
\end{itemize}

\begin{definition}\label{def:sth E6}
A $12$-gon $\hgon$ of type $E_6$ is \emph{stable} if it is positively convex and
both its ice and fire cores are inside the level-3 diagonal-gon.
\end{definition}
Denote by $\Sth(E_6)$ the moduli space of stable $12$-gon of type $E_6$ up to translation.
By Lemma~\ref{lem:E6}, the complex dimension of $\Sth(E_6)$ is $6$.

Using the above observations, one can simplify the stability condition of $12$-gons a bit as in the lemma below.
\begin{lemma}\label{lem:sth equiv}
A positively convex $12$-gon $\hgon$ of type $E_6$ is stable if and only if
its ice/fire core is inside its fire/ice boundary respectively.
\end{lemma}
\begin{proof}
As we know, the positive convexity of $\hgon$ implies the convexity of the narrow hexagons. Then, the ice/fire core is automatically inside ice/fire boundary respectively.
Moreover, the level-3 diagonal-gon is the intersection of the ice and fire boundaries.
Thus, the stability condition of $\hgon$ can be simplified to the condition stated in the lemma,
provided the positive convexity of $\hgon$.
\end{proof}

Now we proceed to show that total stability of $\sigma$ deduces
stability of its far-end $12$-gon.

\begin{proposition}\label{pp:sth E6}
If $\sigma=(Z,\sli)\in\ToSt(E_6)$,
then its far-end $12$-gon is a stable $12$-gon of type $E_6$.
\end{proposition}
\begin{proof}
By Proposition~\ref{pp:sthgon}, the far-end $12$-gon $\hgon$ is positively convex.
Since the equation \eqref{eq:E6-rel} holds due to Lemma~\ref{lem:E6-1},
we can follow Construction~\ref{con:E6} to obtain its ice/fire core with vertices $W_j, j\in\ZZ_{12}$ via triangles (cf. Figure~\ref{fig:AR E6}).
Then we only need to show that each $W_{j}$ are bounded by length-3-diagonals of $\hgon$.
By Lemma~\ref{lem:sth equiv},
it suffices to check that $W_{2j}$ is bounded by the fire boundary and $W_{2j+1}$ is bounded by the ice boundary.
By the observation of the narrow hexagons,
it is equivalent to prove that each $W_{j+1}$ is on the left side of the diagonal $V_{j-1}V_{j+2}$.

Consider the triangle $V_{j-1} W_{j+1} V_{j+2}$ with edges
\[
    V_{j-1}W_{j+1} =Z(P^{j}_{3}),\quad
    W_{j+1} V_{j+2}=Z(P^{j-1}_{3}),\quad
    V_{j-1} V_{j+2}=Z(P^{j-1}_{4}).
\]
These edges correspond to the central charges of the terms in the following AR triangle:
\begin{gather}\label{eq:dis.e6}
 P^{j-1}_{3} \to P^{j-1}_{4} \to P^{j}_{3} \to P^{j-1}_{3}[1].
\end{gather}
See the yellow triangle (and its mirror in lighter yellow) in Figure~\ref{fig:E6 tech} for $j=1$.
And the total stability implies
\[
    \phi_\sigma(P^{j-1}_{3})<
    \phi_\sigma(P^{j-1}_{4})<
    \phi_\sigma(P^{j}_{3})<
    \phi_\sigma(P^{j-1}_{3})+1
\]
By taking $\arg$ in $[\pi\cdot \phi_\sigma(P^{j-1}_{3}), \pi\cdot \phi_\sigma(P^{j-1}_{3})+2\pi)$,
we conclude that each $W_{j+1}$ is indeed on the left side of $Z(P^{j-1}_{4})=V_{j-1}V_{j+2}$.
\end{proof}

\subsection{From stable $12$-gon to total stability conditions for $E_6$}\label{sec:recon:E6}
\begin{theorem}\label{thm:E6}
There is a natural isomorphism
\[ Z_h\colon\ToSt(E_6)/[2]\to\Sth(E_6),\]
sending a total stability condition $\sigma$ to the far-end $12$-gon.
\end{theorem}
\begin{proof}
We can construct a slicing $\sli$ from a stable $12$-gon $\hgon$,
which is compatible with the central charge defined in Theorem~\ref{thm:geo.E6} and the shift $[1]$.
Such a construction is basically taking $\arg$ in appropriate length $2\pi$ intervals for central charges of each indecomposables.
It follows exactly the same line of work as in Section~\ref{sec:recon:Dn}.
Thus, we have the corresponding statement in type $E_6$.
\end{proof}

\section{$\ToSt$ of exceptional type $E_7$}\label{sec:E7}

\subsection{The $h$-gon of type $E_7$}\label{sec:E7-1}
\begin{figure}[hb]\centering
\includegraphics[width=4in]{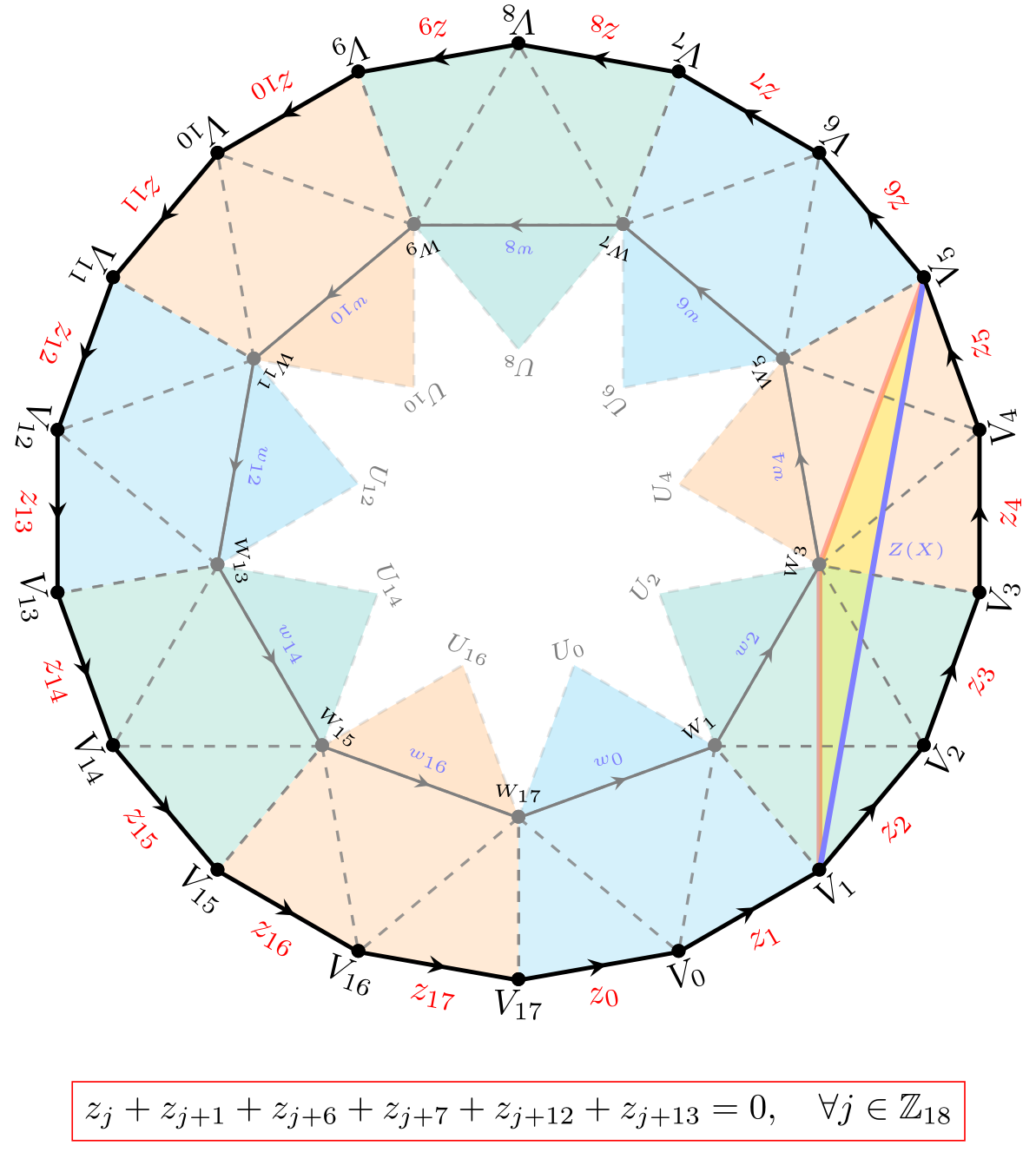}\vskip .5pc
\makebox[\textwidth][c]{
\begin{tikzpicture}[xscale=.64,yscale=.64,rotate=-45, font=\tiny,
arrow/.style={->,>=stealth,thick},font=\tiny]
\clip[rotate=-45](-.8,2)rectangle(7.8,27.5);
\foreach \j in {-3,...,12}{
\foreach \k in {0,...,12}{
    \draw[teal,->,>=stealth,very thick]
        (2*\k,2*\j+.5)to++(0,1);
    \draw[teal,->,>=stealth,very thick]
        (2*\k+.5,2*\j)to++(1,0);
}}
\draw[rotate=-45,white,fill=white]
    (0,2)rectangle(-2,28)(10,2)rectangle(7,28);
\foreach \j in {-3,...,12}{
\foreach \k in {0,...,12}{
    \draw[white](2*\k,2*\j) node[circle,draw=teal] (x\k\j) {$00$};
}}
\foreach \j in {-1,...,9}{
    \draw[](2*\j+6,2*\j-4) node[cyan!15,circle,draw=teal,fill=cyan!15] {$00$} node[font=\scriptsize]{$B_\j$};
  \draw[blue!0](2*\j+2+.005,2*\j-2-.005)
    node[circle,draw=teal,fill=blue!0](x\j) {$00$};
  \draw[teal,->,>=stealth,very thick]
    (2*\j+4+1+.005+.35,2*\j+1-.005+.35)to (2*\j+4+2-.35,2*\j+2-.35);
  \draw[teal,->,>=stealth,very thick] (2*\j+4+.35,2*\j+.35) to
    (2*\j+4+1+.005-.35,2*\j+1-.005-.35);
}\draw[](2*4+2+.005,2*4-2-.005)node[font=\small]{$X$};
\foreach \j/\k in {0,...,8}{
  \draw[red!0](2*\j+2+3+.005,2*\j+1-.005)
    node[circle,draw=teal,fill=yellow!15](m\j) {$00$};
  \draw[](m\j)node{$M_\j$};
}
\foreach \j/\k in {0,...,8}{
    \draw[](2*\j+2,2*\j+2) node[green!15,circle,draw=teal,fill=green!15] {$00$} node {$C_{\j}$};
}
\end{tikzpicture}}
\caption{The $18$-gon of type $E_7$ with its fire core and AR-quiver of $\Dwq(E_7)$}
\label{fig:E7}
\end{figure}

\begin{definition}\label{def:hgon E7}
An \emph{$h$-gon $\hgon$ of type $E_7$} is a symmetric $18$-gon satisfying \eqref{eq:E7-rel},
i.e. the equation in Figure~\eqref{fig:E7}.
\end{definition}

A direct calculation shows that:
\begin{lemma}\label{lem:E7}
After setting $z_{j+9}=-z_j$ for $j\in\ZZ_{18}$,
the set of 3 equations \eqref{eq:E7-rel} has rank 2.
So the space of $18$-gons of type $E_7$ has complex dimension $9-2=7$.
\end{lemma}

\begin{construction}\label{con:E7}
Using the hexagon relations in \eqref{eq:E7-rel},
we can draw hexagons
$$\mathbf{L}_{2j}\colon=V_{2j-1}V_{2j}V_{2j+1} W_{2j+1} U_{2j} W_{2j-1}$$ with edges
\begin{gather*}
    V_{2j-1}V_{2j}=z_{2j},\quad V_{2j}V_{2j+1}=z_{2j+1},\quad
    V_{2j+1} W_{2j+1}=z_{2j+6},\\ W_{2j+1} U_{2j}=z_{2j+7},\quad
    U_{2j} W_{2j-1}=z_{2j+12},\quad W_{2j-1}V_{2j-1}=z_{2j+13}
\end{gather*}
as shown in the upper picture of Figure~\ref{fig:E7}.
Note that we have 3 sets of 3 parallel hexagons
\[
    \{\mathbf{L}_{2j+6k}\mid k=0,1,2\}\quad\text{for $j\in\ZZ_3$},
\]
drawn in blue/green/orange respectively.

Set $w_j$ as in \eqref{eq:w_j} for $n=7$ and $h=18$.
As $z_{j+9}=-z_j$, we also have $w_{j+9}=-w_j$.
Moreover, each hexagon decomposes into 3 triangles corresponding to the above equation and 1 triangle corresponding to the relation
\begin{gather}\label{eq:tri7}
    w_{j}+w_{j+6}+w_{j+12}=0.
\end{gather}
So there is a 9-gon $\fcore$, called the \emph{fire core} of $\hgon$, with vertices $\{W_{2j+1}\mid j\in\ZZ_9\}$ and edges $w_{2j}=W_{2j-1}W_{2j+1}$.
\end{construction}

\begin{remark}
We take nine $\mathbf{L}_{2j}$ with even indices to obtain $\fcore$.
One can also take the other nine $\mathbf{L}_{2j+1}$ with odd indices to obtain another 9-gon $\icore$, called the \emph{ice core} of $\hgon$, with vertices $\{W_{2j}\mid j\in\ZZ_9\}$ and edges $w_{2j+1}=W_{2j}W_{2j+2}$.
However, since $\hgon$ is symmetric,
these two cores are also symmetric to each other.
\end{remark}

\subsection{Geometric model for the root category of type $E_7$}
\begin{setup}\label{setup:E7}
Take the orientation in \eqref{eq:En} with $n=7$.
Note that 1/3/7 are the mid-end/near-end/far-end vertices in this case.
The lower picture in Figure~\ref{fig:E7} is part of $\AR\Dwq(E_7)$ and
for simplicity, we only label the objects in the boundary $\tau$-orbits, i.e.
$C_j/M_j/B_j$ in \eqref{eq:137}.
\end{setup}

Given any central charge $Z\colon K\Dwq(E_7)\to\CC$.
Let $\hgon_Z$ be the far-end $18$-gon of $\Dwq(E_7)$ with edges $z_j=Z(B_j)$ for $1\le j\le 18$.
\begin{lemma}\label{lem:E7-1}
The set $\{[B_{j}] \mid 1\leq j\leq 9\}$ spans $K \Dwq(E_7)$.
Moreover, $\hgon_Z$ is a $18$-gon of type $E_7$ and
its ice/fire core is formed by (the odd/even) half of
$$\{ Z(C_{j}) \mid 1\leq j\leq 18\}$$
for $Z(C_{j})=w_{j}$.
\end{lemma}
\begin{proof}
The first statement can be checked directly, e.g. via dimension vectors.

Now we will show that all $z_{j}$ satisfy \eqref{eq:E7-rel}.
By \eqref{eq:common1} for $n=7$ we have
\[[C_{j+2}]=[B_j]+[B_{j+5}]\]
and thus $Z(C_{j+2})=z_j+z_{j+5}=w_{j+2}$ as we have the triangle relation \eqref{eq:w_j}.

Moreover, triangles in the second row/column of \eqref{eq:common2} for $n=7$ implies that
\[
    [B_{j}]+[C_{j+3}]=[M_{j+1}]
    =[C_{j}]+[B_{j+4}]
\]
and thus $z_{j}+w_{j+3}=w_{j}+z_{j+4}$.
Substitute \eqref{eq:w_j} to kill $w_{j}$,
we obtain
\[
    z_{j}+(z_{j+1}+z_{j+6})=(z_{j-2}+z_{j+3})+z_{j+4}.
\]
Noticing that $z_{k+9}=-z_{k}$, the above equation becomes \eqref{eq:E7-rel}.

Thanks to \eqref{eq:w_j}, the last statement is clear.
\end{proof}

\begin{remark}
The triangle relation \eqref{eq:tri7} corresponds to the triangle in the third row of \eqref{eq:common2}
as $L_j=C_{j+6}$ in this case.
\end{remark}

Now we describe a geometric model in $E_7$ case.

\begin{theorem}\label{thm:geo.E7}
An $18$-gon $\hgon$ of type $E_7$ is a geometric model for the root category $\Dwq(E_7)/[2]$
in the sense that by setting \eqref{eq:commonx} for $n=7$,
we obtain a central charge $Z\colon K\Dwq(E_7)\to\CC$.
\end{theorem}
As in type $E_6$,
for each $M_j$ in the $3^{\mathrm{rd}}$-$\tau$-orbits, it can be realized by \eqref{eq:4ways} for $h=18$.
For instance, the (edges of the) yellow triangle in the top picture of Figure~\ref{fig:E7}
corresponds to the (central charges of) triangle
\begin{equation}\label{eq:00}
M_2\to X\to M_3\to M_2[1]
\end{equation}
in $\Dwq(E_7)$, similar to the yellow triangle in the case of $E_6$.
\subsection{Stability of $18$-gon for $E_7$}\label{sec:sth E7}
\begin{definition}\label{def:sth E7}
An $h$-gon $\hgon$ of type $E_7$ is \emph{stable} if it is positively convex and
its ice/fire core is inside the level-4 diagonal-gon.
\end{definition}
Note that the positive convexity of $\hgon$ will be inherited by its ice and fire cores
as in the $E_6$ case.
Denote by $\Sth(E_7)$ the moduli space of stable $18$-gon of type $E_7$ up to translation.
By Lemma~\ref{lem:E7}, the complex dimension of $\Sth(E_7)$ is $7$.

\begin{proposition}\label{pp:sth E7}
If $\sigma=(Z,\sli)\in \ToSt(E_7)$,
then its far-end $18$-gon is a stable $18$-gon of type $E_7$.
\end{proposition}
\begin{proof}
To prove the statement, one needs to use the AR triangle \eqref{eq:00} similar as in type $E_6$, whose central charges of its terms form the yellow triangle in Figure~\ref{fig:E7}.
Then the rest argument follows the same way as in Proposition~\ref{pp:sth E6}.
\end{proof}

\begin{theorem}\label{thm:E7}
There is a natural isomorphism
\[ Z_h\colon\ToSt(E_7)/[2]\to\Sth(E_7).\]
sending a total stability condition $\sigma$ to the far-end $18$-gon.
\end{theorem}

\section{$\ToSt$ of exceptional type $E_8$}\label{sec:E8}
\begin{figure}[h]\centering
\makebox[\textwidth][c]{
\begin{tikzpicture}[xscale=.6,yscale=.6,rotate=-45, font=\tiny,
arrow/.style={->,>=stealth,thick},font=\tiny]
\clip[rotate=-45](-1,2)rectangle(9,24.5);
\foreach \j in {-3,...,20}{
\foreach \k in {0,...,20}{
    \draw[teal,->,>=stealth,thin]
        (2*\k,2*\j+.5)to++(0,1);
    \draw[teal,->,>=stealth,thin]
        (2*\k+.5,2*\j)to++(1,0);
}}
\draw[rotate=-45,white,fill=white]
    (0,2)rectangle(-2,46)(10,2)rectangle(8.5,46);
\foreach \j in {-3,...,18}{
\foreach \k in {0,...,18}{
    \draw[white](2*\k,2*\j) node[circle,draw=teal] (x\k\j) {$_x$};
}}
\foreach \j in {-2,...,15}{
  \draw[yellow!15](2*\j+4+1+.005,2*\j+1-.005)
    node[circle,draw=teal,fill=yellow!15](m\j) {$_x$};
  \draw[](m\j)node{$M_{\j}$};
  \draw[blue!0](2*\j+4+.005,2*\j-.005)
    node[circle,draw=blue!50,fill=blue!0](x\j) {$_x$};
  \draw[teal,->,>=stealth,very thin]
    (2*\j+4+3+.005+.35,2*\j+3-.005+.35)to (2*\j+4+4-.35,2*\j+4-.35);
  \draw[teal,->,>=stealth,very thin] (2*\j+6+.35,2*\j+2+.35) to
    (2*\j+4+3+.005-.35,2*\j+3-.005-.35);
}
\draw[](2*4+2+.005,2*4-2-.005)node[font=\small]{$X$};
\foreach \j/\k in {1,...,15}{
    \draw[](2*\j+6,2*\j-6) node[cyan!15,circle,draw=teal,fill=cyan!15] {$_x$} node[font=\tiny]{$B_{\j}$};
}
\foreach \j/\k in {0,...,14}{
    \draw[](2*\j+2,2*\j+2) node[cyan!15,circle,draw=teal,fill=green!10] {$_x$} node {$C_{\j}$};
}
\draw[white,fill=white,rotate=-45](-9,48)rectangle(9.3,44.5)
    (-1,2)rectangle(-.8,24.5);
\end{tikzpicture}\hskip 2cm}

\makebox[\textwidth][c]{\hskip 2cm
\begin{tikzpicture}[xscale=.6,yscale=.6,rotate=-45, font=\tiny,
arrow/.style={->,>=stealth,thick},font=\tiny]
\clip[rotate=-45](-.8,23.5)rectangle(9.2,44.5);
\foreach \j in {-3,...,20}{
\foreach \k in {0,...,20}{
    \draw[teal,->,>=stealth,thin]
        (2*\k,2*\j+.5)to++(0,1);
    \draw[teal,->,>=stealth,thin]
        (2*\k+.5,2*\j)to++(1,0);
}}
\draw[rotate=-45,white,fill=white]
    (0,2)rectangle(-2,46)(10,2)rectangle(8.5,46);
\foreach \j in {-3,...,18}{
\foreach \k in {0,...,18}{
    \draw[white](2*\k,2*\j) node[circle,draw=teal] (x\k\j) {$_x$};
}}
\foreach \j in {-2,...,15}{
  \draw[yellow!15](2*\j+4+1+.005,2*\j+1-.005)
    node[circle,draw=teal,fill=yellow!15](m\j) {$_x$};
  \draw[](m\j)node{$M_{\j}$};
  \draw[blue!0](2*\j+4+.005,2*\j-.005)
    node[circle,draw=blue!50,fill=blue!0](x\j) {$_x$};
  \draw[teal,->,>=stealth,very thin]
    (2*\j+4+3+.005+.35,2*\j+3-.005+.35)to (2*\j+4+4-.35,2*\j+4-.35);
  \draw[teal,->,>=stealth,very thin] (2*\j+6+.35,2*\j+2+.35) to
    (2*\j+4+3+.005-.35,2*\j+3-.005-.35);
}
\draw[](2*4+2+.005,2*4-2-.005)node[font=\small]{$X$};
\foreach \j/\k in {1,...,15}{
    \draw[](2*\j+6,2*\j-6) node[cyan!15,circle,draw=teal,fill=cyan!15] {$_x$} node[font=\tiny]{$B_{\j}$};
}
\foreach \j/\k in {0,...,14}{
    \draw[](2*\j+2,2*\j+2) node[cyan!15,circle,draw=teal,fill=green!10] {$_x$} node {$C_{\j}$};
}
\end{tikzpicture}}
\caption{The AR-quiver of $\Dwq(E_8)$}
\label{fig:AR E8}
\end{figure}


\subsection{The $h$-gon of type $E_8$}\label{sec:E8-1}
\begin{figure}[h]\centering
\includegraphics[width=5in]{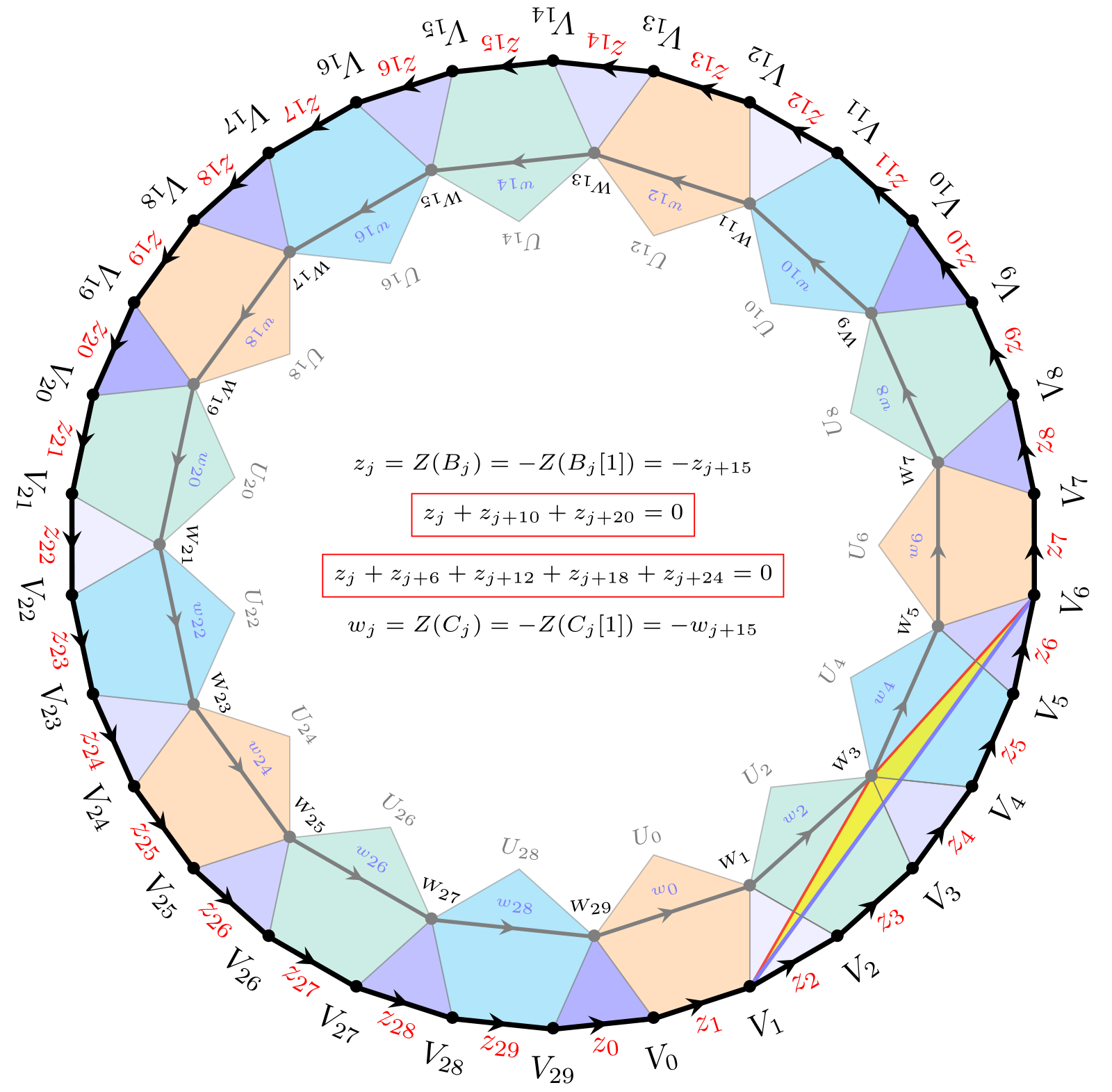}
\caption{The $30$-gon of type $E_8$ with its fire core}
\label{fig:E8-sthgon}
\end{figure}

\begin{definition}\label{def:hgon E8}
An \emph{$h$-gon $\hgon$ of type $E_8$} is a symmetric $30$-gon satisfying \eqref{eq:E8-rel},
i.e. the equation in Figure~\eqref{fig:E8-sthgon}.
\end{definition}

A direct calculation shows that:
\begin{lemma}\label{lem:E8}
After setting $z_{j+15}=-z_j$ for $j\in\ZZ_{30}$,
the set of 8 equations \eqref{eq:E8-rel} has rank 7.
So the space of $30$-gons of type $E_8$ has complex dimension $15-7=8$.
\end{lemma}

\begin{construction}\label{con:E8}
Using the triangle/pentagon relations in \eqref{eq:E8-rel},
we can draw pentagons
\[\mathbf{P}_{2j}\colon=V_{2j}V_{2j+1}W_{2j+1} U_{2j} W_{2j-1}\quad \forall j\in\ZZ_{15}\]
with edges
\[
 \begin{cases}
      V_{2j}V_{2j+1}=z_{2j+1},\; V_{2j+1}W_{2j+1}=z_{2j+7},\;
      W_{2j+1}U_{2j}=z_{2j+13},\\  U_{2j}W_{2j-1}=z_{2j+19},\;
      W_{2j-1} V_{2j}=z_{2j+25}
    \end{cases}
\]
and triangles
\[ \TT_{2j}\colon=V_{2j-1}V_{2j}W_{2j-1}\quad\forall j\in\ZZ_{15}\]
with edges
\[V_{2j-1}V_{2j}=z_{2j},\; V_{2j}W_{2j-1}=z_{2j+10},\;
    W_{2j-1}V_{2j-1}=z_{2j+20},\]
as shown in Figure~\ref{fig:E8-sthgon}.
Note that we have 3 sets of 5 parallel pentagons
\[
    \{\mathbf{P}_{2j+6k}\mid k=0,1,2,3,4\},\quad \text{for $j\in\ZZ_3$}
\]
drawn in orange/green/blue respectively
and 5 sets of 3 parallel triangles
\[
    \{\TT_{2j+10k}\mid k=0,1,2\},\quad \text{for $j\in\ZZ_5$}
\]
drawn in violet with different opacity respectively.

Set $w_j$ as in \eqref{eq:w_j} for $n=8$ and $h=30$.
As $z_{j+15}=-z_j$, we also have $w_{j+15}=-w_j$.
Note that as $W_{2j+1}U_{2j}=z_{2j+13}$ and $U_{2j}W_{2j-1}=z_{2j+19}$,
we have
$$W_{2j-1}W_{2j+1}=-(z_{2j+13}+z_{2j+19})=-w_{2j+15}=w_{2j}.$$

So there is a 15-gon $\fcore$, called the \emph{fire core} of $\hgon$ with vertices $\{W_{2j+1}\mid j\in\ZZ_{15}\}$ and edges $w_{2j}=W_{2j-1}W_{2j+1}$.
\end{construction}

\begin{remark}
We take fifteen $\mathbf{P}_{2j}$ and fifteen $\TT_{2j}$ with even indices to obtain $\fcore$.
One can also take the other fifteen $\mathbf{P}_{2j+1}$ and fifteen $\TT_{2j+1}$ with odd indices to obtain another 15-gon $\icore$,
called the \emph{ice core} of $\hgon$, with vertices $\{W_{2j}\mid j\in\ZZ_{15}\}$ and edges $w_{2j+1}=W_{2j}W_{2j+2}$.
However, since $\hgon$ is symmetric, these two cores are also symmetric to each other.
\end{remark}

\subsection{Geometric model for the root category of type $E_8$}
\begin{setup}\label{setup:E8}
Take the orientation in \eqref{eq:En} with $n=8$.
Note that 1/3/8 are the mid-end/near-end/far-end vertices in this case.
The picture in Figure~\ref{fig:AR E8} is part of $\AR\Dwq(E_8)$ and
for simplicity, we only label the objects in the boundary $\tau$-orbits, i.e.
$C_j/M_j/B_j$ in \eqref{eq:137}.
\end{setup}


Given any central charge $Z\colon K\Dwq(E_8)\to\CC$.
Let $\hgon_Z$ be the far-end $30$-gon of $\Dwq(E_8)$ with edges $z_j=Z(B_j)$ for $1\le j\le 30$.
\begin{lemma}\label{lem:E7-2}
The set $\{[B_{j}] \mid 1\leq j\leq 15\}$ spans $K \Dwq(E_8)$.
Moreover, $\hgon_Z$ is a $30$-gon of type $E_8$ and
its ice/fire core is formed by (the odd/even) half of
$$\{Z(C_{j})\mid 1\leq j\leq 30\}$$
for $Z(C_{j})=w_{j}$.
\end{lemma}
\begin{proof}
The first statement can be checked directly, e.g. via dimension vectors.

Now we will show that all $z_{j}$ satisfy \eqref{eq:E8-rel}.
By \eqref{eq:common1} for $n=8$ we have
\[[C_{j+2}]=[B_j]+[B_{j+6}]\]
and thus $Z(C_{j+2})=z_j+z_{j+6}=w_{j+2}$ as we have the triangle relation \eqref{eq:w_j}.

Moreover, triangles in the third row/column of \eqref{eq:common2} for $n=8$ implies that
\begin{gather}\label{eq:2triangle.e8}
\begin{cases}
    z_j-z_{j+5}+z_{j+10}=0,\\
    w_{j}-w_{j+3}+z_{j+10}=0,
\end{cases}\quad\forall j\in\ZZ_{30}\end{gather}
as $L_j=B_{j+10}$ in this case.
Noticing $-z_{j+5}=z_{j+20}$, the first equation of \eqref{eq:2triangle.e8} becomes the triangle relation in \eqref{eq:E8-rel}.
Substitute \eqref{eq:w_j} to kill $w_{}$'s in the second equation of \eqref{eq:2triangle.e8},
we have
\[z_{j-2}+z_{j+4}-z_{j+1}-z_{j+7}+z_{j+10}=0, \quad\forall j\in\ZZ_{30}.\]
Noticing $z_k=-z_{k+15}$, the above equation becomes the pentagon relation in \eqref{eq:E8-rel}.
\end{proof}

Now we describe a geometric model in type $E_8$.

\begin{theorem}\label{thm:geo.E8}
A $30$-gon $\hgon$ of type $E_8$ is a geometric model for the root category $\Dwq(E_8)/[2]$
in the sense that by setting \eqref{eq:commonx} for $n=8$
we obtain a central charge $Z\colon K\Dwq(E_8)\to\CC$.
\end{theorem}
As in types $E_6$ and $E_7$,
for each $M_j$ in the $3^{\mathrm{rd}}$-$\tau$-orbits, it can be realized by \eqref{eq:4ways} for $h=30$.
For instance, the (edges of the) yellow triangle in the top picture of Figure~\ref{fig:E8-sthgon}
corresponds to the (central charges of) triangle $M_2\to X\to M_3\to M_2[1]$ in $\Dwq(E_8)$.
\subsection{Stability of $30$-gon for $E_8$}\label{sec:sth E8}

\begin{definition}\label{def:sth E8}
A $30$-gon $\hgon$ of type $E_8$ is \emph{stable} if it is positively convex and
its ice/fire core is inside the level-5 diagonal-gon.
\end{definition}
Denote by $\Sth(E_8)$ the moduli space of stable $30$-gon of type $E_8$ up to translation.
By Lemma~\ref{lem:E8}, the complex dimension of $\Sth(E_8)$ is $8$.

As above in type $E_7$, we have the following proposition and theorem for $E_8$.
\begin{proposition}\label{pp:sth E8}
If $\sigma=(Z,\sli)\in\ToSt(E_8)$,
then its far-end $30$-gon is a stable $30$-gon of type $E_8$.
\end{proposition}

\begin{theorem}\label{thm:E8}
There is a natural isomorphism
\[ Z_h\colon\ToSt(E_8)/[2]\to\Sth(E_8).\]
sending a total stability condition $\sigma$ to the far-end $30$-gon.
\end{theorem}

\begin{remark}\label{rem:s.s.}
An interesting numerical observation is that
the triangle/square relations in type $E_6$ correspond to $x^4$ and $y^3$,
and triangle/pentagon relations in type $E_6$ correspond to $x^5$ and $y^3$, respectively,
in the corresponding simple (surface) singularity
 \[
    f(x,y,z)=\begin{cases}
        x^{n+1}+y^2+z^2,& Q=A_n (n\ge1);\\
        x^{n-1}+x y^2+z^2,& Q=D_n, (n\ge4);\\
        x^4+y^3+z^2,& Q=E_6;\\
        x^3 y+y^3+z^2,& Q=E_7;\\
        x^5+y^3+z^2,& Q=E_8.
   \end{cases}
\]
These equations are used to defined $\Dwq(Q)$ via (graded) matrix factorizations in \cite{KST}.
The correspondence in type $E_7$ is not as nice as the other two.
\end{remark}


\section{Non simply-laced types via folding}\label{sec:NON}
The non-simply-laced case is obtained by folding the corresponding simply-laced case, see Section~\ref{sec:KOT-Q}.
The derived category/stability conditions for $\S$ are the $\iota$-stable part of the one for $Q$ (cf. \cite{CQ} for more details),
Thus, the stable $h_Q$-gons for $\S$ are just $\iota$-stable $h$-gons, i.e. satisfying extra symmetry.
More precisely, we have the following
\begin{description}
\item[$B_n$] It is obtained by folding the corresponding quiver of type $D_{n+1}$.
Thus the extra symmetry condition on the corresponding doubly punctured $2n$-gon is the two punctures coincide (and thus at the geometric center of the $2n$-gon).
\item[$C_n$] It is obtained by folding the corresponding quiver of type $A_{2n-1}$.
Thus the extra symmetry condition on the corresponding $2n$-gon is central symmetry.
\item[$F_4$] It is obtained by folding the corresponding quiver of type $E_6$.
Thus the extra symmetry condition on the corresponding $12$-gon is central symmetry.
\item[$G_2$] It is obtained by folding the corresponding quiver of type $D_4$.
Thus the extra condition is: concision of the punctures, central symmetry and \eqref{eq:G2-rel}.
\end{description}
So Theorem~\ref{thm:0} also holds in these cases.

\begin{figure}\centering
\includegraphics[width=\textwidth]{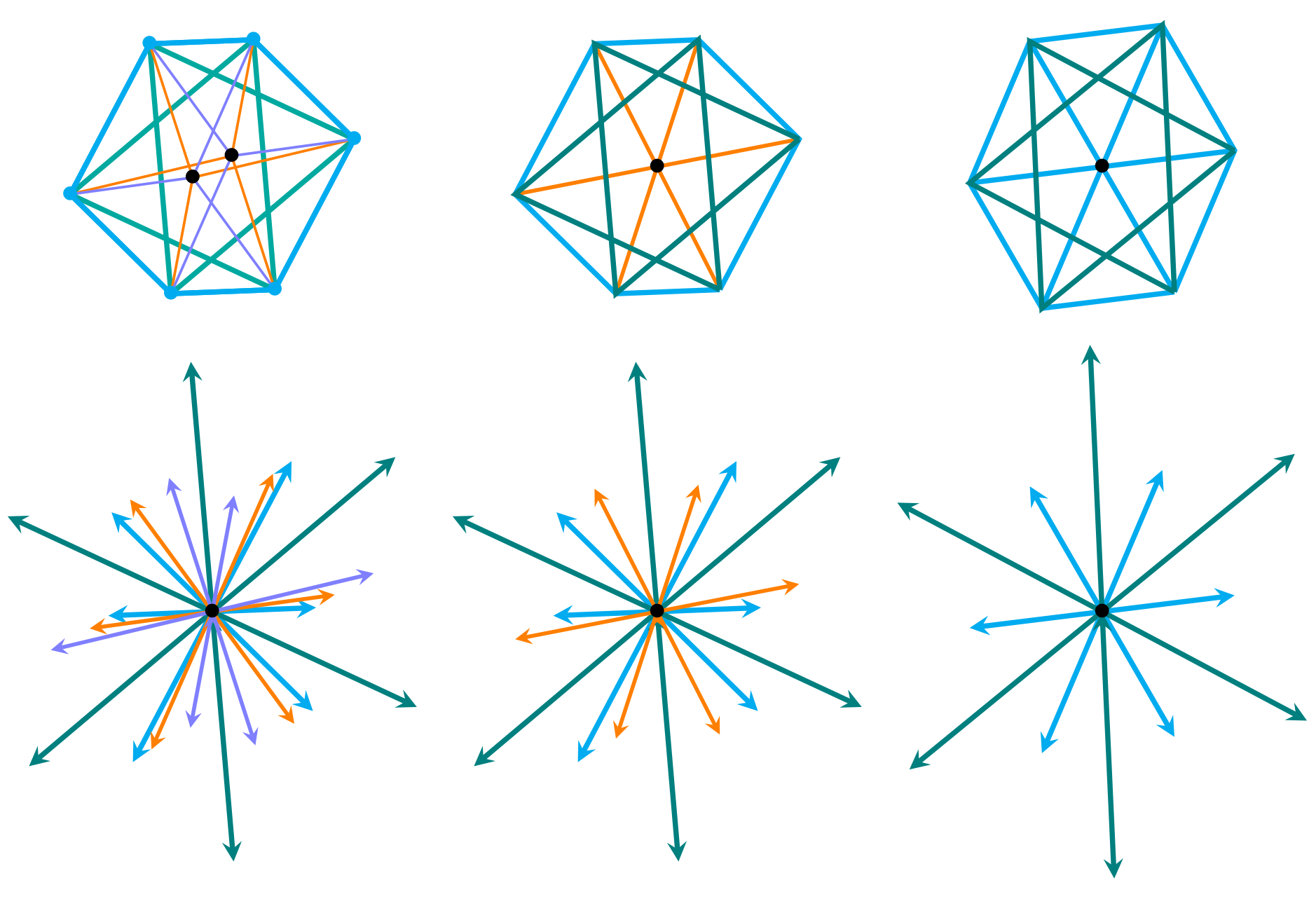}
\caption{Deforming a stable $h$-gon and central charge of type $D_4$ (left) to $B_3$ (middle) to $G_2$ (right)}\label{fig:DBG}
\end{figure}

\begin{example}
In Figure~\ref{fig:DBG}, we show that how to deform the central charge and associated far-end stable $h_Q$-gon of type $D_4$ into
the ones of type $B_3$ and of type $G_2$, respectively.

From $D_4$ to $B_3$, it requires the (orange/violet) central charges of objects in the two boundary (besides the chosen far-end) $\tau$-orbits coincides correspondingly.
Equivalently, the two punctures coincide as mentioned above.

From $D_4$ to $G_2$, it requires the (orange/violet/blue) central charges of objects in all boundary $\tau$-orbits coincides correspondingly.
Equivalently, the long diagonals of the $h_Q=6$-gon intersect at its geometric center that divide it into six parallel/anti-parallel triangles.
\end{example}

\addresseshere\newpage

\begin{appendix}
\section{Tikz Art Gallery by Qy}\label{app}
In this appendix, we collection figures of stable $h$-gons at Gepner points that interact with the corresponding projection of the root systems in the Coxeter plans, as well as a tilting of Ice and Fire of type $E_6$.

\begin{figure}[hb]\centering
\includegraphics[width=3in]{d5-intro.png}
\includegraphics[width=3in]{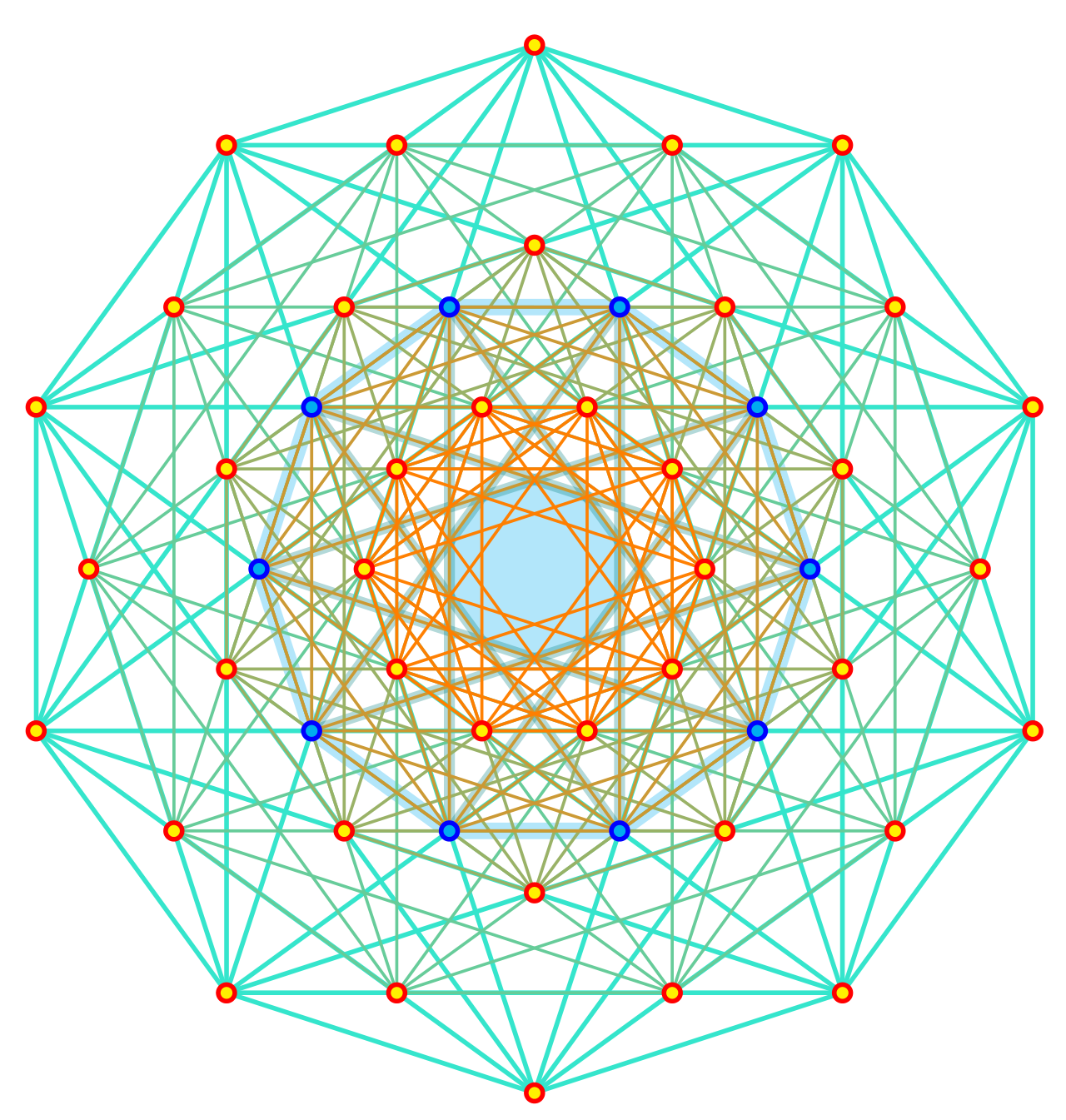}
\caption{The (Coxeter projection of the) root system of type $D_6$ and the stable $h$-gon}
\label{fig:D5}
\end{figure}

\begin{figure}[hb]\centering
\includegraphics[width=3.5in]{E6-1.png}
\includegraphics[width=3.5in]{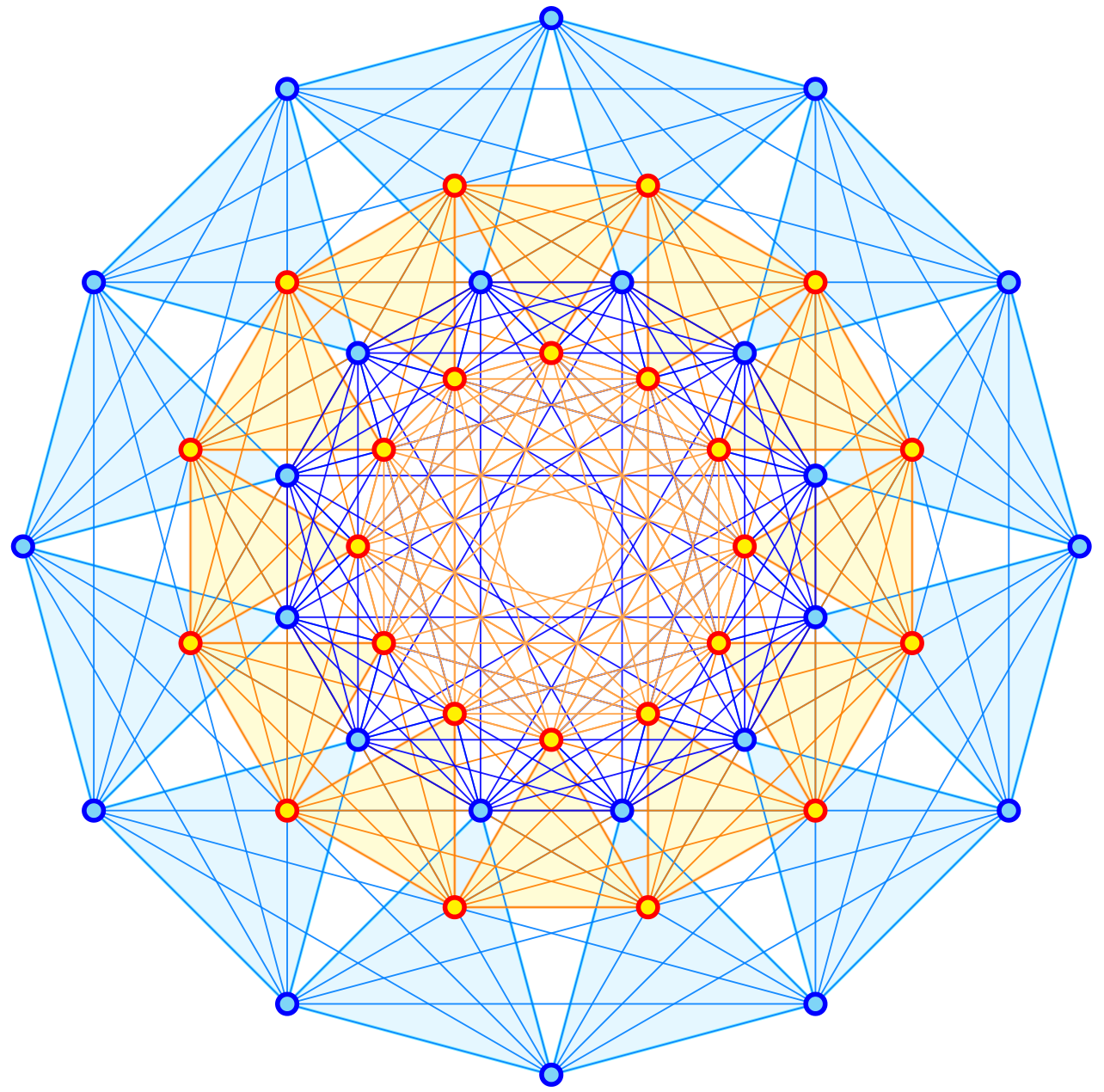}
\caption{The central charges, root system and stable $h$-gons of type $E_6$}
\label{fig:E6-RS}
\end{figure}

\begin{figure}[hb]\centering
\includegraphics[width=4in]{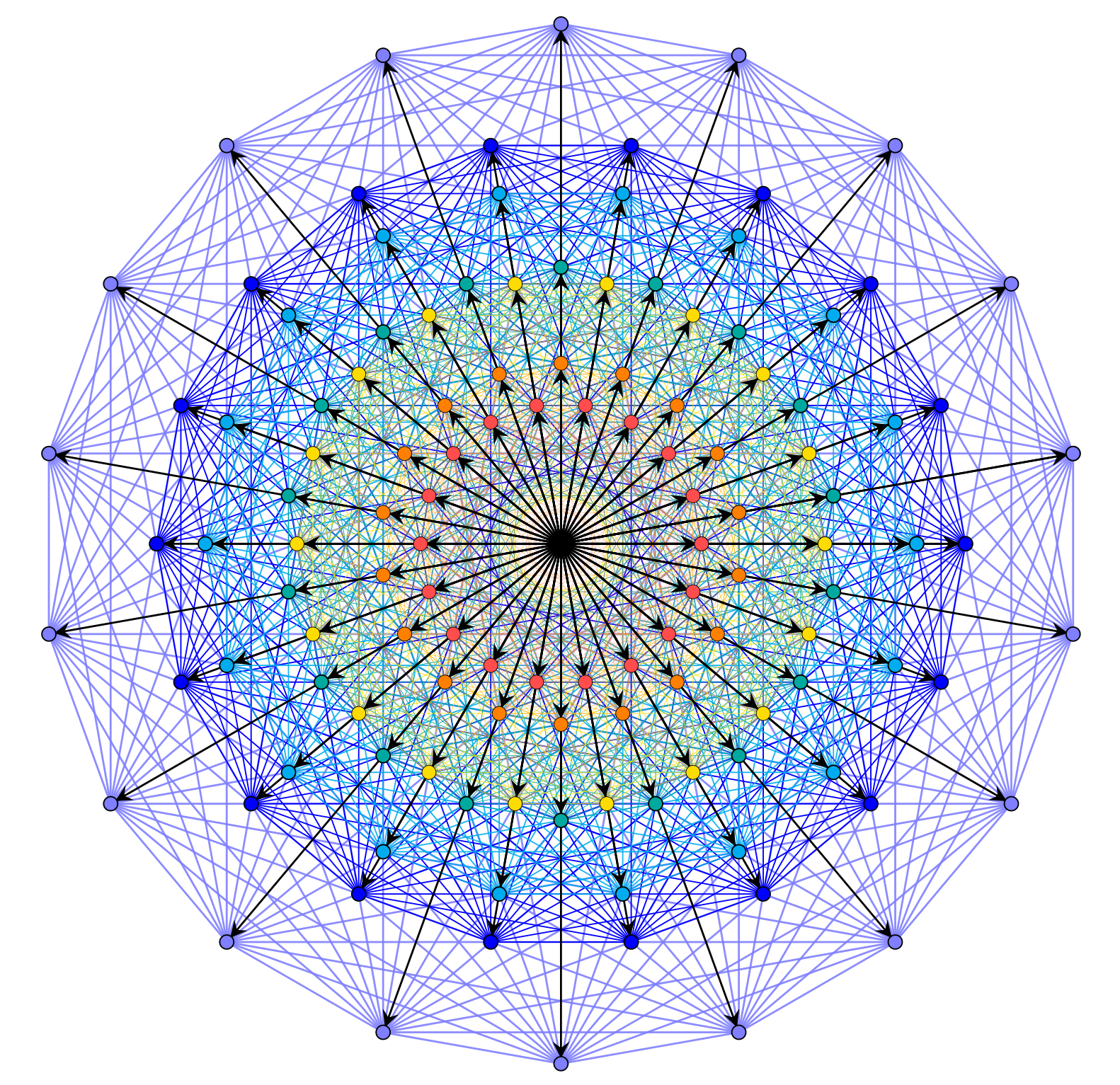}
\includegraphics[width=4in]{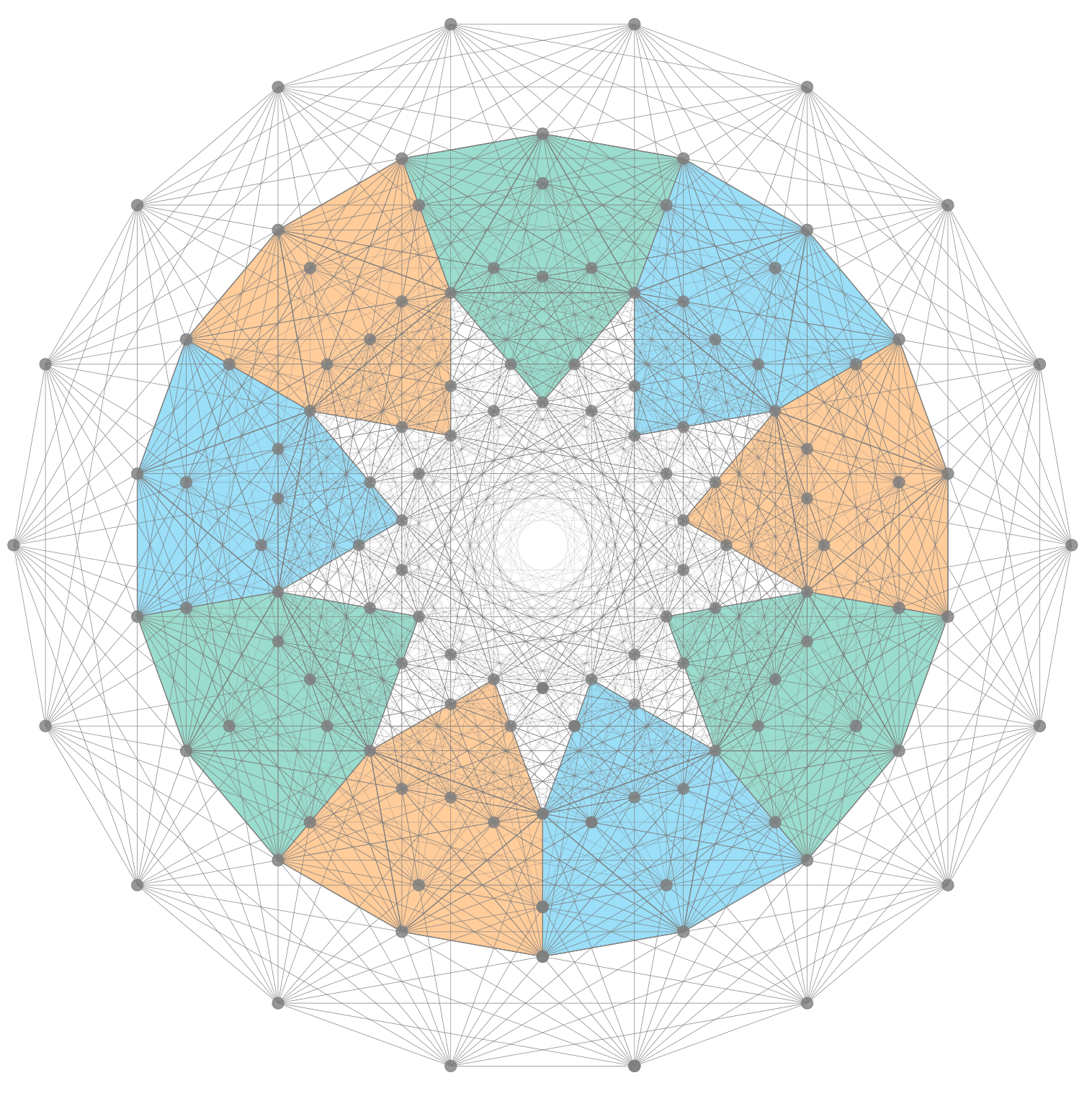}
\caption{The central charges, root system and stable $h$-gons of type $E_7$}
\label{fig:E7-1}
\end{figure}


\begin{figure}[hb]\centering
\makebox[\textwidth][c]{
\includegraphics[width=7in]{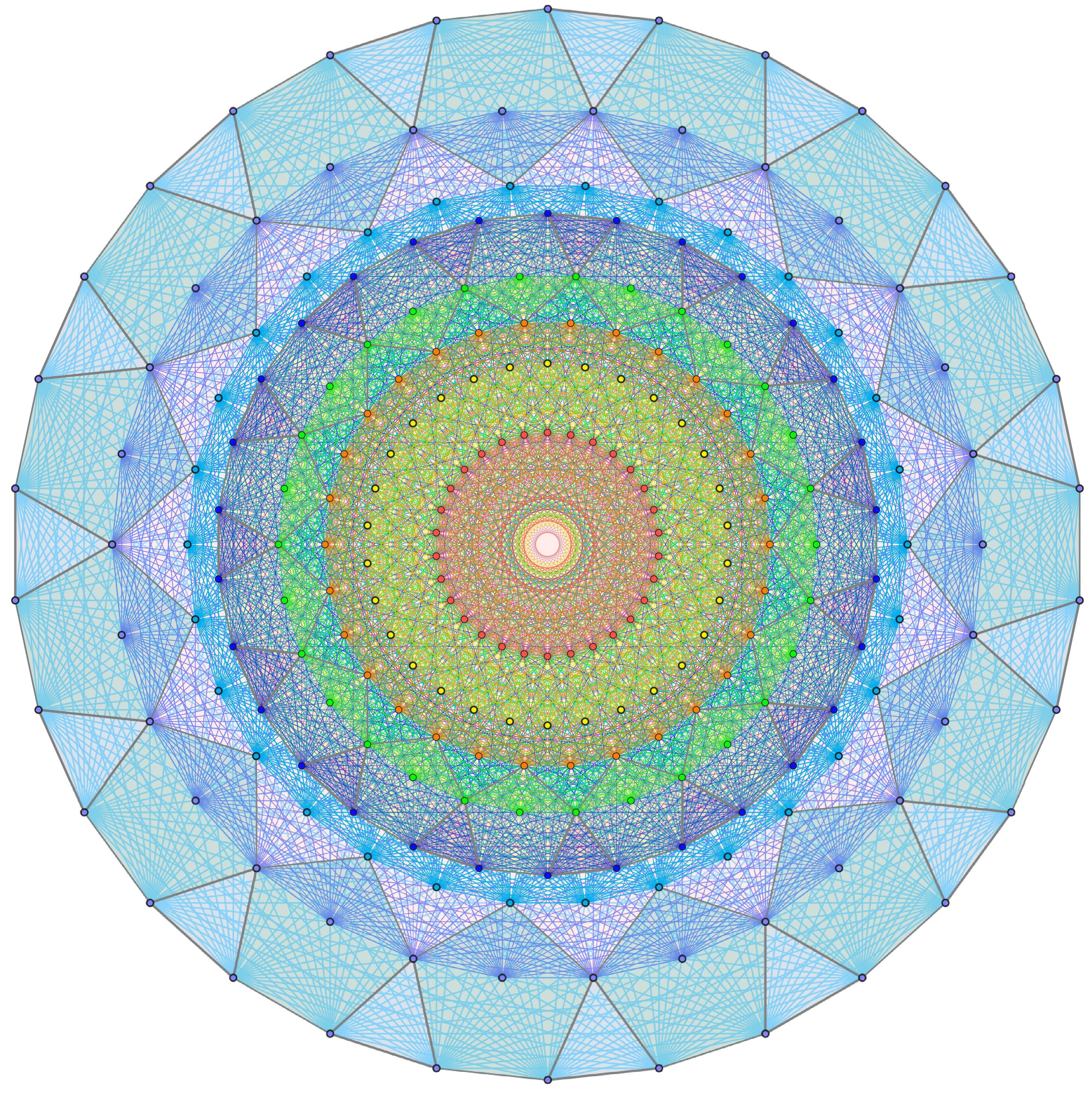}}
\caption{The root system and (two) stable $h$-gons of type $E_8$}
\label{fig:E8-RS}
\[\]
\end{figure}
\[\]
\[\]
\[\]
\end{appendix}

\begin{figure}[b]\centering
\vspace*{-4.5cm}
\makebox[\textwidth][c]{
\includegraphics[width=14in]{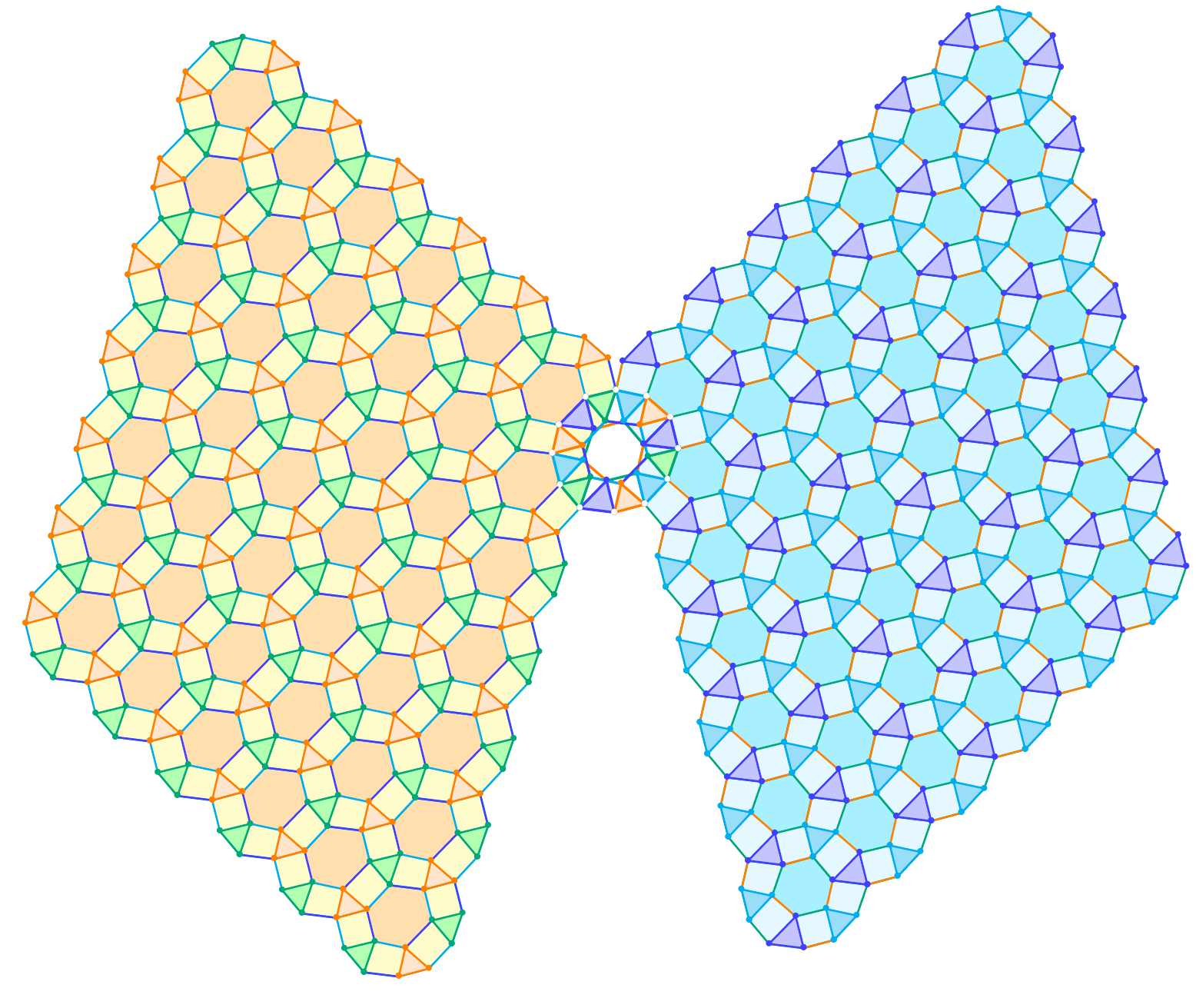}\quad}
\vskip -60pc
\caption{\textsf{A Tiling of Ice and Fire}}
\label{fig:A Tiling}
\end{figure}


\begin{thebibliography}{999999}
\newcommand{\au}[1]{\textrm{#1},}
\newcommand{\ti}[1]{\textrm{#1},}
\newcommand{\jo}[1]{\textit{#1}}
\newcommand{\vo}[1]{\textbf{#1}}
\newcommand{\yr}[1]{(#1)}
\newcommand{\pp}[2]{#1--#2.}
\newcommand{\arxiv}[1]{\href{http://arxiv.org/abs/#1}{arXiv:#1}}

\bibitem[AB]{AB}
\au{C.~Amiot, \and T.~Br\"{u}stle}
\ti{Derived equivalences between skew-gentle algebras}
\jo{Doc. Math.} \vo{27} \yr{2022} \pp{933}{982}
(\arxiv{1912.04367})

\bibitem[AP]{AP}
\au{C.~Amiot, \and P-G.~Plamondon}
\ti{The cluster category of a surface with punctures via group actions}
\jo{Adv. Math.} \vo{389} \yr{2021} Paper No. 107884, 63 pp.

\bibitem[BGMS]{BBMS}
  \au{E.~Barnard, E.~Gunawan, E.~Meehan \and R.~Schiffler}
  \ti{Cambrian combinatorics on quiver representations (type $A_n$)}
  \arxiv{1912.02840}

\bibitem[BM]{BM}
  \au{K.~Baur \and R.-Marsh}
  \ti{Categorification of a frieze pattern determinant}
  \arxiv{1008.5329}

\bibitem[Bri1]{B1}
  \au{T.~Bridgeland}
  \ti{Stability conditions on triangulated categories}
  \jo{Ann. Math.} \vo{166} \yr{2007} \pp{317}{345}
  (\href{http://arxiv.org/abs/math/0212237}{arXiv:0212237})

\bibitem[Bri2]{B2}
  \au{T.~Bridgeland}
  \ti{Stability conditions and Kleinian singularities}
  \jo{Int. Math. Res. Not.} \vo{21} \yr{2009} \pp{4142}{4157}
  (\href{http://arxiv.org/abs/math/0508257}{arXiv:0508257})

\bibitem[BQS]{BQS}
  \au{T.~Bridgeland, Y.~Qiu, and T.~Sutherland}
  \ti{Stability conditions and the ${A}_2$ quiver}
  \jo{Adv. Math.} \vo{365} \yr{2020} 107049.
  (\arxiv{1406.2566})

\bibitem[BS]{BS}
  \au{T.~Bridgeland \and I.~Smith}
  \ti{Quadratic differentials as stability conditions}
  \jo{Publ. Math. de l'IH\'{E}S}
  \vo{121} \yr{2015} \pp{155}{278}
  (\arxiv{1302.7030})

\bibitem[BQ]{BQ}
  \au{T.~Brustle \and Y.~Qiu}
  \ti{Tagged maing class group: Auslander-Reiten translations}
  \jo{Math. Zeit.} \vo{279} \yr{2015} \pp{1103}{1120}
  (\arxiv{1212.0007})

\bibitem[CQ]{CQ}
  \au{W. Chang \and Y.~Qiu}
  \ti{Folding quivers and numerical stability conditions}
  \jo{Publ. Res. Inst. Math. Sci.} \vo{60} \yr{2024} \pp{271}{303}
  (\arxiv{1210.0243})

\bibitem[CQZ]{CQZ}
  \au{W. Chang, Y.~Qiu and X.~Zhang}
  \ti{Geometric model for module categories of Dynkin quivers via hearts of total stability conditions}
  \jo{J. Algebra} \vo{638} \yr{2024} \pp{57}{89}
  (\arxiv{2208.00073})

\bibitem[DGK]{DGK}
  \au{Y.~Diaz, C.~Gilbert and R.~Kinser}
  \ti{Total stability and Auslander-Reiten theory for Dynkin quivers}
  \arxiv{2208.02445}.

\bibitem[FLLQ]{FLLQ}
  \au{Y. Fan, C. Li, W. Liu \and Y. Qiu}
  \ti{Contractibility of space of stability conditions on the projective plane via global dimension function}
  \jo{Math. Res. Lett.} \vo{30} \yr{2023} \pp{51}{87}
  (\arxiv{2001.11984})

\bibitem[FST]{FST}
  \au{S.~Fomin, M.~Shapiro \and D.~Thurston}
  \ti{Cluster algebras and triangulated surfaces, part I: Cluster complexes}
  \jo{Acta Math.} \vo{201} \yr{2008} \pp{83}{146}

\bibitem[Gab]{Gab}
  \au{A. Gabrièlov}
  \title{Dynkin diagrams of unimodal singularities. (Russian)}
  \jo{Funkcional. Anal. i Priložen.} \vo{8} \yr{1974} \pp{1}{6}

\bibitem[HKK]{HKK}
   \au{F.~Haiden, L.~Katzarkov and M.~Kontsevich}
   \ti{Stability in Fukaya categories of surfaces}
   \jo{Publ. Math. de l'IH\'{E}S}
   \vo{126} \yr{2017} \pp{247}{318}
   (\arxiv{1409.8611})

\bibitem[Hil]{H}
  \au{L.~Hille},
  \ti{Tilting modules over the path algebra of type A, polytopes, and Catalan numbers}. (Lie algebras and related topics)
  \jo{Contemp. Math.} \vo{652} \pp{91}{101} \yr{2015}
  (\arxiv{1505.06011})

\bibitem[IQ1]{IQ1}
  \au{A.~Ikeda \and Y.~Qiu}
  \ti{$q$-Stability conditions on Calabi-Yau-$\XX$ categories and twisted periods}
  \jo{Compos. Math.} \vo{159} \yr{2023} \pp{1347}{1386}
  )\arxiv{1807.00469})

\bibitem[IQ2]{IQ2}
  \au{A.~Ikeda \and Y.~Qiu}
  \ti{$q$-Stability conditions via $q$-quadratic differentials}
  \jo{Memoirs of Amer. Math. Soc.} to appear.
  (\arxiv{1812.00010})

\bibitem[KST]{KST}
  \au{H.~Kajiura, K.~Saito \and A.~Takahashi}
  \ti{Matrix factorizations and representations of quivers II:Type ADE case}
  \jo{Adv. Math. } \vo{211} \yr{2007} \pp{327}{362}
  (\href{https://arxiv.org/abs/math/0511155}{arxiv:0511155})

\bibitem[Kel]{K11}
  B.~Keller,
  On cluster theory and quantum dilogarithm.
  \href{http://arxiv.org/abs/1102.4148}{arXiv:1102.4148v4}.

\bibitem[KOT]{KOT}
  \au{K.~Kikuta, G.~Ouchi \and A.~Takahashi}
  \ti{Serre dimension and stability conditions}
  \arxiv{1907.10981}.

\bibitem[Kin]{Ki}
  \au{R.~Kinser}
  \ti{Total stability functions for type A quivers}
  \arxiv{2002.12396}

\bibitem[Len]{L}
  \au{H.~Lenzing}
  \ti{Coxeter Transformations associated with Finite Dimensional Algebras}
  \jo{Computational methods for representations of groups and algebras} \yr{1999} \pp{287}{308}

\bibitem[Qiu1]{Q11}
  \au{Y.~Qiu}
  \ti{Stability conditions and quantum dilogarithm identities for Dynkin quivers}
  \jo{Adv. Math.} \vo{269} \yr{2015} \pp{220}{264}
  (\arxiv{1111.1010})

\bibitem[Qiu2]{Q18}
  \au{Y.~Qiu}
  \ti{Global dimension function on stability conditions and Gepner equations}
  \jo{Math. Zeit.} \vo{303} \yr{2023} No. 11.
  (\arxiv{1807.00010})

\bibitem[Qiu3]{Q20}
  \au{Y.~Qiu}
  \ti{Contractible flow of stability conditions via global dimension function}
  \jo{J. Diff. Geom.} to appear.
  (\arxiv{2008.00282})

\bibitem[QW]{QW}
  \au{Y.~Qiu \and J.~Woolf}
  \ti{Contractible stability spaces and faithful braid group actions}
  \jo{Geom. Topol.} \vo{22} \yr{2018} \pp{3701}{3760}
  (\arxiv{1407.5986})

\bibitem[QZ]{QZ}
  \au{Y.~Qiu \and Y.~Zhou}
  \ti{Cluster categories for marked surfaces: punctured case}
  \jo{Compos. Math.} \vo{153} \yr{2017} \pp{1779}{1819}
  (\arxiv{1311.0010})

\bibitem[QZZ]{QZZ}
  \au{Y.~Qiu, C. Zhang \and Y. Zhou}
  \ti{Two geometric models for graded skew-gentle algebras}
in preparation.

\bibitem[Rei]{R}
  M.~Reineke, The Harder-Narasimhan system in quantum groups and cohomology of quiver moduli, \emph{Invent. Math.} 152 (2003), no. 2, 349-368.
  (\href{http://arxiv.org/abs/math/0204059}{arXiv:math/0204059v1})

\bibitem[Sch]{S}
  \au{R.~Schiffler}
  \ti{A geometric model for cluster categories of type $D_n$}
  \jo{J. Algebraic Combin.} \vo{27} \yr{2008} \pp{1}{21}

\end{thebibliography}
\end{document}